\documentclass{article}
\usepackage{amsfonts}
\usepackage{amsmath}
\usepackage{graphicx}
\usepackage{multicol}
\newtheorem{theorem}{Theorem}

\newtheorem{corollary}[theorem]{Corollary}

\newtheorem{lemma}[theorem]{Lemma}

\newtheorem{proposition}[theorem]{Proposition}
\newtheorem{remark}[theorem]{Remark}

\newenvironment{proof}[1][Proof]{\noindent \textbf{#1.} }{\  \rule{0.5em}{0.5em}}
\input{tcilatex}
\begin{document}

\title{\textbf{Robust parameter estimation of the log-logistic distribution
based on density power divergence estimators}}
\date{}
\date{}
\author{Felipe, A.; Jaenada, M., Miranda, P. and Pardo, L. \\
$^{1}${\small Department of Statistics and O.R., Complutense University of
Madrid, Spain}}

\maketitle

\begin{abstract}
	Robust inferential methods based on divergences measures have shown an appealing trade-off between efficiency and robustness in many different statistical models.
	In this paper, minimum density power divergence estimators (MDPDEs) for the scale and shape parameters of the log-logistic distribution are considered.
	The log-logistic is a versatile distribution modeling lifetime data which is  commonly adopted in survival analysis and reliability engineering studies when the hazard rate is initially increasing but then it decreases after some point.
	Further, it is shown that the classical estimators based on maximum likelihood (MLE) are included as a particular case of the MDPDE family. Moreover, the corresponding influence function of the MDPDE is obtained, and its boundlessness is proved, thus leading to robust estimators. A simulation study is carried out to illustrate  the slight loss in efficiency of MDPDE with respect to MLE and, at besides, the considerable gain in robustness.
\end{abstract}

\noindent \underline{\textbf{AMS 2001 Subject Classification}}\textbf{: }%
62F35, 62J12\medskip

\noindent \underline{\textbf{Keywords and phrases}}\textbf{: }Log-logistic
distribution, Minimum density power divergence estimator, robustness.

\section{Introduction}

In this paper we deal with the log-logistic distribution. This distribution arises as the logarithmic transformation of the logistic distribution, similar to how the log-normal distribution is derived from a normal distribution. The corresponding probability density function is given by
\begin{equation}
f_{\alpha ,\beta }(x)=\frac{\beta \alpha ^{\beta }x^{\beta -1}}{\left(
x^{\beta }+\alpha ^{\beta }\right) ^{2}}, \quad x>0,\, (\alpha
,\beta >0),  \label{9.1}
\end{equation}%
where $\alpha >0$ is the scale parameter representing the median of the distribution and $\beta >0$ is the shape parameter that controls the form of the distribution.
Moreover, the distribution is uni-modal for all $\beta >1$.
The log-logistic distribution shares the same form as the log-normal distribution but exhibits heavier tails. It is a special case of the Burr type-II distribution family \cite{bur42} as well as  the Kappa distribution family \cite{mijo73}. Furthermore, for the log-logistic distribution, the $k$-th moment exists only when $k < \beta$, where $\beta$ is the distribution parameter, and is given by
\begin{equation*}
E\left[ X^{k}\right] =\alpha ^{k}B\left( 1-\frac{k}{\beta },1-\frac{k}{\beta
}\right) ,
\end{equation*}%
where by $B\left( a,b\right) $ we are denoting the function beta of
parameters $a$ and $b,$ $i.e,$%
\begin{equation}\label{beta}
B\left( a,b\right) =\int_{0}^{1}x^{a-1}\left( 1-x\right) ^{b-1}dx.
\end{equation}

The log-logistic distribution finds diverse applications across various fields. Some appealing examples of application include its use in Engineering \cite{asma03}; Survival analysis, where it serves as a parametric model for events with hazard rates that initially increase and then decrease—such as in the study of mortality from cancer following diagnosis \cite{guaklv99}; Hydrology, for example, \cite{shmitr88} the log-logistic distribution has been applied to model stream flow rates and precipitation in Canada; Economics, particularly in the study of wealth or income distributions, where the log-logistic distribution is known as the Fisk distribution \cite{fis61} and is considered an equivalent alternative to a lognormal distribution. Additionally, in Networking, the log-logistic distribution is employed to model the time a set of data is processed sequentially by several computers and then returns to the initial computer.

Let us now discuss the estimation of the log-logistic model parameters, $\alpha$ and $\beta.$ These two parameters are typically
estimated via classical methods such as maximum likelihood and least squares estimation approaches. From a frequentist perspective, \cite{bama87} obtained the best linear unbiased estimation (BLUE) of the two model parameters of the log-logistic distribution.
Later, \cite{kasr02} developed the maximum likelihood estimator (MLE) of the scale parameter $\alpha$ for known shape parameter $\beta.$ In \cite{redowa18}, unbiased or nearly unbiased estimators for the parameters of the log-logistic distribution are proposed.
Additional inferential methods for parameter estimation of the log-logistic distribution can be found in \cite{shmitr88, guaklv99, asma03}, among others. In the meantime, Bayesian inference of the log-logistic parameters has also received attention in the literature. Recently, \cite{abta16} considered objective Bayesian estimation of this distribution using reference prior and Jeffreys' prior.
In \cite{hechqi20}, MLE based on simple random sampling and ranked set sampling techniques are reconsidered. Moreover, in \cite{hechya21}, the modified best linear unbiased estimator for estimating $\beta$ is studied.

The aforementioned existing estimators perform well for complete set of data without any contamination; however, they may become unreliable and result in severely distorted estimates when the data under study contain outlier observations. In other words, the previous estimators are asymptotically efficient (indeed, they are BAN (Best Asymptotically Normal Estimators)), but non-robust. Despite the lack of robustness of classical estimators, robust estimation of the log-logistic parameter has not been widely considered in the statistical literature yet. To the best of our knowledge, the only reference we have for robust estimation is \cite{mawapa23}.

To address this gap in the literature, this paper introduces a new family of estimators for the parameters of the log-logistic distribution based on Density Power Divergences (DPD), first presented in \cite{bahahjjo98}. The selection of this family is justified by its demonstrated efficiency and robustness in numerous statistical models. Indeed, we prove below that the  influence function of the MDPDE  is bounded, thus demonstrating their robustness against individual data points. Furthermore, we derive the asymptotic distribution of the estimators, showing that it generalizes the corresponding distribution of the classical estimators obtained via maximum likelihood.

The rest of the paper goes as follows.
In Section 2, we introduce key concepts related to the DPD, required for the subsequent discussion. In Section 3 we present the estimators of the parameters of a log-logistic distribution obtained through the DPD approach, hereafter referred to as MDPDEs. The asymptotic distribution of MDPDEs is also derived. Section 4 provides the expression of the influence function for these estimators. In Section 5, we carry out a simulation study to empirically evaluate the performance of MDPDEs in the presence of outliers, as well as to compare our proposed estimators with others proposals in the literature. Section 6 discusses a real dataset regarding the breakdown point of an insulating fluid. The paper concludes with a revision of findings and open problems for future research. The proofs of the results established in the paper are included in an Appendix.

\section{The minimum density power divergence estimator}

Let $\mathcal{G}$ denote the set of all distributions having densities with respect to a dominating measure (generally the Lebesgue measure or the
counting measure). Given any two densities $g$ and $f$ in $\mathcal{G}$, the
{\bf density power divergence} (DPD) between them is defined, as a function depending on a
non-negative tuning parameter $\tau,$ by
\begin{equation}
d_{\tau }(g,f)=\int \left \{ f^{1+\tau }(x)-\left( 1+\frac{1}{\tau }\right)
f^{\tau }(x)g(x)+\frac{1}{\tau }g^{1+\tau }(x)\right \} dx.  \label{10.1}
\end{equation}%

The previous definition can be extended to the limiting value $\tau =0$  by taking the continuous limit as $\tau \rightarrow 0$. Then, at $\tau=0,$ the DPD is defined by
\begin{equation*}
d_{0}(g,f)=\lim_{\tau \rightarrow 0}d_{\tau }(g,f)=\int g(x)\log f(x)dx,
\end{equation*}%
i.e., it coincides with the well-known Kullback-Leibler divergence. For $\tau =1$ the square of
the standard $L_{2}$ distance between $g$ and $f$ is obtained. The quantity
defined in equation (\ref{10.1}) is a genuine divergence in the sense that $d_{\tau }(g,f)\geq 0$ for all $g,f\in \mathcal{G}$ with $\tau \geq 0$, and $d_{\tau }(g,f)$ is equal to zero if and only if the densities $g$ and $f$
are identically equal. Intuitively, the DPD quantifies the statistical closeness between two distributions, so the more similar two distributions are, the smaller the DPD value between them. More details about power divergence measures can be found in \cite{bahahjjo98}.

Let us consider a parametric family of densities depending on a parameter  $\boldsymbol{\theta} \in \Theta,$ denoted by $\mathcal{F}_{\boldsymbol{\theta}} = \{f_{\boldsymbol{\theta}}, \boldsymbol{\theta} \in \Theta \}.$
Additionally, let $G$ denote another distribution function corresponding to the density $g,$ that is modeled through the parametric family $\mathcal{F}_{\boldsymbol{\theta}}$.
Because the DPD measures the similarity between two distributions, the best parameter value $\boldsymbol{\theta}$ approximating the distribution $G$ should minimize the divergence between the true and assumed densities, $g$ and $f_{\boldsymbol{\theta}}.$ Therefore, the minimum density power divergence functional at the distribution $G$, denoted by $\boldsymbol{T}_{\tau}(G)$, is defined as
\begin{equation}
d_{\tau }(g,f_{\boldsymbol{T}_{\tau }(G)})=\min_{\boldsymbol{\theta }\in \boldsymbol{\Theta }} d_{\tau}(g,f_{\boldsymbol{\theta }}).  \label{10.2}
\end{equation}%
where $\boldsymbol{T}_{\tau }(G) \in \Theta$ denotes the best parameter so the distributions $G$ and $F_{\boldsymbol{T}_{\tau }(G)}$ are as close as possible. Indeed, if the distribution $G$ belongs to the parametric family $\mathcal{F}_{\boldsymbol{\theta}}$ with parameter $\boldsymbol{\theta}_0,$ then the minimum DPD functional coincides with the true parameter underlying, i.e. $\boldsymbol{T}_{\tau }(G)= \boldsymbol{\theta}_0.$

In practice, distribution $G$ underlying the data is unknown and we only have a sample $(X_1, ..., X_n)$ of size $n$ of values from the population of interest with
distribution function $G$ (and density $g$). Hence, $G$ need to be estimated by the empirical distribution of the data, that we will denote $G_n$.
Therefore the {\bf minimum density power divergence estimator} (MDPDE) of $\boldsymbol{\theta }$ can defined by
\begin{equation}
\widehat{\boldsymbol{\theta }}_{\tau }=\boldsymbol{T}_{\tau }(G_{n}).
\label{10.3}
\end{equation}%

Remark that minimize $d_{\tau }(g_n, f_{\boldsymbol{\theta }})$ with respect to $\boldsymbol{\theta }$ is equivalent to maximize the following alternative function
\begin{equation*}
h_{n,\tau }(\boldsymbol{\theta )=}\left( 1+\frac{1}{\tau }\right) \frac{1}{n}%
\sum_{i=1}^{n}f_{\boldsymbol{\theta }}^{\tau }(X_{i})-\int f_{\boldsymbol{%
\theta }}^{1+\tau }(x)dx,
\end{equation*}%
and similarly, maximize $h_{n,\tau }(\boldsymbol{\theta )}$ with respect to $%
\boldsymbol{\theta }$ is equivalent to maximize
\begin{equation}
H_{n,\tau }(\boldsymbol{\theta )} = \left( 1+\frac{1}{\tau }\right) \frac{1}{n}%
\sum_{i=1}^{n}f_{\boldsymbol{\theta }}^{\tau }(X_{i})-\int f_{\boldsymbol{%
\theta }}^{1+\tau }(x)dx+\frac{1}{\tau },  \label{10.5}
\end{equation}%
as $\frac{1}{\tau }$ does not depend on $\boldsymbol{\theta .}$ We will adopt this last expression as the DPD-objective function.

Note that at $\tau=0$ it can be easily shown that
\begin{equation*}
H_{n,0}(\boldsymbol{\theta })=\lim_{\tau \rightarrow 0}H_{n,\tau }(%
\boldsymbol{\theta })=\frac{1}{n}\sum_{i=1}^{n}\log f_{\boldsymbol{\theta }%
}(X_{i}).
\end{equation*}
thus recovering the expression of the  log-likelihood function, as discussed before. That it, the MLE is included on the MDPDE family at the limiting case with $\tau=0.$

Summarizing the above, the MDPDE of $\boldsymbol{\theta}$ with tuning parameter $\tau \geq 0$, $\widehat{\boldsymbol{\theta }}_{\tau}, $ is given by
\begin{equation*}
\begin{aligned}
\widehat{\boldsymbol{\theta }}_{\tau} &= \arg \max_{\boldsymbol{\theta }\in
	\Theta }H_{n,\tau }(\boldsymbol{\theta }) \\
	&=
	\begin{cases}
	{\displaystyle \arg \max_{\boldsymbol{\theta }\in \Theta } 	\left( 1+\frac{1}{\tau }\right) \frac{1}{n}%
	\sum_{i=1}^{n}f_{\boldsymbol{\theta }}^{\tau }(X_{i})-\int f_{\boldsymbol{%
			\theta }}^{1+\tau }(x)dx+\frac{1}{\tau } }& \tau > 0\\
		\arg \max_{\boldsymbol{\theta }\in \Theta }
	{\displaystyle \frac{1}{n}\sum_{i=1}^{n}\log f_{\boldsymbol{\theta }}(X_{i})} & \tau = 0
\end{cases}.
\end{aligned}
 \label{10.6}
\end{equation*}

Building upon the preceding expressions, we now proceed to derive explicit expressions for the log-logistic model.
To compute the objective function $H_{n,\tau }(\boldsymbol{\theta}) = H_{n,\tau }(\beta, \alpha)$ for the log-logistic distribution, we will require the following technical lemma and a proposition.
\begin{lemma}
\label{lema1} Let $m(\beta )$ be a real function for which there exists the first
derivative. Then,
\begin{equation*}
I_{1}=\int_{0}^{\infty }\frac{t^{m\left( \beta \right) }}{\left( 1+t\right)
^{s}}dt=B\left( s-m\left( \beta \right) -1,m\left( \beta \right) +1\right) ,
\end{equation*}%
where by $B\left( a,b\right) $ we are denoting the function beta of
parameters $a$ and $b,$ as defined in (\ref{beta}).
\end{lemma}

\begin{proof}
See Appendix.
\end{proof}

\begin{proposition}
\label{proposition1}For the log-logistic distribution we get,
\begin{equation*}
\int_{0}^{\infty }f_{\alpha ,\beta }^{1+\tau }(x)dx\boldsymbol{=}\left(
\frac{\beta }{\alpha }\right) ^{\tau }B\left( \frac{\beta \tau +\tau +\beta
}{\beta },\frac{\beta \tau -\tau +\beta }{\beta }\right) .
\end{equation*}
\end{proposition}

\begin{proof}
See Appendix.
\end{proof}

Now, taking into account that
\begin{equation*}
f_{\alpha ,\beta }^{\tau }(x)=\frac{\beta ^{\tau }\alpha ^{\tau \beta
}x^{\tau \left( \beta -1\right) }}{\left( x^{\beta }+\alpha ^{\beta }\right)
^{2\tau }},
\end{equation*}%
we have that, for $\tau >0,$ the DPD-based objective function is given by
\begin{eqnarray*}
H_{n,\tau }(\beta ,\alpha \boldsymbol{)} &\boldsymbol{=}&\left( 1+\frac{1}{%
\tau }\right) \frac{1}{n}\sum_{i=1}^{n}f_{\boldsymbol{\theta }}^{\tau
}(X_{i})-\int_{0}^{\infty }f_{\boldsymbol{\theta }}^{1+\tau }(x)dx +\frac{1}{%
\tau } \\
&=&\frac{1}{n}\dsum\limits_{i=1}^{n}\left\{ \left( 1+\frac{1}{\tau }\right)
\frac{\beta ^{\tau }\alpha ^{\tau \beta }X_{i}^{\tau \left( \beta -1\right) }%
}{\left( X_{i}^{\beta }+\alpha ^{\beta }\right) ^{2\tau }} \right\} -\left( \frac{%
\beta }{\alpha }\right) ^{\tau }B\left( \frac{\beta \tau +\tau +\beta }{%
\beta },\frac{\beta \tau -\tau +\beta }{\beta }\right) +\frac{1}{\tau }%
\end{eqnarray*}%
and for $\tau =0,$%
\begin{equation*}
H_{n,0}(\beta ,\alpha \boldsymbol{)=}\frac{1}{n}\sum_{i=1}^{n}\left\{ \log
\beta +\beta \log \frac{X_{i}}{\alpha }-2\log \left( 1+\left( \frac{X_{i}}{%
\alpha }\right) ^{\beta }\right) +d\right\} ,
\end{equation*}%
being $d$ a constant that does not depend on $\beta $ and $\alpha .$
The MDPDE of tuning parameter $\tau$ for $\alpha $ and $\beta $ will denoted in the rest of the paper by $\widehat{\alpha }_{\tau }$ and $\widehat{\beta }_{\tau },$ respectively.

Because the MDPDEs maximize a differentiable function, they must annuls its first derivatives. That is, the MDPDEs can be alternatively computed as a solution of the system of equations
$$ {\partial H_{n,\tau }(\beta ,\alpha )\over \partial \alpha }= 0 \hspace{0.3cm} \text{and} \hspace{0.3cm} \quad  {\partial H_{n,\tau }(\beta ,\alpha )\over \partial \beta }= 0.$$

For any positive value of $\tau$, the derivative of the objective function $H_{n,\tau }(\beta ,\alpha)$ with respect to $\alpha $ is given by
\begin{equation*}
\dsum\limits_{i=1}^{n}\left\{ \left( 1+\frac{1}{\tau }\right) \beta ^{\tau
}X_{i}^{\tau \left( \beta -1\right) }\left( \frac{\tau \beta \alpha ^{-\tau
\beta -1}}{\left( X_{i}^{\beta }+\alpha ^{\beta }\right) ^{2\tau }} -\frac{%
2\tau \beta \alpha ^{\beta \left( \tau +1\right) }}{\alpha \left(
X_{i}^{\beta }+\alpha ^{\beta }\right) ^{2\tau +1}}\right) \right.
\end{equation*}

\begin{equation*}
\left. -\frac{\tau }{%
\alpha }\left( \frac{\beta }{\alpha }\right)^{\tau } B\left( \frac{\beta \tau +\tau
+\beta }{\beta },\frac{\beta \tau -\tau +\beta }{\beta }\right) \right\} =0
\end{equation*}%
and for $\tau =0,$%
\begin{equation*}
\dsum\limits_{i=1}^{n}\left\{ -\frac{\beta }{\alpha }+\frac{2\beta }{\alpha }%
\left( \frac{X_{i}^{\beta }}{\alpha ^{\beta }}\right) \left( \frac{%
X_{i}^{\beta }}{\alpha ^{\beta }}+1\right) ^{-1}\right\} =0.
\end{equation*}%

On the other hand, to obtain the derivative with respect to $\beta$ for  $\tau >0$ we may rewrite $H_{n,\tau }(\beta ,\alpha)$ as follows:
\begin{equation*}
H_{n,\tau }(\beta ,\alpha \boldsymbol{)=}\frac{1}{n}\tsum\limits_{i=1}^{n}%
\left( 1+\frac{1}{\tau }\right) f(X_{i},\beta )-\alpha ^{-\tau }h(\beta ) + {1\over \tau },
\end{equation*}%
being
\begin{equation*}
f(X_{i},\beta )=\frac{\beta ^{\tau }\alpha ^{\tau \beta }X_{i}^{\tau \left(
\beta -1\right) }}{\left( x_{i}^{\beta }+\alpha ^{\beta }\right) ^{2\tau }}%
=\alpha ^{-\tau }\frac{\beta ^{\tau }\left( \frac{X_{i}}{\alpha }\right)
^{\tau \left( \beta -1\right) }}{\left( \frac{X_{i}^{\beta }}{\alpha ^{\beta
}}+1\right) ^{2\tau }}
\end{equation*}%
and
\begin{equation*}
h(\beta )=\beta ^{\tau }B\left( \frac{\beta \tau +\tau +\beta }{\beta },%
\frac{\beta \tau -\tau +\beta }{\beta }\right) .
\end{equation*}%

Now we have,
\begin{eqnarray*}
\frac{\partial f(X_{i},\beta )}{\partial \beta } &=&\alpha ^{-\tau }\tau
\beta ^{\tau }\left( \frac{X_{i}}{\alpha }\right) ^{\tau \left( \beta
-1\right) }\left( \left( \frac{1}{\beta }+\log \frac{X_{i}}{\alpha }\right)
\left( \frac{X_{i}^{\beta }}{\alpha ^{\beta }}+1\right) ^{-2\tau }\right. \\
&&\left. -2\frac{X_{i}^{\beta }}{\alpha ^{\beta }}\left( \log \frac{X_{i}}{%
\alpha }\right) \left( \frac{X_{i}^{\beta }}{\alpha ^{\beta }}+1\right)
^{-2\tau -1}\right)
\end{eqnarray*}%
and
\begin{equation*}
\frac{\partial h(\beta )}{\partial \beta }=\tau \beta ^{\tau -1}B\left( \frac{\beta \tau +\tau
+\beta }{\beta },\frac{\beta \tau -\tau +\beta }{\beta }\right) \left( 1+%
\frac{\tau }{\beta }\left( \Psi \left( \frac{\beta \tau -\tau +\beta }{\beta
}\right) -\Psi \left( \frac{\beta \tau +\tau +\beta }{\beta }\right) \right)
\right) .
\end{equation*}

Hence, the equation for $\tau >0$ with respect to $\beta $ is then given by
\begin{eqnarray*}
\frac{\partial H_{n,\tau }(\beta ,\alpha \boldsymbol{)}}{\partial \beta }
&=&\dsum\limits_{i=1}^{n}\left( 1+\frac{1}{\tau }\right) \left\{ \alpha
^{-\tau }\tau \beta ^{\tau }\left( \frac{X_{i}}{\alpha }\right) ^{\tau
\left( \beta -1\right) }\left( \left( \frac{1}{\beta }+\log \frac{X_{i}}{%
\alpha }\right) \left( \frac{X_{i}^{\beta }}{\alpha ^{\beta }}+1\right)
^{-2\tau }\right. \right. \\
&&\left. \left. -2\frac{X_{i}^{\beta }}{\alpha ^{\beta }}\left( \log \frac{%
X_{i}}{\alpha }\right) \left( \frac{X_{i}^{\beta }}{\alpha ^{\beta }}%
+1\right) ^{-2\tau -1}\right) -\alpha ^{-\tau }\tau \beta ^{\tau -1}B\left(
\frac{\beta \tau +\tau +\beta }{\beta },\frac{\beta \tau -\tau +\beta }{%
\beta }\right) \right. \\
&&\left. \times \left( 1+\frac{\tau }{\beta }\right) \left( \Psi \left( \frac{\beta
\tau -\tau +\beta }{\beta }\right) -\Psi \left( \frac{\beta \tau +\tau
+\beta }{\beta }\right) \right) \right\} =0,
\end{eqnarray*}%
and for $\tau =0,$%
\begin{equation*}
\dsum\limits_{i=1}^{n}\left\{ \frac{1}{\beta }-\log \alpha +\log X_{i}-2%
\frac{X_{i}^{\beta }}{\alpha ^{\beta }}\left(
\frac{X_{i}^{\beta }}{\alpha ^{\beta }}+1\right) ^{-1} \log \frac{X_{i}}{\alpha }\right\} =0.
\end{equation*}%

\section{Asymptotic distribution of MDPDE}

In this section we will compute the asymptotic distribution of the MDPDEs for both model parameters, marginally for each parameter (assuming fixed the other) and jointly. To facilitate the calculations in this section, we state the following two lemmas.
\begin{lemma}
\label{lema2}Let $m(\beta )$ be a real function for which there exists the first
derivative. Then,
\begin{eqnarray*}
I_{2} &=&\int_{0}^{\infty }\left( \log t\right) \frac{t^{m\left( \beta
\right) }}{\left( t+1\right) ^{s}}dt \\
&=&B\left( s-m(\beta )-1,m(\beta )+1\right) \left \{ \Psi \left( m(\beta
)+1\right) -\Psi \left( s-m(\beta )-1\right) \right \} ,
\end{eqnarray*}%
where and $\Psi \left( x\right) $ is the digamma function defined as the
logarithmic derivative of the gamma function.
\end{lemma}
\begin{proof}
See Appendix.
\end{proof}

\begin{lemma}
\label{lema3}Let $m(\beta )$ be a real function for which there exists the first
derivative. Then,
\begin{eqnarray*}
I_{3} &=&\int_{0}^{\infty }\left( \log t\right) ^{2}\frac{t^{m\left( \beta
\right) }}{\left( t+1\right) ^{s}}dt \\
&=&B(s-m\left( \beta \right) -1,m\left( \beta \right) +1)\left \{ \left(
\Psi \left( m\left( \beta \right) +1\right) -\Psi \left( s-m\left( \beta
\right) -1\right) \right) ^{2}\right. \\
&&\left. +\left( \Psi ^{\prime }\left( m\left( \beta \right) +1\right) +\Psi
^{\prime }\left( s-m\left( \beta \right) -1\right) \right) \right \} .
\end{eqnarray*}
\end{lemma}
\begin{proof}
See Appendix.
\end{proof}


\subsection{Asymptotic distribution for $\protect\widehat{\protect\alpha }_{%
\protect\tau }$}

In relation to the asymptotic distribution of $\widehat{\alpha }_{\tau }$,
assuming $\beta $ known, we have \cite{bahahjjo98},

\begin{equation*}
\sqrt{n}\left( \widehat{\alpha }_{\tau }-\alpha \right) \underset{%
n\longrightarrow \infty }{\overset{\mathcal{L}}{\longrightarrow }}\mathcal{N}%
\left( 0,J_{\tau }^{-1}(\alpha )K_{\tau }(\alpha )J_{\tau }^{-1}(\alpha
)\right) ,
\end{equation*}%
being%
\begin{equation}\label{J}
J_{\tau }(\alpha )=\int_{0}^{\infty }\left( \frac{\partial \log
f_{\alpha }(x)}{\partial \alpha }\right) ^{2}f_{\alpha }(x)^{\tau +1}dx\text{
and }K_{\tau }(\alpha )=J_{2\tau }(\alpha )-\xi _{\tau }\left( \alpha
\right) ^{2},
\end{equation}%
with

\begin{equation*}
\xi _{\tau }\left( \alpha \right) =\int_{0}^{\infty }\frac{\partial
\log f_{\alpha }(x)}{\partial \alpha }f_{\alpha }(x)^{\tau +1}dx.
\end{equation*}%

In the next theorem we are going to get the expression of $J_{\tau }(\alpha
).$

\begin{theorem}
\label{Theorem1}For the log-logistic distribution given in (\ref{9.1}) we
have%
\begin{equation}\label{J-alpha}
J_{\tau }(\alpha )=\left( \frac{\beta }{\alpha }\right) ^{\tau +2}B\left(
\frac{\beta \tau +\tau +\beta }{\beta },\frac{\tau \beta -\tau +\beta }{%
\beta }\right) \left[ 2\frac{\left( \beta \tau +\tau +\beta \right) \left(
-\tau \beta -\beta +\tau \right) }{\beta ^{2}\left( \tau +1\right) \left(
2\tau +3\right) }+1\right] .
\end{equation}
\end{theorem}

\begin{proof}
See Appendix.
\end{proof}

\begin{remark}
If we observe the expression in Eq. (\ref{J-alpha}), it follows that for $\tau =0$ we recover the Fisher information for $\alpha .$ Therefore, we have

\begin{equation*}
I_{F}\left( \alpha \right) =J_{\tau =0}(\alpha )=\left( \frac{\beta }{\alpha }\right)
^{2}(1+\frac{4}{3}-2)=\frac{\beta ^{2}}{3\alpha ^{2}}.
\end{equation*}
\end{remark}

Let us now turn to $K_{\tau }(\alpha ).$

\begin{theorem}
\label{Theorem1a}For the log-logistic distribution given in (\ref{9.1}) we
have,
\begin{equation*}
K_{\tau }(\alpha )=J_{2\tau }(\alpha )-\xi _{\tau }\left( \alpha \right)
^{2}.
\end{equation*}%
being
\begin{equation}\label{xi-alpha}
\xi _{\tau }\left( \alpha \right) =\xi _{\tau }\left( \alpha \right) =\left(
\frac{\beta }{\alpha }\right) ^{\tau +1}B\left( \frac{\tau \beta +\beta
+\tau }{\beta },\frac{\tau \beta +\beta -\tau }{\beta }\right) \left( \frac{%
-\tau }{\beta +\tau \beta }\right) .
\end{equation}
\end{theorem}

\begin{proof}
See Appendix.
\end{proof}

\begin{corollary}
We can observe that for Eq. (\ref{xi-alpha}) and $\tau =0$ we have
\begin{equation*}
\xi _{\tau =0}\left( \alpha \right) =0, \quad \text{and}\quad K_{\tau =0}(\alpha
)=I_{F}\left( \alpha \right) .
\end{equation*}
\end{corollary}

Based on the previous results, we obtain the classical results for the MLE, $\hat{\alpha }_{\tau=0} = \hat{\alpha },$ i.e.

$$ \sqrt{n}\left( \hat{\alpha } - \alpha \right) \rightarrow {\cal N}(0, I_F\left( \alpha \right)^{-1} ).$$

\subsection{Asymptotic distribution for $\protect\widehat{\protect\beta }_{%
\protect\tau }$}

In relation to the asymptotic distribution of $\widehat{\beta }_{\tau }$,
assuming $\alpha $ known, we have again
\begin{equation*}
\sqrt{n}\left( \widehat{\beta }_{\tau }-\beta \right) \underset{%
n\longrightarrow \infty }{\overset{\mathcal{L}}{\longrightarrow }}\mathcal{N}%
\left( 0,J_{\tau }^{-1}(\beta )K_{\tau }(\beta )J_{\tau }^{-1}(\beta )\right) ,
\end{equation*}%
being

\begin{equation*}
J_{\tau }(\beta )=\int_{0}^{\infty }\left( \frac{\partial \log
f_{\beta }(x)}{\partial \beta }\right) ^{2}f_{\beta }(x)^{\tau +1}dx\text{
and }K_{\tau }(\beta )=J_{2\tau }(\beta )-\xi _{\tau }\left( \beta \right) ,
\end{equation*}
with

\begin{equation*}
\xi _{\tau }\left( \beta \right) =\int_{0}^{\infty }\frac{\partial
\log f_{\beta }(x)}{\partial \beta }f_{\beta }(x)^{\tau +1}dx.
\end{equation*}

\begin{theorem}
\label{Theorem2}For the log-logistic distribution given in (\ref{9.1}) we
have%
\begin{equation*}
J_{\tau }(\beta )=N_{1}+N_{2}+N_{3}+N_{4}+N_{5}+N_{6},
\end{equation*}%
where $N_{1},N_{2},N_{3},N_{4},N_{5}$ and $N_{6}$ are given in (\ref{N1}), (%
\ref{N2}), (\ref{N3}),(\ref{N4}), (\ref{N5}) and (\ref{N6}), respectively.
\end{theorem}

\begin{proof}
See Appendix.
\end{proof}

\begin{corollary}
The Fisher information for $\beta $ is given by

\begin{equation*}
I_{F}\left( \beta \right) = J_{\tau =0}(\beta
)= \frac{3+\pi ^{2}}{9 \beta ^{2}}
\end{equation*}
\end{corollary}

\begin{theorem}
\label{Theorem3}For the log-logistic distribution given in (\ref{9.1}) we
have,
\begin{equation*}
K_{\tau }(\beta )=J_{2\tau }(\beta )-\xi _{\tau }\left( \beta \right) ^{2}.
\end{equation*}%
being

\begin{eqnarray}\label{xi-beta}
\xi _{\tau }(\beta ) &=& \frac{\beta ^{\tau -1}}{\alpha ^{\tau }}{\tau \over \tau +1} B\left( \frac{\tau \beta
+\beta +\tau }{\beta },\frac{\tau \beta +\beta -\tau }{\beta }\right) \\
& & \left\{ 1+\Psi \left( \frac{\tau \beta +\beta -\tau }{\beta }\right) -\Psi
\left( \frac{\tau \beta +\beta +\tau }{\beta }\right) \right\} \nonumber
\end{eqnarray}%

\end{theorem}

\begin{proof}
See Appendix.
\end{proof}

\begin{corollary}
We can observe that for Eq. (\ref{xi-beta}) and $\tau =0$ we have
\begin{equation*}
\xi _{\tau =0}\left( \beta \right) =0, \quad \text{and}\quad K_{\tau =0}(\beta
)=I_{F}\left( \beta \right) .
\end{equation*}
\end{corollary}

Based on the previous results, we obtain again that for $\hat{\beta }_{\tau=0} = \hat{\beta },$ it follows

$$ \sqrt{n}\left( \hat{\beta } - \beta \right) \rightarrow {\cal N}(0, I_F\left( \beta \right)^{-1} ).$$

\subsection{Asymptotic distribution for $\left( \protect\widehat{\protect\alpha
}_{\protect\tau },\protect\widehat{\protect\beta }_{\protect\tau }\right) $}

In relation to the asymptotic distribution of $\left( \widehat{\alpha }%
_{\tau },\widehat{\beta }_{\tau }\right) $ , we have

\begin{equation*}
\sqrt{n}\left( \left( \widehat{\alpha }_{\tau },\widehat{\beta }_{\tau
}\right) ^{T}-\left( \alpha ,\beta \right) ^{T}\right) \underset{%
n\longrightarrow \infty }{\overset{\mathcal{L}}{\longrightarrow }}\mathcal{N}%
\left( \boldsymbol{0},\boldsymbol{J}_{\tau }^{-1}(\alpha ,\beta )\boldsymbol{%
K}_{\tau }(\alpha ,\beta )\boldsymbol{J}_{\tau }^{-1}(\alpha ,\beta )\right) ,
\end{equation*}%
being%
\begin{equation*}
\boldsymbol{J}_{\tau }(\alpha ,\beta )=\left(
\begin{array}{cc}
J_{\tau }^{11}\left( \alpha ,\beta \right) & J_{\tau }^{12}\left( \alpha
,\beta \right) \\
J_{\tau }^{12}\left( \alpha ,\beta \right) & J_{\tau }^{22}\left( \alpha
,\beta \right)%
\end{array}%
\right)
\end{equation*}%
and
\begin{equation*}
\boldsymbol{K}_{\tau }(\alpha ,\beta )=\boldsymbol{J}_{2\tau }(\alpha ,\beta
)-\xi _{\tau }\left( \alpha ,\beta \right) ^{T}\xi _{\tau }\left( \alpha
,\beta \right) .
\end{equation*}%

Here,

\begin{equation*}
J_{\tau }^{11}\left( \alpha ,\beta \right) =J_{\tau }\left( \alpha \right) ,%
\text{ }J_{\tau }^{22}\left( \alpha ,\beta \right) =J_{\tau }\left( \beta
\right) ,\text{ }J_{\tau }^{12}\left( \alpha ,\beta \right)
=\int_{0}^{\infty }\left( \frac{\partial \log f_{\alpha ,\beta }(x)}{%
\partial \alpha }\right) \left( \frac{\partial \log f_{\alpha ,\beta }(x)}{%
\partial \beta }\right) f_{\alpha ,\beta }(x)^{\tau +1}dx
\end{equation*}%
and

\begin{equation*}
\xi _{\tau }\left( \alpha ,\beta \right) =\left( \xi _{\tau }\left( \alpha
\right) ,\xi _{\tau }\left( \beta \right) \right) .
\end{equation*}

Hence, it suffices to compute $J_{\tau }^{12}\left( \alpha ,\beta \right).$ This is done in next theorem.

\begin{theorem}
\label{Theorem4}For the log-logistic distribution given in (\ref{9.1}) we
have%
\begin{equation*}
J_{\tau }^{12}\left( \alpha ,\beta \right)
=B_{1}+B_{2}+B_{3}+B_{4}+B_{5}+B_{6},
\end{equation*}%
where $B_{1},B_{2},B_{3},B_{4},B_{5}$ and $B_{6}$ are given in (\ref{B1}), (%
\ref{B2}), (\ref{B3}), (\ref{B4}), (\ref{B5}) and (\ref{B6}), respectively.

\end{theorem}

\begin{proof}
See Appendix.
\end{proof}

\begin{corollary}
For $\tau =0$ it can be seen that $J_{\tau
=0}^{12}\left( \alpha ,\beta \right) =0.$ Therefore, the Fisher information
matrix is given by%
\begin{equation*}
\boldsymbol{I}_{F}\left( \alpha ,\beta \right) =\left(
\begin{array}{cc}
J_{\tau =0}\left( \alpha \right) & \text{ }J_{\tau =0}^{12}\left( \alpha
,\beta \right) \\
\text{ }J_{\tau =0}^{12}\left( \alpha ,\beta \right) & J_{\tau =0}\left(
\beta \right)%
\end{array}%
\right) =\left(
\begin{array}{cc}
\frac{\beta ^{2}}{3\alpha ^{2}} & \text{ }0 \\
\text{ }0 & \frac{3+\pi ^{2}}{\beta ^{2}9}%
\end{array}%
\right) .
\end{equation*}%
\end{corollary}

A directly obtention of the Fisher Information matrix can be seen in \cite{hechqi20}.

\section{Influence function}

The influence function (IF) is a classical tool for measuring (local) robustness which indicates the possible asymptotic bias in the estimation functional due to an infinitesimal contamination. Intuitively, the IF quantifies the impact that entails a infinitesimal point contamination in the sample on the parameter estimates. Then, an estimator is said to be locally robust if its associated IF is bounded. In this Section, we obtain the IF of the MDPDE for  $\alpha$ and $\beta,$ and prove its boundlessness.

As discussed in Section 2, the MDPDEs  are obtained by solving the derivative of the DPD loss with respect to both model parameters.
Obtaining the model parameters by solving a system of equations characterizes the MDPDE as an M-estimator. Thus, its influence function can be obtained by means of the theory of M-estimators, widely discussed in the robust inference literature.
%
%
After some algebra, the influence function for the MDPDE of $\alpha ,$ with $\tau >0,$
is given by
\begin{eqnarray*}
IF(x,\alpha ,\tau ) &=&\frac{1}{J\left( \alpha \right) }\left( \left( 1+%
\frac{1}{\tau }\right) \beta ^{\tau }x^{\tau \left( \beta -1\right) }\left(
\dfrac{\tau \beta \alpha ^{\tau \beta -1}}{\left( x^{\beta }+\alpha ^{\beta
}\right) ^{2\tau }}-\dfrac{2\tau \beta \alpha ^{\beta \left( \tau +1\right) }%
}{\alpha \left( x^{\beta }+\alpha ^{\beta }\right) ^{2\tau +1}}\right)
\right. \\
&&\left. -\frac{\tau }{\alpha }\left( \frac{\beta }{\alpha }\right) ^{\tau
}B\left( \frac{\beta \tau +\tau +\beta }{\beta },\frac{\beta \tau -\tau
+\beta }{\beta }\right) \right)
\end{eqnarray*}%
and for $\tau =0,$%
\begin{equation*}
IF(x,\alpha ,0)=\left( \frac{\beta ^{2}}{3\alpha ^{2}}\right) ^{-1}\left( -%
\frac{\beta }{\alpha }+\frac{2\beta }{\alpha }\left( \frac{x^{\beta }}{%
\alpha ^{\beta }}\right) \left( \frac{x^{\beta }}{\alpha ^{\beta }}+1\right)
^{-1}\right) .
\end{equation*}%

\begin{figure}
	\centering
	\includegraphics[height=6cm, width=12cm]{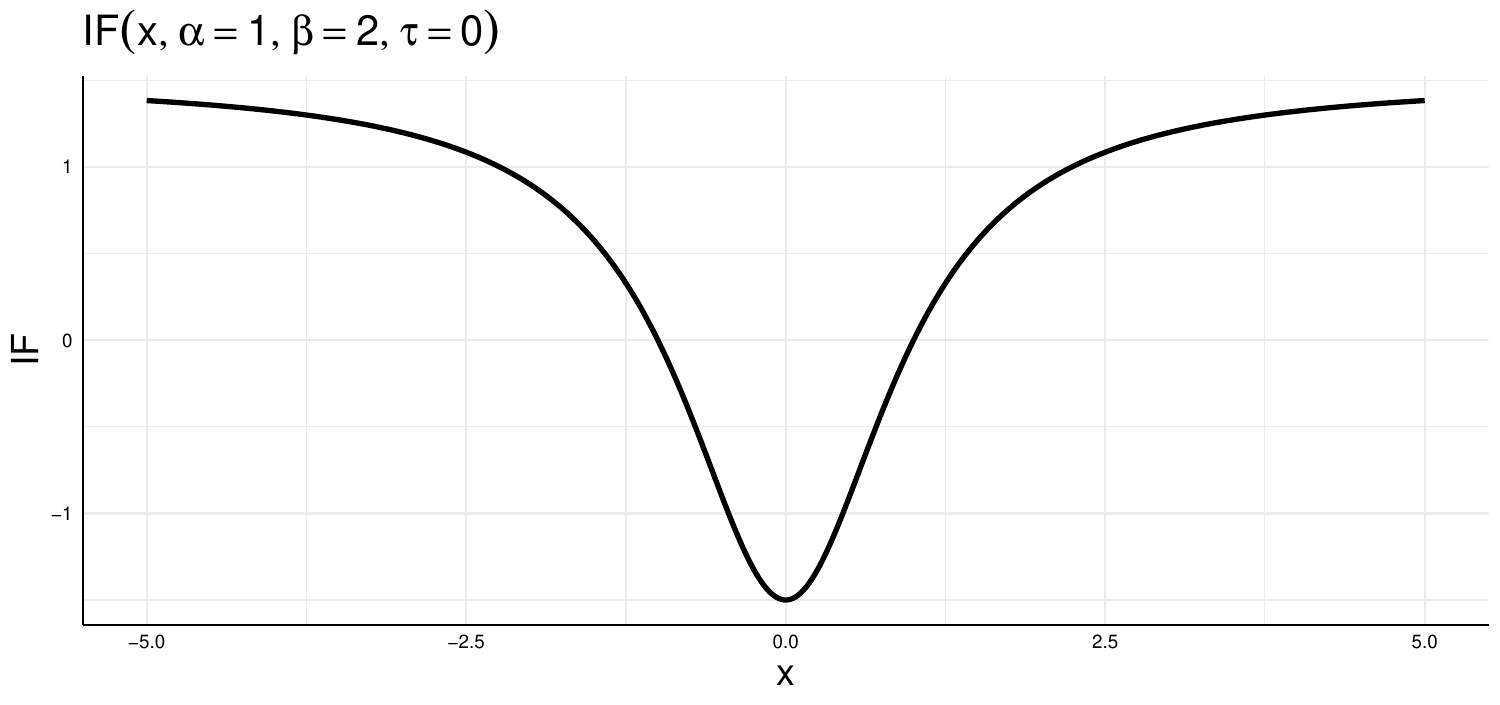}
	\caption{Influence function for the MLE of $\alpha$ }
	\label{fig:IFMLEalpha}
\end{figure}

\begin{figure}
	\centering
	\includegraphics[height=6cm, width=12cm]{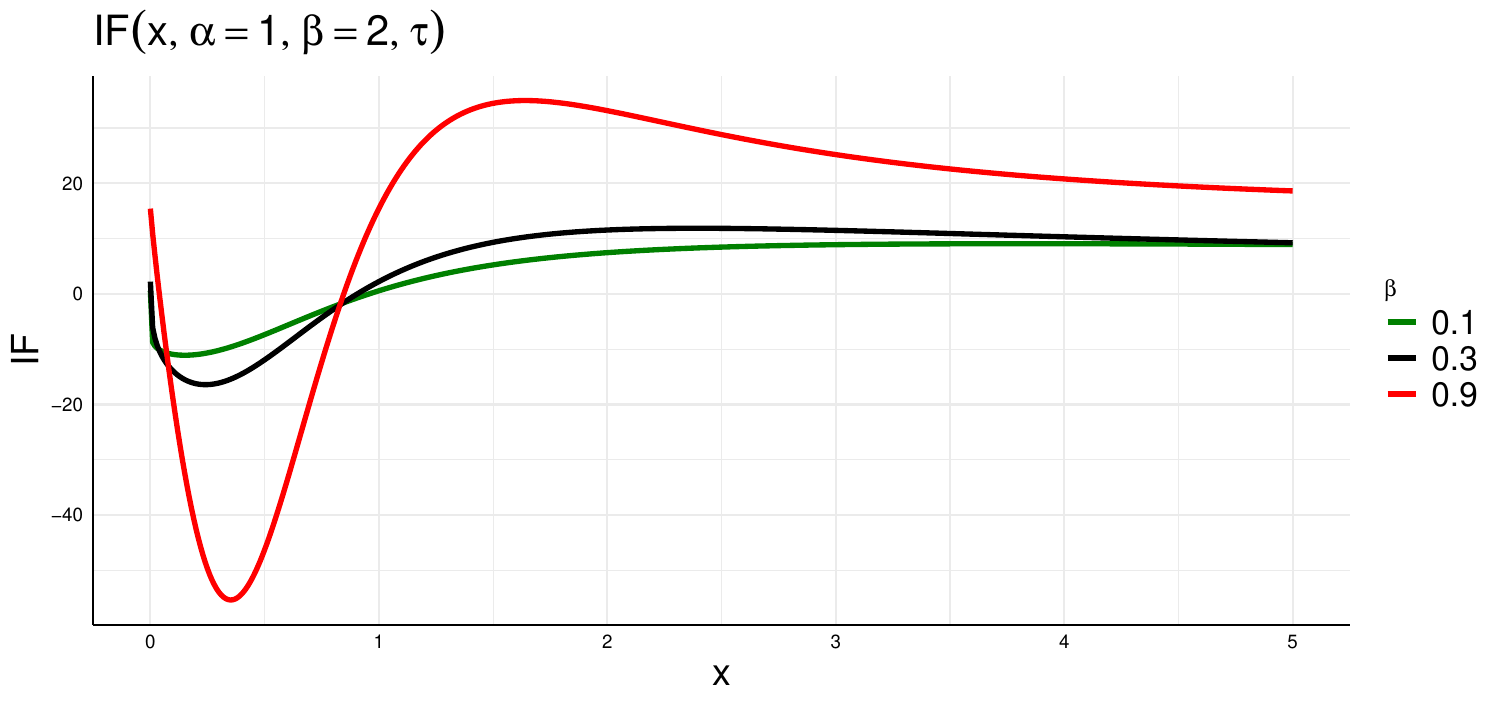}
	\caption{Influence function for the MDPDE of $\alpha$. Black($\tau =0.3$), Red ($0.9$), Green (0.1)}
	\label{fig:IFalpha}
\end{figure}

In Figure \ref{fig:IFMLEalpha} we present the influence function for the MLE of $\alpha ,$ for $%
\alpha =1$ and $\beta =2$ and in Figure \ref{fig:IFalpha} the influence function, of the MDPDEs
of $\alpha ,$ for $\alpha =1$ and $\beta =2$ for different values of the tuning parameter. As shown, the IF of the MLE is unbounded, revealing a lack of robustness. In contrast, the IF of the MDPDE wih positive values of the tuning parameter $\tau$ is bounded, and so MDPDEs of $\alpha$ with $\tau >0$ are locally robust.

Let us now derive the IF for the second parameter $\beta $. The  estimating equation of $\beta $ has been obtained in Section 2. As for the case of $\alpha ,$ we can observe that for $\tau \geq 0$ we have that MDPDEs of $\beta $ are
M-estimators. Therefore, the influence function for the MDPDEof $\beta ,$ with $\tau >0,$
is given by
\begin{eqnarray*}
IF(x,\beta ,\tau ) &=&\frac{1}{J_{\tau }\left( \beta \right) }\left(
1+\frac{1}{\tau }\right) \left\{ \alpha ^{-\tau }\tau \beta ^{\tau } x^{\tau \left( \beta -1\right) }\left( \left( \frac{1}{\beta }+\log
\frac{x}{\alpha }\right) \left( \frac{x^{\beta }}{\alpha ^{\beta }}+1\right)
^{-2\tau }\right. \right. \\
&&\left. \left. -2\frac{x^{\beta }}{\alpha ^{\beta }}\left( \log \frac{x}{\alpha }\right) \left( \frac{x^{\beta }}{\alpha ^{\beta }}+1\right)
^{-2\tau -1}\right) -\alpha ^{-\tau }\tau \beta ^{\tau -1}B\left( \frac{%
\beta \tau +\tau +\beta }{\beta },\frac{\beta \tau -\tau +\beta }{\beta }%
\right) \right. \\
&&\left. \times \left( 1+\frac{\tau }{\beta }\right) \left( \Psi \left( \frac{\beta \tau
-\tau +\beta }{\beta }\right) -\Psi \left( \frac{\beta \tau +\tau +\beta }{%
\beta }\right) \right) \right\} ,
\end{eqnarray*}%
and for $\tau =0,$%
\begin{equation*}
IF(x,\beta ,0)=\frac{\beta ^{2}}{3+\pi ^{2}}\left( \frac{1}{\beta }+\ln \frac{x%
}{\alpha }-2\frac{1}{1+\left( \frac{x}{\alpha }\right) ^{\beta }}\left(
\frac{x}{\alpha }\right) ^{\beta }\ln \frac{x}{\alpha }\right) .
\end{equation*}%

We can observe that

\begin{equation*}
\lim_{x\rightarrow \infty }IF(x,\beta ,0)=-\infty .
\end{equation*}

\begin{figure}
	\centering
	\includegraphics[height=6cm, width=12cm]{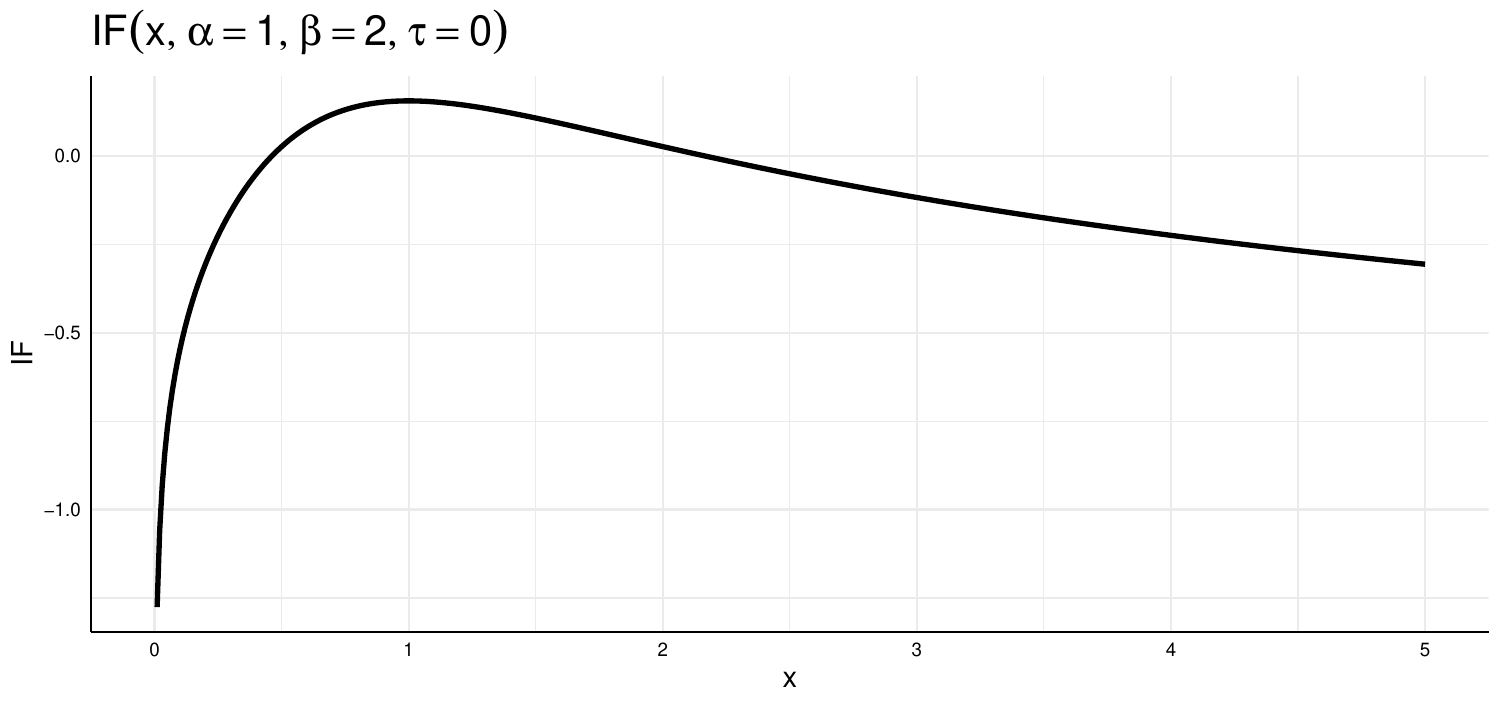}
	\caption{Influence function for the MLE of $\beta$ }
	\label{fig:IFMLEbeta}
\end{figure}

\begin{figure}
	\centering
	\includegraphics[height=6cm, width=12cm]{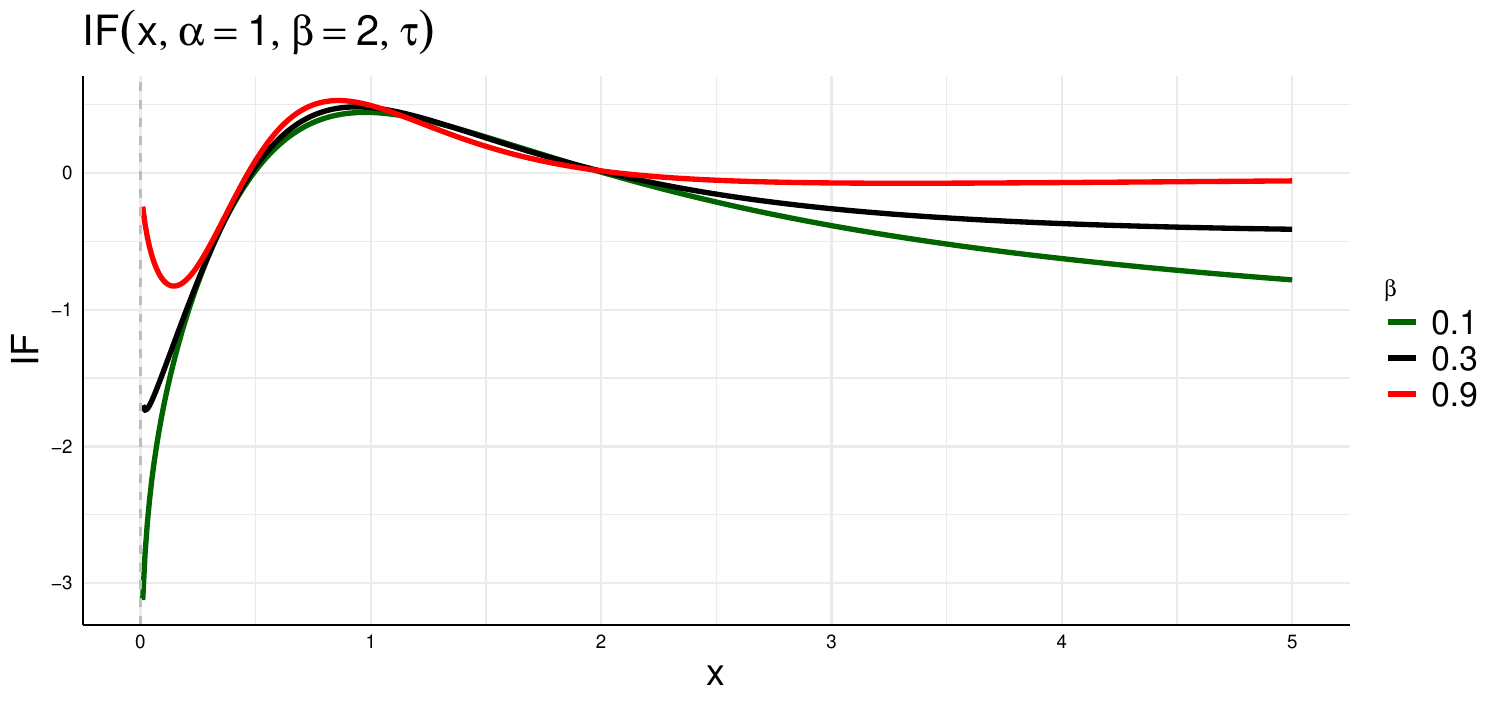}
	\caption{Influence function for the MDPDE of $\beta$. Black($\tau =0.3$), Red ($0.9$), Green ($0.1$)}
	\label{fig:IFbeta}
\end{figure}

In Figures \ref{fig:IFMLEbeta} and \ref{fig:IFbeta} we present the influence functions of the MDPDEs of $\beta
, $ for $\alpha =1$ and $\beta =2$ for $\tau =0$ (Fig. \ref{fig:IFMLEbeta} ) and for different values of the tuning parameter (Fig.  \ref{fig:IFbeta}).

\section{Simulation study}

To assess the performance of DPD in estimating log-logistic parameters, we have conducted a simulation study following the one presented in \cite{mawapa23}. For our simulations, we have considered various log-logistic models and different sample sizes.
For each combination of model ans sample size, we estimate the log-logistic parameters using MLE and MDPDE approaches, with different values of the tuning parameter, namely  $\tau = 0.1, 0.2, ..., 1.0$, as well as other existing estimators that have been proposed in the literature. These include the repeated median estimator (RM), sample median estimator (SM), and the Hodge-Lehmann and Shamos estimator (HL).
Without delving into details, which can be found in \cite{mawapa23}, let us briefly outline the formulas for these estimators.

Consider a sample $(x_1, ..., x_n)$ and assume the log-logistic model holds. Let us denote the ordered sample by $(x_{(1)}, ..., x_{(n)}).$
%
%
For the RM-estimator, we consider the values
\begin{equation*}
y_i= log \left[ {\frac{1}{1-F(x_{(i)})}} -1\right] ,\quad z_i=log (x_{(i)}).
\end{equation*}
Next, we compute
\begin{equation*}
b_1= Med_{1\leq i\leq n} Med_{j\ne i} \left( {\frac{y_i- y_j}{z_i- z_j}}
\right) ,\quad b_0= Med_{1\leq i\leq n} (y_i - b_1 z_i),
\end{equation*}
and finally the RM-estimations for $\alpha $ and $\beta $ are given by
\begin{equation*}
\hat{\alpha }_{RM} = exp (-{\frac{b_0}{b_1}}),\quad \hat{\beta}_{RM}= b_1.
\end{equation*}

To compute the SM and HL estimators, we define $
z_i= log (x_i).
$
Next, we compute
\begin{equation*}
\hat{\mu }_{SM} = Med (z_1, ..., z_n),\quad \hat{\mu }_{HL} =
Med_{i<j}\left( {\frac{z_i + z_j}{2}}\right) ,
\end{equation*}

\begin{equation*}
\hat{s}_{SM} = {\frac{Med (|z_i - \hat{\mu }_{SM}|)}{\Phi^{-1}(3/4)}}, \quad
\hat{s}_{HL} = {\frac{Med_{i<j} (|z_i - z_j|)}{\sqrt{2} \Phi^{-1}(3/4)}},
\end{equation*}
with $\Phi^{-1}$ the inverse of the distribution function of a standard
normal. Then, the SM- and HL- estimations for $\alpha $ and $\beta $ are given, respectively, by

\begin{equation*}
\hat{\alpha }_{SM} = exp (\hat{\mu }_{SM}), \hat{\beta}_{SM}= {\frac{1}{\hat{s}_{SM}}},\quad \quad \hat{\alpha }_{HL} = exp (\hat{%
\mu }_{HL}), \hat{\beta}%
_{HL}= {\frac{1}{\hat{s}_{HL}}}.
\end{equation*}

All estimators are obtained under all considered scenarios. In particular, $%
\alpha $ is fixed to be 1 and $\beta =1.5, 2.5, 5.0, 10.$ The sample sizes range from $n=10$ to $100,$ $n=10, 25, 50, 75,
100.$ For each pair of $(\alpha, \beta )$ we generate a sample of size $n$ for a log-logistic distribution and we estimate the model parameters
applying the different estimators. We repeat the simulations over $M= 10 000$ repetitions. To measure the goodness-of-fit we consider the following measures:
\begin{equation*}
Bias = |\hat{\alpha} - \alpha | + |\hat{\beta } - \beta | ,\, MSE = (%
\hat{\alpha} - \alpha )^2 + (\hat{\beta } - \beta )^2.
\end{equation*}
averaged over the $M=10 000$ simulations 
\begin{equation*}
Bias = {\frac{1}{M}} \sum_{i=1}^M Bias_i ,\, RMSE = \sqrt{{\frac{1}{M}}
\sum_{i=1}^M MSE_i}.
\end{equation*}

Firstly, we consider an uncontaminated scenario. The corresponding error measures are presented in Tables 1 and 2. We have included the results for $\beta = 2.5$ since conclusions for other values of $\beta$ are similar, and the corresponding tables can be found in the Supplementary Material.

From these results, it can be shown that the MLE performs quite well in the absence of contamination. However, note that DPD-based estimators are competitive with the MLE, especially for small values of the tuning parameter $\tau$. This behaviour was expected, as the smaller the tuning parameter, the more efficient the estimator. Also, observe that the behaviour of DPD-based estimators for small values of $\tau$ is better than that of other competing estimators, with only RM being able to compete with them. Moreover, the performance of the estimators improves as the sample size increases.

%
\begin{table}
\begin{center}
\begin{tabular}{|c|cccc|}
\multicolumn{5}{c}{$n=10$}\\
\hline
& Bias & RMSE & $\hat{\alpha }$ & $\hat{\beta }$ \\ \hline
MLE & 0.86117 & 0.99160 & 1.02635 & 2.88258 \\
$DPD_{0.1}$ & 0.87040 & 1.00759 & 1.02438 & 2.88073 \\
$DPD_{0.2}$ & 0.91207 & 1.07743 & 1.02174 & 2.91365 \\
$DPD_{0.3}$ & 0.97792 & 1.19125 & 1.01872 & 2.97222 \\
$DPD_{0.4}$ & 1.06071 & 1.33208 & 1.01545 & 3.04856 \\
$DPD_{0.5}$ & 1.14342 & 1.46068 & 1.01236 & 3.12639 \\
$DPD_{0.6}$ & 1.22792 & 1.58548 & 1.00932 & 3.20736 \\
$DPD_{0.7}$ & 1.30270 & 1.68709 & 1.00638 & 3.27959 \\
$DPD_{0.8}$ & 1.36963 & 1.77269 & 1.00368 & 3.34440 \\
$DPD_{0.9}$ & 1.42582 & 1.84076 & 1.00110 & 3.40043 \\
$DPD_{1.0}$ & 1.47734 & 1.90194 & 0.99867 & 3.45194 \\
RM          & 0.90345 & 1.09261 & 1.01045 & 2.66286 \\
SM          & 1.27049 & 1.21828 & 1.12259 & 1.67069 \\
HL          & 1.21681 & 1.13798 & 1.02703 & 1.51554 \\
%
\hline
\multicolumn{5}{c}{$n=25$}\\
\hline
MLE & 0.48028 & 0.50752 & 1.00983 & 2.63727 \\
$DPD_{0.1}$ & 0.48763 & 0.51759 & 1.00919 & 2.63747 \\
$DPD_{0.2}$ & 0.50554 & 0.54230 & 1.00828 & 2.64756 \\
$DPD_{0.3}$ & 0.52864 & 0.57471 & 1.00722 & 2.66305 \\
$DPD_{0.4}$ & 0.55465 & 0.61302 & 1.00605 & 2.68205 \\
$DPD_{0.5}$ & 0.58227 & 0.65862 & 1.00485 & 2.70358 \\
$DPD_{0.6}$ & 0.60767 & 0.69806 & 1.00369 & 2.72405 \\
$DPD_{0.7}$ & 0.63134 & 0.73463 & 1.00257 & 2.74390 \\
$DPD_{0.8}$ & 0.65404 & 0.77150 & 1.00151 & 2.76356 \\
$DPD_{0.9}$ & 0.67422 & 0.80133 & 1.00055 & 2.78127 \\
$DPD_{1.0}$ & 0.69244 & 0.82689 & 0.99964 & 2.79784 \\
RM & 0.51604 & 0.55662 & 1.01006 & 2.57113 \\
SM & 0.99659 & 0.95598 & 1.01341 & 1.67301 \\
HL & 1.14041 & 1.07281 & 1.01168 &  1.47365 \\
%
\hline
\multicolumn{5}{c}{$n=50$}\\
\hline
& Bias & RMSE & $\hat{\alpha }$ & $\hat{\beta }$ \\ \hline
MLE & 0.32909 & 0.33628 & 1.00623 & 2.56531 \\
$DPD_{0.1}$ & 0.33391 & 0.34216 & 1.00580 & 2.56618 \\
$DPD_{0.2}$ & 0.34529 & 0.35621 & 1.00524 & 2.57183 \\
$DPD_{0.3}$ & 0.35938 & 0.37378 & 1.00463 & 2.57971 \\
$DPD_{0.4}$ & 0.37397 & 0.39244 & 1.00399 & 2.58854 \\
$DPD_{0.5}$ & 0.38798 & 0.41086 & 1.00336 & 2.59760 \\
$DPD_{0.6}$ & 0.40103 & 0.42831 & 1.00275 & 2.60644 \\
$DPD_{0.7}$ & 0.41301 & 0.44449 & 1.00217 & 2.61487 \\
$DPD_{0.8}$ & 0.42391 & 0.45940 & 1.00164 & 2.62283 \\
$DPD_{0.9}$ & 0.43381 & 0.47280 & 1.00114 & 2.63015 \\
$DPD_{1.0}$ & 0.44288 & 0.48480 & 1.00068 & 2.63692 \\
RM & 0.35286 & 0.36535 & 1.00431 & 2.51693 \\
SM & 1.03277 & 0.98680 & 1.02396 & 1.56206 \\
HL & 1.10907 & 1.05263 & 1.00620 & 1.46905 \\ \hline
\end{tabular}%
\end{center}
\par
\label{tablasimsincont50-2}
\caption{Results for $n=10,25$ and $50$ and $\protect\beta = 2.5.$}
\end{table}

\begin{table}
\begin{center}
\begin{tabular}{|c|cccc|}
	\multicolumn{5}{c}{$n=75$}\\
	\hline
& Bias & RMSE & $\hat{\alpha }$ & $\hat{\beta }$ \\ \hline
MLE & 0.26557 & 0.26803 & 1.00389 & 2.54115 \\
$DPD_{0.1}$ & 0.27000 & 0.27299 & 1.00360 & 2.54085 \\
$DPD_{0.2}$ & 0.27886 & 0.28332 & 1.00325 & 2.54332 \\
$DPD_{0.3}$ & 0.28899 & 0.29529 & 1.00286 & 2.54706 \\
$DPD_{0.4}$ & 0.29927 & 0.30739 & 1.00245 & 2.55142 \\
$DPD_{0.5}$ & 0.30906 & 0.31896 & 1.00204 & 2.55603 \\
$DPD_{0.6}$ & 0.31804 & 0.32973 & 1.00162 & 2.56064 \\
$DPD_{0.7}$ & 0.32625 & 0.33960 & 1.00122 & 2.56513 \\
$DPD_{0.8}$ & 0.33368 & 0.34861 & 1.00084 & 2.56944 \\
$DPD_{0.9}$ & 0.34047 & 0.35687 & 1.00048 & 2.57356 \\
$DPD_{1.0}$ & 0.34668 & 0.36449 & 1.00014 & 2.57749 \\
RM & 0.28623 & 0.29136 & 1.00393 & 2.51036 \\
SM & 1.00466 & 0.96143 & 1.00495 & 1.57055 \\
HL & 1.10116 & 1.05026 & 1.00395 & 1.46391 \\
%
\hline
\multicolumn{5}{c}{$n=100$}\\
\hline
MLE & 0.22742 & 0.22921 & 1.00262 & 2.53247 \\
$DPD_{0.1}$ & 0.23106 & 0.23338 & 1.00241 & 2.53353 \\
$DPD_{0.2}$ & 0.23835 & 0.24227 & 1.00215 & 2.53658 \\
$DPD_{0.3}$ & 0.24682 & 0.25275 & 1.00186 & 2.54046 \\
$DPD_{0.4}$ & 0.25548 & 0.26342 & 1.00157 & 2.54465 \\
$DPD_{0.5}$ & 0.26389 & 0.27363 & 1.00129 & 2.54888 \\
$DPD_{0.6}$ & 0.27175 & 0.28309 & 1.00102 & 2.55299 \\
$DPD_{0.7}$ & 0.27900 & 0.29171 & 1.00077 & 2.55689 \\
$DPD_{0.8}$ & 0.28556 & 0.29950 & 1.00053 & 2.56055 \\
$DPD_{0.9}$ & 0.29153 & 0.30655 & 1.00032 & 2.56397 \\
$DPD_{1.0}$ & 0.29700 & 0.31298 & 1.00012 & 2.56718 \\
RM & 0.24559 & 0.25133 & 1.00259 & 2.50775 \\
SM & 1.01447 & 0.97193 & 1.01121 & 1.55085 \\
HL & 1.09114 & 1.04550 & 1.00277 & 1.46494 \\ \hline
\end{tabular}%
\end{center}
\par
\label{tablasimsincont100-2}
\caption{Results for $n=75$ and $100$ and $\protect\beta = 2.5.$}
\end{table}

Because the results remain consistent across various sizes and sample sizes, to investigate performance of the estimators in the presence of contamination, we repeat the previous simulation design with $\alpha=1$, $\beta=10$, and $n=25$ under five different conditions:

\begin{itemize}
\item \textbf{Case 1:} No contamination.

\item \textbf{Case 2:} For the first three elements (12$\% $ of contamination), we change the value in the
sample for a new value coming from a log-logistic distribution with $\alpha
=1, \beta =0.2.$

\item \textbf{Case 3:} For the first three elements (12$\% $ of contamination), we change the value in the
sample for a new value coming from a log-logistic distribution with $\alpha
=4, \beta =10.$

\item \textbf{Case 4:} For the first three elements (12$\% $ of contamination), we change the value in the
sample for a new value coming from a uniform distribution $\mathcal{U}%
(0,20). $

\item \textbf{Case 5:} For the first three elements (12$\% $ of contamination), we change the value in the
sample for a new value of 50.
\end{itemize}

We use the same error measures as before to assess the goodness-of-fit for each estimator. The corresponding values can be found in Tables 3 and 4. The results illustrate that while the MLE performs best in the absence of contamination, small values of $\tau$, specifically 0.1 to 0.3, yield similar results in terms of efficiency. In contrast, in the presence of contamination, the performance of MLE dramatically decreases, whereas estimators based on the DPD are much more robust in all contamination scenarios; the estimations, bias, and RMSE are almost unaffected by the presence of outliers.
Comparing with other competing estimators, namely RM, SM, and HL estimators, we observe that MDPDEs provide superior estimations for both parameters across all scenarios. The only estimator exhibiting similar good properties is RM, but it is outperformed by MDPDEs for small values of $\tau$ (below 0.2-0.3).

In conclusion, these simulations demonstrate that DPD is a highly appealing alternative to MLE. It remains competitive in the absence of contamination and exhibits high stability when contamination arises.

In contrast, one natural challenge lies in determining the optimal value for the tuning parameter $\tau$. Since the best trade-off between efficiency and robustness would depend on the amount of contamination in data, there is no overall optimal value of $\tau.$
 It has been observed in many fields that the best values for $\tau$ are those near 0, and the MDPDE performance tends to worsen significantly for $\tau > 1$. For the log-logistic model, based on our simulations, we would recommend a moderate value or $\tau$, over $0.2$ to $0.3.$

\begin{table}
\begin{center}
\begin{tabular}{|c|cccc|}
	\multicolumn{5}{c}{Case 1}\\
\hline
& Bias & RMSE & $\hat{\alpha }$ & $\hat{\beta }$ \\ \hline
MLE & 1.49171 & 1.95834 & 1.00066 & 10.55203 \\
$DPD_{0.1}$ & 1.51250 & 1.98858 & 1.00061 & 10.55324 \\
$DPD_{0.2}$ & 1.57362 & 2.07727 & 1.00055 & 10.59465 \\
$DPD_{0.3}$ & 1.66093 & 2.21255 & 1.00047 & 10.66236 \\
$DPD_{0.4}$ & 1.76018 & 2.38247 & 1.00039 & 10.74650 \\
$DPD_{0.5}$ & 1.86387 & 2.55809 & 1.00030 & 10.83777 \\
$DPD_{0.6}$ & 1.96730 & 2.73449 & 1.00022 & 10.93175 \\
$DPD_{0.7}$ & 2.06721 & 2.90189 & 1.00015 & 11.02490 \\
$DPD_{0.8}$ & 2.16069 & 3.06012 & 1.00008 & 11.11490 \\
$DPD_{0.9}$ & 2.24385 & 3.19376 & 1.00003 & 11.19776 \\
$DPD_{1.0}$ & 2.31952 & 3.31277 & 0.99998 & 11.27629 \\
RM & 1.61824 & 2.10729 & 1.00065 & 10.14714 \\
SM & 3.48972 & 3.76874 & 1.00064 & 6.72057 \\
HL & 4.11864 & 4.23298 & 1.00061 & 5.91665 \\
\hline
%
\multicolumn{5}{c}{Case 1}\\
\hline
& Bias & RMSE & $\hat{\alpha }$ & $\hat{\beta }$ \\ \hline
MLE  & 8.52751 & 8.46992 & 1.00903 & 1.57094 \\
$DPD_{0.1}$ & 2.07685 & 2.56554 & 1.00035 & 9.08418 \\
$DPD_{0.2}$ & 1.67618 & 2.13926 & 1.00115 & 10.02165 \\
$DPD_{0.3}$ & 1.69258 & 2.19620 & 1.00155 & 10.15721 \\
$DPD_{0.4}$ & 1.74641 & 2.29910 & 1.00196 & 10.18910 \\
$DPD_{0.5}$ & 1.80740 & 2.41018 & 1.00236 & 10.19156 \\
$DPD_{0.6}$ & 1.86526 & 2.50595 & 1.00274 & 10.18188 \\
$DPD_{0.7}$ & 1.91600 & 2.58332 & 1.00317 & 10.16605 \\
$DPD_{0.8}$ & 1.95985 & 2.65006 & 1.00343 & 10.15033 \\
$DPD_{0.9}$ & 1.99572 & 2.69727 & 1.00378 & 10.13157 \\
$DPD_{1.0}$ & 2.02827 & 2.74217 & 1.00436 & 10.11390 \\
RM          & 2.31909 & 2.61352 & 1.00114 & 8.08162 \\
SM          & 4.42249 & 4.60468 & 1.00085 & 5.65700 \\
HL          & 5.74297 & 5.77044 & 1.00186 &  4.29279 \\ \hline


\multicolumn{5}{c}{Case 3}\\
\hline
& Bias & RMSE & $\hat{\alpha }$ & $\hat{\beta }$ \\ \hline
MLE         & 5.45171 & 5.39410 & 1.07107 & 4.62023 \\
$DPD_{0.1}$ & 3.57797 & 3.69628 & 1.02991 & 6.49055 \\
$DPD_{0.2}$ & 2.04492 & 2.43551 & 1.00962 & 9.03039 \\
$DPD_{0.3}$ & 1.82260 & 2.30697 & 1.00482 & 9.93559 \\
$DPD_{0.4}$ & 1.80705 & 2.34653 & 1.00380 & 10.16775 \\
$DPD_{0.5}$ & 1.83445 & 2.41908 & 1.00375 & 10.21575 \\
$DPD_{0.6}$ & 1.87490 & 2.49402 & 1.00399 & 10.21057 \\
$DPD_{0.7}$ & 1.91800 & 2.56804 & 1.00430 & 10.19080 \\
$DPD_{0.8}$ & 1.95858 & 2.63341 & 1.00463 & 10.16763 \\
$DPD_{0.9}$ & 1.99393 & 2.68120 & 1.00494 & 10.14435 \\
$DPD_{1.0}$ & 2.02407 & 2.71758 & 1.00522 & 10.12323 \\
RM          & 2.24717 & 2.54459 & 1.04463 & 8.20681 \\
SM          & 4.49541 & 4.66776 & 1.02943 & 5.58761 \\
HL          & 5.74186 & 5.75108 & 1.04833 & 4.31049 \\ \hline
\end{tabular}%
\end{center}
\par
\label{tablasimconcont3}
\caption{Results for Cases 1,2 and 3.}
\end{table}


\begin{table}
\begin{center}
\begin{tabular}{|c|cccc|}
	\multicolumn{5}{c}{Case 4}\\
\hline
& Bias & RMSE & $\hat{\alpha }$ & $\hat{\beta }$ \\ \hline
MLE         & 6.53474 & 6.49281 & 1.08250 & 3.54993 \\
$DPD_{0.1}$ & 2.83919 & 3.14015 & 1.01501 & 7.52329 \\
$DPD_{0.2}$ & 1.74435 & 2.16942 & 1.00459 & 9.84963 \\
$DPD_{0.3}$ & 1.71037 & 2.19858 & 1.00351 & 10.16774 \\
$DPD_{0.4}$ & 1.75382 & 2.28980 & 1.00346 & 10.22819 \\
$DPD_{0.5}$ & 1.81074 & 2.38837 & 1.00364 & 10.23436 \\
$DPD_{0.6}$ & 1.86913 & 2.49315 & 1.00388 & 10.22698 \\
$DPD_{0.7}$ & 1.91726 & 2.56938 & 1.00414 & 10.21038 \\
$DPD_{0.8}$ & 1.95772 & 2.62786 & 1.00440 & 10.19199 \\
$DPD_{0.9}$ & 1.99289 & 2.68134 & 1.00463 & 10.17533 \\
$DPD_{1.0}$ & 2.02073 & 2.71779 & 1.00484 & 10.15880 \\
RM          & 2.20910 & 2.50561 & 1.04017 & 8.25661 \\
SM          & 4.44278 & 4.61701 & 1.02622 & 5.63958 \\
HL          & 5.70811 & 5.72190 & 1.04321 & 4.34100 \\ \hline
%
%
\multicolumn{5}{c}{Case 5}\\
\hline
& Bias & RMSE & $\hat{\alpha }$ & $\hat{\beta }$ \\ \hline
MLE         & 8.13066 & 7.98137 & 1.15094 & 2.02027 \\
$DPD_{0.1}$ & 1.94249 & 2.38775 & 1.00301 & 9.64431 \\
$DPD_{0.2}$ & 1.65182 & 2.13934 & 1.00147 & 10.41309 \\
$DPD_{0.3}$ & 1.69372 & 2.20955 & 1.00185 & 10.36985 \\
$DPD_{0.4}$ & 1.75495 & 2.31111 & 1.00228 & 10.32227 \\
$DPD_{0.5}$ & 1.81834 & 2.41712 & 1.00270 & 10.28239 \\
$DPD_{0.6}$ & 1.87848 & 2.51433 & 1.00310 & 10.24787 \\
$DPD_{0.7}$ & 1.93483 & 2.60722 & 1.00347 & 10.21949 \\
$DPD_{0.8}$ & 1.98050 & 2.67430 & 1.00380 & 10.19105 \\
$DPD_{0.9}$ & 2.01923 & 2.73096 & 1.00410 & 10.16611 \\
$DPD_{1.0}$ & 2.05149 & 2.77271 & 1.00436 & 10.14442 \\
RM          & 2.25891 & 2.55123 & 1.04410 & 8.21468 \\
SM          & 4.49541 & 4.67479 & 1.02891 & 5.58856 \\
HL          & 5.69134 & 5.70403 & 1.04770 & 4.36077 \\ \hline
\end{tabular}%
\end{center}
\par
\label{tablasimconcont5}
\caption{Results for cases 4 and 5.}
\end{table}


\section{Real data example}

%
%
%
%

In this final section, we apply the estimators developed in the previous sections based on DPD to a situation with real data. For this purpose, we have considered an example considered in \cite{zhchtswa21} and referred to the Scotland's annual maximum flood frequency series in $m^3/s$ from 1952 to 1982. This data are given in Table \ref{tabladatosreales} and were discused in \cite{acsi86} among others.

\begin{table}
\begin{center}
\begin{tabular}{|cccccccc|}
\hline
89.8  & 109.1 & 202.2 & 146.3 & 212.3 & 116.7 & 109.1 & 80.7  \\
127.4 & 138.8 & 283.5 & 85.6  & 105.5 & 118   & 387.8 & 80.7  \\
165.7 & 111.6 & 134.4 & 131.5 & 102   & 104.3 & 242.5 & 214.8 \\
144.6 & 114.2 & 98.3  & 102.8 & 104.3 & 196.2 & 143.7 & \\
\hline
\end{tabular}%
\end{center}
\par
\label{tabladatosreales}
\caption{Annual maximum flood in $m^3/s$ from 1952 to 1982.}
\end{table}

In this set of data, there is a value, namely 387.8, that seems to be an outlier. Assuming the log-logistic distribution, we have obtained the MDPDEs of the parameters $\alpha , \beta $ for several values of the tuning parameter ($\tau $=0, 0.1, ..., 1). To show the robustness of MDPDEs with respect to MLE, we have obtained the estimations in two other situations. In the first one, we have removed datum 387.8, so that no outliers arise in the data. In the second situation, we have substituted datum 387.8 by $387.8\times 5,$ so that a very extreme outlier arise in the data. The results for the three situations can be found in Table \ref{tablasolreal}.

\begin{table}
\begin{center}
\begin{tabular}{|c|cc|cc|cc|}
\hline & \multicolumn{2}{c|}{Original data} & \multicolumn{2}{c|}{Data without outlier} & \multicolumn{2}{c|}{Data with extreme outlier}\\
\hline Method & $\alpha $ & $\beta $ & $\alpha $ & $\beta $  & $\alpha $ & $\beta $ \\
\hline MLE    & 128.59299 & 4.81482  & 125.57231 & 5.38407   & 129.58311 & 4.07479 \\
$DPD_{0.1}$   & 126.41016 & 4.96074  & 124.03480 & 5.43479   & 125.43259 & 4.95591 \\
$DPD_{0.2}$   & 124.22583 & 5.15710  & 122.45253 & 5.54203   & 122.92421 & 5.39308 \\
$DPD_{0.3}$   & 122.13689 & 5.39338  & 120.87111 & 5.69915   & 121.16535 & 5.61685 \\
$DPD_{0.4}$   & 120.27189 & 5.64689  & 119.37588 & 5.89073   & 119.69423 & 5.80014 \\
$DPD_{0.5}$   & 118.73857 & 5.88462  & 118.06536 & 6.09034   & 118.44261 & 5.97323 \\
$DPD_{0.6}$   & 117.56760 & 6.07831  & 117.00301 & 6.26885   & 117.42982 & 6.12397 \\
$DPD_{0.7}$   & 116.71259 & 6.21776  & 116.19007 & 6.40834   & 116.65135 & 6.23973 \\
$DPD_{0.8}$   & 116.09615 & 6.30865  & 115.58392 & 6.50643   & 116.06922 & 6.31889 \\
$DPD_{0.9}$   & 115.64679 & 6.36291  & 115.13113 & 6.57022   & 115.63486 & 6.36764 \\
$DPD_{1.0}$   & 115.31082 & 6.39223  & 114.78594 & 6.60927   & 115.30545 & 6.39440 \\
\hline
\end{tabular}%
\end{center}
\par
\label{tablasolreal}
\caption{Estimating values of the parameters based on DPD for annual maximum flood data in Scotland \protect\cite{zhchtswa21}.}
\end{table}

From this table, we can see that $MDPDEs$ are very robust, especially for big values of the tuning parameter. This was expected, as it is in accordance to the results form the simulation study and it is consistent with other studies based on DPD. Note however, that even for small values of $\tau $, the estimations based on DPD are more stable than the corresponding estimations based on ML.

\section{Conclusions and open problems}

In this paper we have developed a new family of estimators for the parameters of the log-logistic distribution. This new family is based on DPD, a family of divergence that have proved itself to lead to robust estimators in many different situations. This family extends the classical MLE. We have derived the asymptotic distribution of the estimators in three different cases, two of them for one of the parameters and another one for both parameters simultaneously. Besides, we obtain the influence function for each parameter. Finally, we have carried out an extensive simulation study to show the behavior of these new estimators in the presence of contamination. From this simulation study, it seems that MDPDEs are very robust, especially for small values (0.2, 0.3) of the tuning parameter.

Several future research problems merit attention. One such problem involves testing hypotheses on the parameters of the log-logistic distribution. Typically, these hypotheses are tested using the likelihood ratio test based on MLE. However, tests based on non-robust estimators are generally non-robust themselves. Therefore, a promising avenue for future work is the development of test statistics based on MDPDE. We plan to begin by considering simple null hypotheses of the form $H_{0}:\beta =\beta_{0}$ (assuming $\alpha$ is known), $H_{0}:\alpha =\alpha_{0}$ (assuming $\beta$ is known), and $H_{0}:\beta =\beta _{0}, \alpha =\alpha _{0}.$ In this context, we will explore robust Wald-type tests, as considered in \cite{bamamapa16, baghmamapa17, baghmapa18, bamamapa19, baghmamapa21, ghmamapa16, ghbapa21, ghmabapa18}, along with divergence-type tests, as discussed in \cite{bamamapa15, bamamapa18, baghmapa22}.
In a subsequent step, we will investigate composite null hypotheses. In this case, we assume that the null hypothesis can be defined as
\begin{equation}
\Theta _{0}=\left\{ \left( \beta ,\alpha \right) /\boldsymbol{m}(\beta
,\alpha )=\boldsymbol{0}_{r}\right\} , \label{A}
\end{equation}%
where $\boldsymbol{m:}\Theta \rightarrow \mathbb{R}^{r}$ $\left( r\leq
2\right) $ is a differentiable function and its derivative $\boldsymbol{M}%
(\beta ,\alpha )=\frac{\partial \boldsymbol{m}(\beta ,\alpha )}{\partial
\boldsymbol{\theta }}$ ($\boldsymbol{\theta =}(\beta ,\alpha ))$ exists and
is a continuous function in $\boldsymbol{\theta }$ and has $\dim \left(
\boldsymbol{M}(\beta ,\alpha )\right) =2\times r$ and $rank\left(
\boldsymbol{M}(\beta ,\alpha )\right) =r.$
For testing composite null hypotheses as given in (\ref{A}), we will explore Wald-type tests based on MDPDE and divergence-type tests, as discussed in \cite{bamamapa16, baghmamapa17, baghmapa18, bamamapa19, baghmamapa21, ghmamapa16, ghbapa21, ghmabapa18}. In the case of composite null hypotheses, we can also consider Rao-type tests, following the approach outlined in \cite{baghmapa22, mar20, mar20b}. To develop these tests, it will be necessary, beforehand, to introduce MDPDE restricted to the null hypothesis and study its asymptotic distribution, in line with the methods presented in \cite{bamamapa18}. In this context, developing restricted minimum distance estimators is a longstanding problem, initially addressed in \cite{papazo02} using Phi-divergence measures (see also \cite{alsi66, csi63}).

\section*{Acknowledgements}

This work was supported by the Spanish Grant PID2021-124933NB-I00.

\newpage

\section{APPENDIX: Proofs of the results}

\subsection{Proof Lemma \protect\ref{lema1}.}

If we consider
\begin{equation*}
\frac{1}{1+t}=u \Leftrightarrow -\frac{1}{\left( 1+t\right) ^{2}}dt=du \Leftrightarrow dt=-\left(
1+t\right) ^{2}du \Leftrightarrow dt=-u^{-2}du.
\end{equation*}%

On the other hand,

$$1+t=\frac{1}{u} \Leftrightarrow t=\frac{1}{u}-1 \Leftrightarrow t=\frac{1-u}{u}.$$

We have also
that if $t=0$ then $u=1$ and if $t=\infty $ then $u=0.$ Therefore,
\begin{eqnarray*}
I_{1} &=&\int_{0}^{1}\frac{\left( 1-u\right) ^{m(\beta )}}{u^{m(\beta )}}%
u^{s}\frac{1}{u^{2}}du=\int_{0}^{1}\left( 1-u\right) ^{\left( m(\beta
)+1\right) -1}u^{\left( s-m(\beta )-1\right) -1}du \\
&=&B\left( s-m(\beta )-1,m(\beta )+1\right) .
\end{eqnarray*}

\subsection{Proof of Proposition \protect\ref{proposition1}}

We have

\begin{eqnarray*}
\int_{0}^{\infty }f_{\alpha ,\beta }^{1+\tau }(x)dx &=&\int_{0}^{\infty }%
\frac{\beta ^{\tau +1}\alpha ^{\beta \left( \tau +1\right) }x^{\left( \beta
-1\right) \left( \tau +1\right) }}{\left( x^{\beta }+\alpha ^{\beta }\right)
^{2\left( \tau +1\right) }}dx \\
&=&\beta ^{\tau +1}\alpha ^{\beta \left( \tau +1\right) }\alpha ^{\left(
\beta -1\right) \left( \tau +1\right) }\alpha ^{-2\beta \left( \tau
+1\right) }\int_{0}^{\infty }\frac{\left( \frac{x^{\beta }}{\alpha ^{\beta }}%
\right) ^{\frac{\left( \beta -1\right) \left( \tau +1\right) }{\beta }}}{%
\left( \frac{x^{\beta }}{\alpha ^{\beta }}+1\right) ^{2\left( \tau +1\right)
}}dx.
\end{eqnarray*}%

If we consider that
\begin{equation*}
t=\frac{x^{\beta }}{\alpha ^{\beta }}
\end{equation*}%
we have,

$$ dt= \frac{\beta }{\alpha ^{\beta }}x^{\beta -1}dx \Leftrightarrow dx=\frac{\alpha }{\beta }\frac{x^{1-\beta }}{\alpha ^{\beta -1}}dt \Leftrightarrow dx=\frac{\alpha }{\beta }\left( \frac{x^{\beta }}{\alpha
^{\beta }}\right) ^{\frac{1-\beta }{\beta }}dt.$$

Then,
\begin{eqnarray*}
\int_{0}^{\infty }f_{\alpha ,\beta }^{1+\tau }(x)dx &=&\beta ^{\tau
+1}\alpha ^{\beta \left( \tau +1\right) }\alpha ^{\left( \beta -1\right)
\left( \tau +1\right) -2\beta \left( \tau +1\right) }\frac{\alpha }{\beta }%
\int_{0}^{\infty }\frac{t^{\frac{1-\beta }{\beta }\tau }}{\left( 1+t\right)
^{2\left( \tau +1\right) }}dt \\
&=&\left( \frac{\beta }{\alpha }\right) ^{\tau }\int_{0}^{\infty }\frac{t^{%
\frac{\beta -1}{\beta }\tau }}{\left( 1+t\right) ^{2\left( \tau +1\right) }}%
dt.
\end{eqnarray*}%

Based on Lemma \ref{lema1} we have,%
\begin{equation*}
\int_{0}^{\infty }f_{\alpha ,\beta }^{1+\tau }(x)dx=\left( \frac{\beta }{%
\alpha }\right) ^{\tau }B\left( \frac{\beta \tau +\tau +\beta }{\beta },%
\frac{\beta \tau -\tau +\beta }{\beta }\right) .
\end{equation*}

\subsection{Proof Lemma \protect\ref{lema2}}

If we consider the function $g(\beta )=t^{m(\beta )}$ we have
\begin{equation*}
\frac{1}{m^{\prime }(\beta )}\frac{\partial t^{m(\beta )}}{\partial \beta }%
=t^{m(\beta )}\log t.
\end{equation*}%

Therefore,%
\begin{equation*}
I_{2}=\int_{0}^{\infty }\left( \log t\right) \frac{t^{m\left( \beta \right) }%
}{\left( t+1\right) ^{s}}dt=\frac{1}{m^{\prime }\left( \beta \right) }\frac{%
\partial }{\partial \beta }\int_{0}^{\infty }\frac{t^{m\left( \beta \right) }%
}{\left( t+1\right) ^{s}}dt.
\end{equation*}%

If we consider
\begin{equation*}
\frac{1}{1+t}=u
\end{equation*}%
we have,

$$-\frac{1}{\left( 1+t\right) ^{2}}dt=du \Leftrightarrow dt=-\left(
1+t\right) ^{2}du \Leftrightarrow dt=-u^{-2}du.$$

On the other hand,

$$1+t=\frac{1}{u} \Leftrightarrow t=\frac{1}{u}-1 \Leftrightarrow t=\frac{1-u}{u}.$$

We have also
that if $t=0$ then $u=1$ and if $t=\infty $ then $u=0.$ Therefore, by Lemma %
\ref{lema1}
\begin{eqnarray*}
I_{2} &=&\frac{1}{m^{\prime }\left( \beta \right) }\frac{\partial }{\partial
\beta }\int_{0}^{1}(1-u)^{\left( m\left( \beta \right) +1\right)
-1}u^{\left( s-m\left( \beta \right) -1\right) -1}du \\
&=&\frac{1}{m^{\prime }\left( \beta \right) }\frac{B(s-m\left( \beta \right)
-1,m\left( \beta \right) +1)}{B(s-m\left( \beta \right) -1,m\left( \beta
\right) +1)}\frac{\partial }{\partial \beta }B(s-m\left( \beta \right)
-1,m\left( \beta \right) +1) \\
&=&\frac{1}{m^{\prime }\left( \beta \right) }B(s-m\left( \beta \right)
-1,m\left( \beta \right) +1)\frac{\partial \log B(s-m\left( \beta \right)
-1,m\left( \beta \right) +1)}{\partial \beta } \\
&=&\frac{1}{m^{\prime }\left( \beta \right) }B(s-m\left( \beta \right)
-1,m\left( \beta \right) +1)\left[ \frac{\partial }{\partial \beta }\left(
\log \Gamma \left( s-m\left( \beta \right) -1\right) +\log \Gamma \left(
m\left( \beta \right) +1\right) -\log \Gamma \left( s\right) \right) \right]
\\
&=&\frac{1}{m^{\prime }\left( \beta \right) }B(s-m\left( \beta \right)
-1,m\left( \beta \right) +1)\left\{ \frac{\partial \log \Gamma \left(
s-m\left( \beta \right) -1\right) }{\partial \left( s-m\left( \beta \right)
-1\right) }\frac{\partial \left( s-m\left( \beta \right) -1\right) }{%
\partial \beta }\right. \\
&&\left. +\frac{\partial \log \Gamma \left( m\left( \beta \right) +1\right)
}{\partial \left( m\left( \beta \right) +1\right) }\frac{\partial \left(
m\left( \beta \right) +1\right) }{\partial \beta }\right\} \\
&=&\frac{1}{m^{\prime }\left( \beta \right) }B(s-m\left( \beta \right)
-1,m\left( \beta \right) +1)\left\{ -\Psi \left( s-m\left( \beta \right)
-1\right) m^{\prime }\left( \beta \right) +\Psi \left( m\left( \beta \right)
+1\right) m^{\prime }\left( \beta \right) \right\} \\
&=&\frac{1}{m^{\prime }\left( \beta \right) }B(s-m\left( \beta \right)
-1,m\left( \beta \right) +1)\left\{ \Psi \left( m\left( \beta \right)
+1\right) m^{\prime }\left( \beta \right) -\Psi \left( s-m\left( \beta
\right) -1\right) m^{\prime }\left( \beta \right) \right\} \\
&=&B(s-m\left( \beta \right) -1,m\left( \beta \right) +1)\left\{ \Psi \left(
m\left( \beta \right) +1\right) -\Psi \left( s-m\left( \beta \right)
-1\right) \right\} .
\end{eqnarray*}

\subsection{Proof Lemma \protect\ref{lema3}}

We know that%
\begin{equation*}
\frac{\partial t^{m(\beta )}}{\partial \beta }=m^{\prime }(\beta )t^{m(\beta
)}\log t.
\end{equation*}%

Therefore,
\begin{eqnarray*}
\frac{\partial ^{2}t^{m(\beta )}}{\partial \beta ^{2}} &=&\log t\left\{
m^{\prime \prime }(\beta )t^{m(\beta )}+\left( m^{\prime }%
(\beta )\right) ^{2}t^{m(\beta )}\log t\right\} \\
&=&m^{\prime \prime }(\beta )t^{m(\beta )}\log t+\left( \log t\right)
^{2}\left( m%
^{\prime }%
(\beta )\right) ^{2}t^{m(\beta )}.
\end{eqnarray*}%

Hence,
\begin{eqnarray*}
\frac{\partial ^{2}}{\partial \beta ^{2}}\int_{0}^{\infty }\frac{t^{m\left(
\beta \right) }}{\left( t+1\right) ^{s}}dt &=&\left( m%
^{\prime }%
(\beta )\right) ^{2}\int_{0}^{\infty }\left( \log t\right) ^{2}t^{m\left(
\beta \right) }\frac{1}{\left( t+1\right) ^{s}}dt+m^{\prime \prime }(\beta
)\int_{0}^{\infty }\frac{\left( \log t\right) t^{m\left( \beta \right) }}{%
\left( t+1\right) ^{s}}dt \\
&=&\left( m%
^{\prime }%
(\beta )\right) ^{2}I_{3}+m^{\prime \prime }(\beta )I_{2},
\end{eqnarray*}%
and
\begin{equation*}
I_{3}=\frac{1}{\left( m%
^{\prime }%
(\beta )\right) ^{2}}\left( \frac{\partial ^{2}}{\partial \beta ^{2}}%
\int_{0}^{\infty }\frac{t^{m\left( \beta \right) }}{\left( t+1\right) ^{s}}%
dt-m^{\prime \prime }(\beta )I_{2}\right) .
\end{equation*}

Now we have,%
\begin{equation*}
\frac{\partial }{\partial \beta }\int_{0}^{\infty }\frac{t^{m\left( \beta
\right) }}{\left( t+1\right) ^{s}}dt=m%
^{\prime }%
(\beta )\left( B(s-m\left( \beta \right) -1,m\left( \beta \right) +1)\right)
\left\{ \Psi \left( m\left( \beta \right) +1\right) -\Psi \left( s-m\left(
\beta \right) -1\right) \right\} ,
\end{equation*}%
and
\begin{eqnarray*}
\frac{\partial ^{2}}{\partial \beta ^{2}}\int_{0}^{\infty }\frac{t^{m\left(
\beta \right) }}{\left( t+1\right) ^{s}}dt &=&m^{\prime \prime }(\beta
)B(s-m\left( \beta \right) -1,m\left( \beta \right) +1)\left\{ \left( \Psi
\left( m\left( \beta \right) +1\right) \right. \right. \\
&&\left. -\Psi \left( s-m\left( \beta \right) -1\right) \right\} \\
&&\left. +\left( m%
^{\prime }%
(\beta )\right) ^{2}B(s-m\left( \beta \right) -1,m\left( \beta \right)
+1)\left( \Psi \left( m\left( \beta \right) +1\right) -\Psi \left( s-m\left(
\beta \right) -1\right) \right) ^{2}\right. \\
&&\left. +\left( m%
^{\prime }%
(\beta )\right) ^{2}B(s-m\left( \beta \right) -1,m\left( \beta \right)
+1)\left( \Psi \left( m\left( \beta \right) +1\right) +\Psi \left( s-m\left(
\beta \right) -1\right) \right) \right\} .
\end{eqnarray*}%

Therefore,%
\begin{eqnarray*}
I_{3} &=&\frac{1}{\left( m%
^{\prime }%
(\beta )\right) ^{2}}\left( \frac{\partial ^{2}}{\partial \beta ^{2}}%
\int_{0}^{\infty }\frac{t^{m\left( \beta \right) }}{\left( t+1\right) ^{s}}%
dt-m^{\prime \prime }(\beta )I_{2}\right) \\
&=&\frac{1}{\left( m%
^{\prime }%
(\beta )\right) ^{2}}\left\{ \left( m%
^{\prime }%
(\beta )\right) ^{2}B(s-m\left( \beta \right) -1,m\left( \beta \right)
+1)\left( \Psi \left( m\left( \beta \right) +1\right) -\Psi \left( s-m\left(
\beta \right) -1\right) \right) ^{2}\right. \\
&&\left. +\left( m%
^{\prime }%
(\beta )\right) ^{2}B(s-m\left( \beta \right) -1,m\left( \beta \right)
+1)\left( \Psi ^{\prime }\left( m\left( \beta \right) +1\right) +\Psi
^{\prime }\left( s-m\left( \beta \right) -1\right) \right) \right\} \\
&=&B(s-m\left( \beta \right) -1,m\left( \beta \right) +1) \\
& & \left\{ \left( \Psi
\left( m\left( \beta \right) +1\right) -\Psi \left( s-m\left( \beta \right)
-1\right) \right) ^{2}+\left( \Psi ^{\prime }\left( m\left( \beta \right)
+1\right) +\Psi ^{\prime }\left( s-m\left( \beta \right) -1\right) \right)
\right\} .
\end{eqnarray*}

\subsection{Proof Theorem \protect\ref{Theorem1}}

It is clear that

\begin{equation*}
\log f_{\alpha }(x)=\log \beta +\beta \log \alpha +\left( \beta -1\right)
\log x-2\log \left( x^{\beta }+\alpha ^{\beta }\right) .
\end{equation*}%

Therefore,

\begin{equation*}
\frac{\partial \log f_{\alpha }(x)}{\partial \alpha }=\frac{\beta }{\alpha }%
-2\frac{\beta \alpha ^{\beta -1}}{\left( x^{\beta }+\alpha ^{\beta }\right) } ,
\end{equation*}%
and
\begin{eqnarray*}
\left( \frac{\partial \log f_{\alpha }(x)}{\partial \alpha }\right) ^{2}
&=&\left( \frac{\beta }{\alpha }\right) ^{2}+4\frac{\beta ^{2}\alpha
^{2\left( \beta -1\right) }}{\left( x^{\beta }+\alpha ^{\beta }\right) ^{2}}%
-4\frac{\beta }{\alpha }\frac{\beta \alpha ^{\beta -1}}{\left( x^{\beta
}+\alpha ^{\beta }\right) } \\
&=&\left( \frac{\beta }{\alpha }\right) ^{2}+4\frac{\beta ^{2}\alpha
^{2\left( \beta -1\right) }}{\left( x^{\beta }+\alpha ^{\beta }\right) ^{2}}%
-4\frac{\beta ^{2}\alpha ^{\beta -2}}{\left( x^{\beta }+\alpha ^{\beta
}\right) }.
\end{eqnarray*}%

Now,%
\begin{eqnarray*}
J_{\tau }(\alpha ) &=&\left( \frac{\beta }{\alpha }\right)
^{2}\int_{0}^{\infty }f_{\alpha }(x)^{\tau +1}dx+4\beta ^{2}\alpha
^{2\left( \beta -1\right) }\int_{0}^{\infty }\frac{1}{\left(
x^{\beta }+\alpha ^{\beta }\right) ^{2}}f_{\alpha }(x)^{\tau +1}dx \\
&&-4\beta ^{2}\alpha ^{\beta -2}\int_{0}^{\infty }\frac{1}{x^{\beta
}+\alpha ^{\beta }}f_{\alpha }(x)^{\tau +1}dx \\
&=&L_{1}+L_{2}-L_{3}.
\end{eqnarray*}

\underline{In relation to $L_{1}$ we have:}%
\begin{eqnarray*}
L_{1} &=&\left( \frac{\beta }{\alpha }\right) ^{2}\int_{0}^{\infty }%
\frac{\beta ^{\tau +1}\alpha ^{\beta \left( \tau +1\right) }x^{\left( \beta
-1\right) \left( \tau +1\right) }}{\left( x^{\beta }+\alpha ^{\beta }\right)
^{^{2\left( \tau +1\right) }}}dx=\left( \frac{\beta }{\alpha }\right)
^{2}\beta ^{\tau +1}\alpha ^{\beta \left( \tau +1\right) } \\
&&\times \int_{0}^{\infty }\frac{x^{\left( \beta -1\right) \left(
\tau +1\right) }}{\left( \frac{x^{\beta }}{\alpha ^{\beta }}+1\right)
^{2\left( \tau +1\right) }\alpha ^{2\beta \left( \tau +1\right) }}dx \\
&=&\left( \frac{\beta }{\alpha }\right) ^{2}\beta ^{\tau +1}\alpha ^{\beta
\left( \tau +1\right) }\alpha ^{-2\beta \left( \tau +1\right)
}\int_{0}^{\infty }\frac{x^{\left( \beta -1\right) \tau }}{\left(
\frac{x^{\beta }}{\alpha ^{\beta }}+1\right) ^{2\left( \tau +1\right) }}%
x^{\beta -1}dx \\
&=&\left( \frac{\beta }{\alpha }\right) ^{2}\beta ^{\tau +1}\alpha ^{\beta
\left( \tau +1\right) }\alpha ^{-2\beta \left( \tau +1\right) }\alpha
^{\left( \beta -1\right) \tau }\int_{0}^{\infty }\frac{\left( \frac{x%
}{\alpha }\right) ^{\frac{\left( \beta -1\right) \tau }{\beta }}}{\left(
\frac{x^{\beta }}{\alpha ^{\beta }}+1\right) ^{2\left( \tau +1\right) }}%
x^{\beta -1}dx.
\end{eqnarray*}%

If we consider
\begin{equation*}
t=\frac{x^{\beta }}{\alpha ^{\beta }}
\end{equation*}%
we have,

$$ dt=\frac{\beta }{\alpha ^{\beta }}x^{\beta -1}dx \Leftrightarrow x^{\beta -1}dx=\frac{\alpha ^{\beta }}{\beta }dt.$$

Then,%
\begin{equation*}
L_{1}=\left( \frac{\beta }{\alpha }\right) ^{\tau
+2}\int_{0}^{\infty }\frac{t^{\frac{\left( \beta -1\right) \tau }{%
\beta }}}{\left( t+1\right) ^{2\left( \tau +1\right) }}dt.
\end{equation*}

By Lemma \ref{lema1} with $m(\beta )=\frac{\left( \beta -1\right) \tau }{%
\beta }$ and $s=2\left( \tau +1\right) ,$ we have
\begin{equation*}
\int_{0}^{\infty }\frac{t^{\frac{\left( \beta -1\right) \tau }{\beta
}}}{\left( t+1\right) ^{2\left( \tau +1\right) }}dt=B\left( \frac{\beta \tau
+\tau +\beta }{\beta },\frac{\tau \beta -\tau +\beta }{\beta }\right) .
\end{equation*}%

Therefore,
\begin{equation*}
L_{1}=\left( \frac{\beta }{\alpha }\right) ^{\tau +2}B\left( \frac{\beta
\tau +\tau +\beta }{\beta },\frac{\tau \beta -\tau +\beta }{\beta }\right) .
\end{equation*}

For $\tau =0$ we have,
\begin{equation*}
L_{1}=\left( \frac{\beta }{\alpha }\right) ^{2}.
\end{equation*}%

\underline{In relation to $L_{2}$ we have:}%
\begin{eqnarray*}
L_{2} &=&4\beta ^{2}\alpha ^{2\left( \beta -1\right)
}\int_{0}^{\infty }\frac{1}{\left( x^{\beta }+\alpha ^{\beta
}\right) ^{2}}f_{\alpha }(x)^{\tau +1}dx \\
&=&4\beta ^{2}\alpha ^{2\left( \beta -1\right) }\int_{0}^{\infty }%
\frac{1}{\left( x^{\beta }+\alpha ^{\beta }\right) ^{2}}\frac{\beta ^{\tau
+1}\alpha ^{\beta \left( \tau +1\right) }x^{\left( \beta -1\right) \left(
\tau +1\right) }}{\left( x^{\beta }+\alpha ^{\beta }\right) ^{^{2\left( \tau
+1\right) }}}dx \\
&=&4\beta ^{2}\alpha ^{2\left( \beta -1\right) }\beta ^{\tau +1}\alpha
^{\beta \left( \tau +1\right) }\int_{0}^{\infty }\frac{x^{\left(
\beta -1\right) \left( \tau +1\right) }}{\left( x^{\beta }+\alpha ^{\beta
}\right) ^{2\tau +4}}dx \\
&=&4\beta ^{2}\alpha ^{2\left( \beta -1\right) }\beta ^{\tau +1}\alpha
^{\beta \left( \tau +1\right) }\alpha ^{-\beta \left( 2\tau +4\right)
}\int_{0}^{\infty }\frac{x^{\left( \beta -1\right) \tau }}{\left(
\frac{x^{\beta }}{\alpha ^{\beta }}+1\right) ^{2\tau +4}}x^{\beta -1}dx \\
&=&4\beta ^{2}\alpha ^{2\left( \beta -1\right) }\beta ^{\tau +1}\alpha
^{\beta \left( \tau +1\right) }\alpha ^{-\beta \left( 2\tau +4\right)
}\alpha ^{\left( \beta -1\right) \tau }\int_{0}^{\infty }\frac{%
\left( \frac{x^{\beta }}{\alpha ^{\beta }}\right) ^{\frac{\left( \beta
-1\right) \tau }{\beta }}}{\left( \frac{x^{\beta }}{\alpha ^{\beta }}%
+1\right) ^{2\tau +4}\alpha ^{\beta \left( 2\tau +4\right) }}x^{\beta -1}dx.
\end{eqnarray*}%

If we consider
\begin{equation*}
t=\frac{x^{\beta }}{\alpha ^{\beta }}
\end{equation*}%
we have,

$$dt= \frac{\beta }{\alpha ^{\beta }}x^{\beta -1}dx \Leftrightarrow x^{\beta -1}dx=\frac{\alpha ^{\beta }}{\beta }dt.$$

Then,%
\begin{equation*}
L_{2}=4\beta ^{2}\alpha ^{2\left( \beta -1\right) }\beta ^{\tau +1}\alpha
^{\beta \left( \tau +1\right) }\alpha ^{-\beta \left( 2\tau +4\right)
}\alpha ^{\left( \beta -1\right) \tau }\frac{\alpha ^{\beta }}{\beta }%
\int_{0}^{\infty }\frac{t^{\frac{\left( \beta -1\right) \tau }{\beta
}}}{\left( t+1\right) ^{2\tau +4}}dt.
\end{equation*}

By Lemma \ref{lema1} with $m(\beta )=\frac{\left( \beta -1\right) \tau }{%
\beta }$ and $s=2\tau +4,$ we have%
\begin{equation*}
\int_{0}^{\infty }\frac{t^{\frac{\left( \beta -1\right) \tau }{\beta
}}}{\left( t+1\right) ^{2\tau +4}}dt=B\left( \frac{\beta \tau +\tau +3\beta
}{\beta },\frac{\tau \beta -\tau +\beta }{\beta }\right) .
\end{equation*}%

Therefore,
\begin{equation*}
L_{2}=4\left( \frac{\beta }{\alpha }\right) ^{\tau +2}B\left( \frac{\beta
\tau +\tau +3\beta }{\beta },\frac{\tau \beta -\tau +\beta }{\beta }\right) .
\end{equation*}%

But
\begin{equation*}
B\left( \frac{\beta \tau +\tau +3\beta }{\beta },\frac{\tau \beta -\tau
+\beta }{\beta }\right) =\frac{1}{2}\frac{\left( \beta \tau +\tau +2\beta
\right) \left( \beta \tau +\tau +\beta \right) }{\beta ^{2}\left( 2\tau
+3\right) \left( \tau +1\right) }B\left( \frac{\beta \tau +\tau +\beta }{%
\beta },\frac{\tau \beta -\tau +\beta }{\beta }\right) .
\end{equation*}%

Hence,
\begin{equation*}
L_{2}=4\left( \frac{\beta }{\alpha }\right) ^{\tau +2}\frac{1}{2}\frac{%
\left( \beta \tau +\tau +2\beta \right) \left( \beta \tau +\tau +\beta
\right) }{\beta ^{2}\left( 2\tau +3\right) \left( \tau +1\right) }B\left(
\frac{\beta \tau +\tau +\beta }{\beta },\frac{\tau \beta -\tau +\beta }{%
\beta }\right) .
\end{equation*}

For $\tau =0$ we have%
\begin{equation*}
L_{2}=\left( \frac{\beta }{\alpha }\right) ^{2}\frac{4}{3}.
\end{equation*}%

\underline{In relation to $L_{3}$ we have:}%
\begin{eqnarray*}
L_{3} &=&4\beta ^{2}\alpha ^{\beta -2}\int_{0}^{\infty }\frac{1}{%
x^{\beta }+\alpha ^{\beta }}f_{\alpha }(x)^{\tau +1}dx \\
&=&4\beta ^{2}\alpha ^{\beta -2}\int_{0}^{\infty }\frac{1}{x^{\beta
}+\alpha ^{\beta }}\frac{\beta ^{\tau +1}\alpha ^{\beta \left( \tau
+1\right) }x^{\left( \beta -1\right) \left( \tau +1\right) }}{\left(
x^{\beta }+\alpha ^{\beta }\right) ^{^{2\left( \tau +1\right) }}}dx \\
&=&4\beta ^{2}\alpha ^{\beta -2}\beta ^{\tau +1}\alpha ^{\beta \left( \tau
+1\right) }\alpha ^{-2\beta (\tau +1)}\int_{0}^{\infty }\frac{%
x^{\left( \beta -1\right) \tau }}{\left( \frac{x^{\beta }}{\alpha ^{\beta }}%
+1\right) ^{2\tau +3}}x^{\beta -1}dx \\
&=&4\beta ^{2}\alpha ^{\beta -2}\beta ^{\tau +1}\alpha ^{\beta \left( \tau
+1\right) }\alpha ^{-2\beta (\tau +1)}\alpha ^{\left( \beta -1\right) \tau
}\int_{0}^{\infty }\frac{\left( \frac{x^{\beta }}{\alpha ^{\beta }}%
\right) ^{\frac{\left( \beta -1\right) \tau }{\beta }}}{\left( \frac{%
x^{\beta }}{\alpha ^{\beta }}+1\right) ^{2\tau +3}}x^{\beta -1}dx.
\end{eqnarray*}%

If we consider
\begin{equation*}
t=\frac{x^{\beta }}{\alpha ^{\beta }}
\end{equation*}%
we have,

$$dt= \frac{\beta }{\alpha ^{\beta }}x^{\beta -1}dx \Leftrightarrow x^{\beta -1}dx=\frac{\alpha ^{\beta }}{\beta }dt.$$

Then,%
\begin{equation*}
L_{3}=4\left( \frac{\beta }{\alpha }\right) ^{\tau
+2}\int_{0}^{\infty }\frac{t^{\frac{\left( \beta -1\right) \tau }{%
\beta }}}{\left( t+1\right) ^{2\tau +3}}dt.
\end{equation*}

By Lemma \ref{lema1} with $m(\beta )=\frac{\left( \beta -1\right) \tau }{%
\beta }$ and $s=2\tau +3,$ we have%
\begin{equation*}
\int_{0}^{\infty }\frac{t^{\frac{\left( \beta -1\right) \tau }{\beta
}}}{\left( t+1\right) ^{2\tau +3}}dt=B\left( \frac{\beta \tau +\tau +2\beta
}{\beta },\frac{\tau \beta -\tau +\beta }{\beta }\right) .
\end{equation*}%

Therefore,%
\begin{equation*}
L_{3}=4\left( \frac{\beta }{\alpha }\right) ^{\tau +2}B\left( \frac{\beta
\tau +\tau +2\beta }{\beta },\frac{\tau \beta -\tau +\beta }{\beta }\right) .
\end{equation*}%

But%
\begin{equation*}
B\left( \frac{\beta \tau +\tau +2\beta }{\beta },\frac{\tau \beta -\tau
+\beta }{\beta }\right) =\frac{1}{2}\frac{\beta \tau +\tau +\beta }{\beta
\left( \tau +1\right) }B\left( \frac{\beta \tau +\tau +\beta }{\beta },\frac{%
\tau \beta -\tau +\beta }{\beta }\right) .
\end{equation*}%

Hence,%
\begin{eqnarray*}
L_{3} & = & 4\left( \frac{\beta }{\alpha }\right) ^{\tau +2}B\left( \frac{\beta
\tau +\tau +2\beta }{\beta },\frac{\tau \beta -\tau +\beta }{\beta }\right) \\
& = & 4\left( \frac{\beta }{\alpha }\right) ^{2+\tau }\frac{1}{2}\frac{\beta \tau
+\tau +\beta }{\beta \left( \tau +1\right) }B\left( \frac{\beta \tau +\tau
+\beta }{\beta },\frac{\tau \beta -\tau +\beta }{\beta }\right) .
\end{eqnarray*}

For $\tau =0,$ we have
\begin{equation*}
L_{3}=2\left( \frac{\beta }{\alpha }\right) ^{2}.
\end{equation*}

Based on the results obtained for $L_{1},$ $L_{2}$ and $L_{3}$ we have,
\begin{eqnarray*}
J_{\tau }(\alpha ) &=&\left( \frac{\beta }{\alpha }\right) ^{\tau +2}B\left(
\frac{\beta \tau +\tau +\beta }{\beta },\frac{\tau \beta -\tau +\beta }{%
\beta }\right) +4\left( \frac{\beta }{\alpha }\right) ^{\tau +2}\frac{1}{2}%
\frac{\left( \beta \tau +\tau +2\beta \right) \left( \beta \tau +\tau +\beta
\right) }{\beta ^{2}\left( 2\tau +3\right) \left( \tau +1\right) } \\
&&\times B\left( \frac{\beta \tau +\tau +\beta }{\beta },\frac{\tau \beta
-\tau +\beta }{\beta }\right) \\
&&-4\left( \frac{\beta }{\alpha }\right) ^{\tau +2}\frac{1}{2}\frac{\beta
\tau +\tau +\beta }{\beta \left( \tau +1\right) }B\left( \frac{\beta \tau
+\tau +\beta }{\beta },\frac{\tau \beta -\tau +\beta }{\beta }\right) \\
&=&\left( \frac{\beta }{\alpha }\right) ^{\tau +2}B\left( \frac{\beta \tau
+\tau +\beta }{\beta },\frac{\tau \beta -\tau +\beta }{\beta }\right) \\
& & \left\{ 1+4\frac{1}{2}\frac{\left( \beta \tau +\tau +2\beta \right) \left(
\beta \tau +\tau +\beta \right) }{\beta ^{2}\left( 2\tau +3\right) \left(
\tau +1\right) }-4\frac{1}{2}\frac{\beta \tau +\tau +\beta }{\beta \left(
\tau +1\right) }\right\} \\
&=&\left( \frac{\beta }{\alpha }\right) ^{\tau +2}B\left( \frac{\beta \tau
+\tau +\beta }{\beta },\frac{\tau \beta -\tau +\beta }{\beta }\right)
\left\{ 1+2\frac{\left( \beta \tau +\tau +2\beta \right) \left( \beta \tau
+\tau +\beta \right) }{\beta ^{2}\left( \tau +1\right) \left( 2\tau
+3\right) }-2\frac{\beta \tau +\tau +\beta }{\beta \left( \tau +1\right) }%
\right\} \\
&=&\left( \frac{\beta }{\alpha }\right) ^{\tau +2}B\left( \frac{\beta \tau
+\tau +\beta }{\beta },\frac{\tau \beta -\tau +\beta }{\beta }\right) \left[
2\frac{\left( \beta \tau +\tau +\beta \right) \left( -\tau \beta -\beta
+\tau \right) }{\beta ^{2}\left( \tau +1\right) \left( 2\tau +3\right) }+1%
\right] .
\end{eqnarray*}

And for $\tau =0,$ we obtain

$$ J_{\tau }(\alpha ) = \left( \frac{\beta }{\alpha }\right)
^{2}(1+\frac{4}{3}-2)=\frac{\beta ^{2}}{3\alpha ^{2}}.$$

\subsection{Proof Theorem \protect\ref{Theorem1a}}

We know
\begin{eqnarray*}
\xi _{\tau }\left( \alpha \right) &=&\int_{0}^{\infty }\frac{%
\partial \log f_{\alpha }(x)}{\partial \alpha }f_{\alpha }(x)^{\tau +1}dx \\
&=&\int_{0}^{\infty }\left( \frac{\beta }{\alpha }-2\frac{\beta
\alpha ^{\beta -1}}{\left( x^{\beta }+\alpha ^{\beta }\right) }\right) \frac{%
\beta ^{\tau +1}\alpha ^{\beta \left( \tau +1\right) }x^{\left( \beta
-1\right) \left( \tau +1\right) }}{\left( x^{\beta }+\alpha ^{\beta }\right)
^{^{2\left( \tau +1\right) }}}dx \\
&=&\frac{\beta }{\alpha }\int_{0}^{\infty }\frac{\beta ^{\tau
+1}\alpha ^{\beta \left( \tau +1\right) }x^{\left( \beta -1\right) \left(
\tau +1\right) }}{\left( x^{\beta }+\alpha ^{\beta }\right) ^{^{2\left( \tau
+1\right) }}}dx-2\beta \alpha ^{\beta -1}\int_{0}^{\infty }\frac{1}{%
\left( x^{\beta }+\alpha ^{\beta }\right) }\frac{\beta ^{\tau +1}\alpha
^{\beta \left( \tau +1\right) }x^{\left( \beta -1\right) \left( \tau
+1\right) }}{\left( x^{\beta }+\alpha ^{\beta }\right) ^{^{2\left( \tau
+1\right) }}}dx \\
&=&M_{1}+M_{2}.
\end{eqnarray*}%

\underline{In relation to $M_{1}$ we have }

\begin{eqnarray*}
M_{1} &=&\frac{\beta }{\alpha }\int_{0}^{\infty }\frac{\beta ^{\tau
+1}\alpha ^{\beta \left( \tau +1\right) }x^{\left( \beta -1\right) \left(
\tau +1\right) }}{\left( x^{\beta }+\alpha ^{\beta }\right) ^{^{2\left( \tau
+1\right) }}}dx \\
&=&\frac{\beta }{\alpha }\beta ^{\tau +1}\alpha ^{\beta \left( \tau
+1\right) }\alpha ^{-2\beta \left( \tau +1\right) }\alpha ^{\left( \beta
-1\right) \tau }\int_{0}^{\infty }\frac{\left( \frac{x^{\beta }}{%
\alpha ^{\beta }}\right) ^{\frac{\left( \beta -1\right) \tau }{\beta }}}{%
\left( \frac{x^{\beta }}{\alpha ^{\beta }}+1\right) ^{2\tau +2}}x^{\beta
-1}dx.
\end{eqnarray*}%

If we consider
\begin{equation*}
t=\frac{x^{\beta }}{\alpha ^{\beta }}
\end{equation*}%
we have,

$$dt=\frac{\beta }{\alpha ^{\beta }}x^{\beta -1}dx \Leftrightarrow x^{\beta -1}dx=\frac{\alpha ^{\beta }}{\beta }dt.$$

Then,

\begin{equation*}
M_{1}=\left( \frac{\beta }{\alpha }\right) ^{\tau
+1}\int_{0}^{\infty }\frac{t^{\frac{\left( \beta -1\right) \tau }{%
\beta }}}{\left( t+1\right) ^{2\tau +2}}dt.
\end{equation*}

By Lemma \ref{lema1} with $m(\beta )=\frac{\left( \beta -1\right) \tau }{%
\beta }$ and $s=2\tau +2,$ we have,%
\begin{equation*}
\int_{0}^{\infty }\frac{t^{\frac{\left( \beta -1\right) \tau }{\beta
}}}{\left( t+1\right) ^{2\tau +2}}dt=B\left( \frac{\tau \beta +\beta +\tau }{%
\beta },\frac{\tau \beta +\beta -\tau }{\beta }\right) .
\end{equation*}%

Therefore,%
\begin{equation*}
M_{1}=\left( \frac{\beta }{\alpha }\right) ^{\tau +1}B\left( \frac{\tau
\beta +\beta +\tau }{\beta },\frac{\tau \beta +\beta -\tau }{\beta }\right) .
\end{equation*}

\underline{In relation to $M_{2}$ we have }%
\begin{eqnarray*}
M_{2} &=&-2\beta \alpha ^{\beta -1}\beta ^{\tau +1}\alpha ^{\beta \left(
\tau +1\right) }\int_{0}^{\infty }\frac{x^{\left( \beta -1\right)
\left( \tau +1\right) }}{\left( x^{\beta }+\alpha ^{\beta }\right) ^{^{2\tau
+3}}}dx \\
&=&-2\beta \alpha ^{\beta -1}\beta ^{\tau +1}\alpha ^{\beta \left( \tau
+1\right) }\alpha ^{-\beta \left( 2\tau +3\right) }\int_{0}^{\infty }%
\frac{x^{\left( \beta -1\right) \tau }}{\left( \frac{x^{\beta }}{\alpha
^{\beta }}+1\right) ^{2\tau +3}}x^{\beta -1}dx.
\end{eqnarray*}%

If we consider
\begin{equation*}
t=\frac{x^{\beta }}{\alpha ^{\beta }}
\end{equation*}%
we have,

$$dt= \frac{\beta }{\alpha ^{\beta }}x^{\beta -1}dx \Leftrightarrow x^{\beta -1}dx=\frac{\alpha ^{\beta }}{\beta }dt.$$

Then,%
\begin{equation*}
M_{2}=-2\left( \frac{\beta }{\alpha }\right) ^{\tau
+1}\int_{0}^{\infty }\frac{t^{\frac{\left( \beta -1\right) \tau }{%
\beta }}}{\left( t+1\right) ^{2\tau +3}}dt.
\end{equation*}

By Lemma \ref{lema1} with $m(\beta )=\frac{\left( \beta -1\right) \tau }{%
\beta }$ and $s=2\tau +3,$ we have,%
\begin{equation*}
\int_{0}^{\infty }\frac{t^{\frac{\left( \beta -1\right) \tau }{\beta
}}}{\left( t+1\right) ^{2\tau +3}}dt=B\left( \frac{\tau \beta +2\beta +\tau
}{\beta },\frac{\tau \beta +\beta -\tau }{\beta }\right) .
\end{equation*}%

Therefore,
\begin{eqnarray*}
M_{2} &=&-2\left( \frac{\beta }{\alpha }\right) ^{\tau +1}B\left( \frac{\tau
\beta +2\beta +\tau }{\beta },\frac{\tau \beta +\beta -\tau }{\beta }\right)
\\
&=&2\left( \frac{\beta }{\alpha }\right) ^{\tau +1}\frac{\tau \beta +\beta
+\tau }{\beta \left( 2\tau +2\right) }B\left( \frac{\tau \beta +\beta +\tau
}{\beta },\frac{\tau \beta +\beta -\tau }{\beta }\right) .
\end{eqnarray*}

Finally,
\begin{eqnarray*}
\xi _{\tau }\left( \alpha \right) &=&\left( \frac{\beta }{\alpha }\right)
^{\tau +1}B\left( \frac{\tau \beta +\beta +\tau }{\beta },\frac{\tau \beta
+\beta -\tau }{\beta }\right) \left( 1-2\frac{\tau \beta +\beta +\tau }{%
\beta \left( 2\tau +2\right) }\right) \\
&=&\left( \frac{\beta }{\alpha }\right) ^{\tau +1}B\left( \frac{\tau \beta
+\beta +\tau }{\beta },\frac{\tau \beta +\beta -\tau }{\beta }\right) \left(
\frac{-\tau }{\beta +\tau \beta }\right) .
\end{eqnarray*}

\subsection{Proof Theorem \protect\ref{Theorem2}}

We have,
\begin{equation*}
\log f_{\beta }(x)=\log \beta +\beta \log \alpha +\left( \beta -1\right)
\log \frac{x}{\alpha }+\left( \beta -1\right) \log \alpha -2\log \left(
1+\left( \frac{x^{\beta }}{\alpha ^{\beta }}\right) \right) -2\beta \log
\alpha ,
\end{equation*}%
and

\begin{eqnarray*}
\frac{\partial \log f_{\beta }(x)}{\partial \beta } &=&\frac{1}{\beta }+\log
\alpha +\log \frac{x}{\alpha }+\log \alpha -2\frac{\left( \frac{x}{\alpha }%
\right) ^{\beta }\log \frac{x}{\alpha }}{1+\left( \frac{x}{\alpha }\right)
^{\beta }}-2\log \alpha \\
&=&\frac{1}{\beta }+\log \frac{x}{\alpha }-2\frac{\left( \frac{x}{\alpha }%
\right) ^{\beta }\log \frac{x}{\alpha }}{1+\left( \frac{x}{\alpha }\right)
^{\beta }}.
\end{eqnarray*}%

Therefore,%
\begin{eqnarray*}
\left( \frac{\partial \log f_{\beta }(x)}{\partial \beta }\right) ^{2} &=&%
\frac{1}{\beta ^{2}}+\frac{2}{\beta }\log \frac{x}{\alpha }-\frac{4}{\beta }%
\frac{\left( \frac{x}{\alpha }\right) ^{\beta }\log \frac{x}{\alpha }}{%
1+\left( \frac{x}{\alpha }\right) ^{\beta }}+\left( \log \frac{x}{\alpha }%
\right) ^{2} \\
&&-4\left( \frac{x}{\alpha }\right) ^{\beta }\frac{\left( \log \frac{x}{%
\alpha }\right) ^{2}}{1+\left( \frac{x}{\alpha }\right) ^{\beta }}+4\frac{%
\left( \frac{x}{\alpha }\right) ^{2\beta }\left( \log \frac{x}{\alpha }%
\right) ^{2}}{\left( 1+\left( \frac{x}{\alpha }\right) ^{\beta }\right) ^{2}}%
,
\end{eqnarray*}%
and \ $J_{\tau }(\beta )$ is given by
\begin{eqnarray*}
J_{\tau }(\beta ) &=&\frac{1}{\beta ^{2}}\int_{0}^{\infty }\frac{%
\beta ^{\tau +1}\alpha ^{\beta \left( \tau +1\right) }x^{\left( \beta
-1\right) \left( \tau +1\right) }}{\left( x^{\beta }+\alpha ^{\beta }\right)
^{^{2\left( \tau +1\right) }}}dx+\frac{2}{\beta }\int_{0}^{\infty
}\log \frac{x}{\alpha }\frac{\beta ^{\tau +1}\alpha ^{\beta \left( \tau
+1\right) }x^{\left( \beta -1\right) \left( \tau +1\right) }}{\left(
x^{\beta }+\alpha ^{\beta }\right) ^{^{2\left( \tau +1\right) }}}dx \\
&&-\frac{4}{\beta }\int_{0}^{\infty }\frac{\left( \frac{x}{\alpha }%
\right) ^{\beta }\log \frac{x}{\alpha }}{1+\left( \frac{x}{\alpha }\right)
^{\beta }}\frac{\beta ^{\tau +1}\alpha ^{\beta \left( \tau +1\right)
}x^{\left( \beta -1\right) \left( \tau +1\right) }}{\left( x^{\beta }+\alpha
^{\beta }\right) ^{^{2\left( \tau +1\right) }}}dx \\
&&+\int_{0}^{\infty }\left( \log \frac{x}{\alpha }\right) ^{2}\frac{%
\beta ^{\tau +1}\alpha ^{\beta \left( \tau +1\right) }x^{\left( \beta
-1\right) \left( \tau +1\right) }}{\left( x^{\beta }+\alpha ^{\beta }\right)
^{^{2\left( \tau +1\right) }}}dx \\
&&-4\int_{0}^{\infty }\left( \frac{x}{\alpha }\right) ^{\beta }\frac{%
\left( \log \frac{x}{\alpha }\right) ^{2}}{1+\left( \frac{x}{\alpha }\right)
^{\beta }}\frac{\beta ^{\tau +1}\alpha ^{\beta \left( \tau +1\right)
}x^{\left( \beta -1\right) \left( \tau +1\right) }}{\left( x^{\beta }+\alpha
^{\beta }\right) ^{^{2\left( \tau +1\right) }}}dx \\
&&+4\int_{0}^{\infty }\frac{\left( \frac{x}{\alpha }\right) ^{2\beta
}\left( \log \frac{x}{\alpha }\right) ^{2}}{\left( 1+\left( \frac{x}{\alpha }%
\right) ^{\beta }\right) ^{2}}\frac{\beta ^{\tau +1}\alpha ^{\beta \left(
\tau +1\right) }x^{\left( \beta -1\right) \left( \tau +1\right) }}{\left(
x^{\beta }+\alpha ^{\beta }\right) ^{^{2\left( \tau +1\right) }}}dx \\
&=&N_{1}+N_{2}+N_{3}+N_{4}+N_{5}+N_{6}.
\end{eqnarray*}%

\underline{In relation to $N_{1}$ we have }%
\begin{eqnarray}
N_{1} &=&\frac{1}{\beta ^{2}}\int_{0}^{\infty }\frac{\beta ^{\tau
+1}\alpha ^{\beta \left( \tau +1\right) }x^{\left( \beta -1\right) \left(
\tau +1\right) }}{\left( x^{\beta }+\alpha ^{\beta }\right) ^{^{2\left( \tau
+1\right) }}}dx  \notag \\
&=&\frac{1}{\beta ^{2}}\beta ^{\tau +1}\alpha ^{\beta \left( \tau +1\right)
}\alpha ^{-2\beta \left( \tau +1\right) }\int_{0}^{\infty }\frac{%
\left( x^{\beta }\right) ^{\frac{\left( \tau +1\right) \left( \beta
-1\right) }{\beta }}}{\left( \frac{x^{\beta }}{\alpha ^{\beta }}+1\right)
^{2\left( \tau +1\right) }}dx  \notag \\
&=&\frac{1}{\beta ^{2}}\beta ^{\tau +1}\alpha ^{\beta \left( \tau +1\right)
}\alpha ^{-2\beta \left( \tau +1\right) }\alpha ^{\tau \left( \beta
-1\right) }\frac{1}{\beta }\alpha ^{\beta }\int_{0}^{\infty }\frac{%
\left( \frac{x^{\beta }}{\alpha ^{\beta }}\right) ^{\frac{\tau \left( \beta
-1\right) }{\beta }}x^{\beta -1}}{\left( \frac{x^{\beta }}{\alpha ^{\beta }}%
+1\right) ^{2\left( \tau +1\right) }}dx  \notag \\
&=&\dfrac{\beta ^{\tau -2}}{\alpha ^{\tau }}\int_{0}^{\infty }\frac{%
t^{\frac{\left( \beta -1\right) \tau }{\beta }}}{\left( t+1\right) ^{2\tau
+2}}dt.  \label{N1}
\end{eqnarray}

By Lemma \ref{lema1} with $m(\beta )=\frac{\left( \beta -1\right) \tau }{%
\beta }$ and $s=2\tau +2,$ we have,%
\begin{equation*}
\int_{0}^{\infty }\frac{t^{\frac{\left( \beta -1\right) \tau }{\beta
}}}{\left( t+1\right) ^{2\tau +2}}dt=B\left( \frac{\tau \beta +\tau +\beta }{%
\beta },\frac{\tau \beta -\tau +\beta }{\beta }\right) .
\end{equation*}

Therefore,
\begin{equation*}
N_{1}=\dfrac{\beta ^{\tau -2}}{\alpha ^{\tau }}B\left( \frac{\tau \beta
+\tau +\beta }{\beta },\frac{\tau \beta -\tau +\beta }{\beta }\right) .
\end{equation*}%

\underline{In relation to $N_{2}$ we have }%
\begin{eqnarray}
N_{2} &=&\frac{2}{\beta }\int_{0}^{\infty }\frac{1}{\beta }\log
\frac{x^{\beta }}{\alpha ^{\beta }}\frac{\beta ^{\tau +1}\alpha ^{\beta
\left( \tau +1\right) }x^{\left( \beta -1\right) \left( \tau +1\right) }}{%
\left( x^{\beta }+\alpha ^{\beta }\right) ^{^{2\left( \tau +1\right) }}}dx
\notag \\
&=&\frac{2}{\beta ^{2}}\beta ^{\tau +1}\alpha ^{\beta \left( \tau +1\right)
}\alpha ^{\left( \beta -1\right) \tau }\alpha ^{-2\beta \left( \tau
+1\right) }\int_{0}^{\infty }\log \frac{x^{\beta }}{\alpha ^{\beta }}%
\frac{\left( \frac{x^{\beta }}{\alpha ^{\beta }}\right) ^{\frac{\tau \left(
\beta -1\right) }{\beta }}}{\left( \frac{x^{\beta }}{\alpha ^{\beta }}%
+1\right) ^{2\left( \tau +1\right) }}x^{\beta -1}dx  \notag \\
&=&2\frac{\beta ^{\tau -2}}{\alpha ^{\tau }}\int_{0}^{\infty }\log t%
\frac{t^{\frac{\left( \beta -1\right) \tau }{\beta }}}{\left( t+1\right)
^{2\tau +2}}dt.  \label{N2}
\end{eqnarray}%

By Lemma \ref{lema2} with $m(\beta )=\frac{\left( \beta -1\right) \tau }{%
\beta }$ and $s=2\tau +2,$ we have,%
\begin{equation*}
\int_{0}^{\infty }\log t\frac{t^{\frac{\left( \beta -1\right) \tau }{%
\beta }}}{\left( t+1\right) ^{2\tau +2}}dt=B(\frac{\tau \beta +\tau +\beta }{%
\beta },\frac{\tau \beta -\tau +\beta }{\beta })\left\{ \Psi \left( \frac{%
\tau \beta -\tau +\beta }{\beta }\right) -\Psi \left( \frac{\tau \beta +\tau
+\beta }{\beta }\right) \right\} .
\end{equation*}%

Therefore,%
\begin{equation*}
N_{2}=2\frac{\beta ^{\tau -2}}{\alpha ^{\tau }}B(\frac{\tau \beta +\tau
+\beta }{\beta },\frac{\tau \beta -\tau +\beta }{\beta })\left\{ \Psi \left(
\frac{\tau \beta -\tau +\beta }{\beta }\right) -\Psi \left( \frac{\tau \beta
+\tau +\beta }{\beta }\right) \right\} .
\end{equation*}%

\underline{In relation to $N_{3}$ we have }%
\begin{eqnarray*}
N_{3} &=&-\frac{4}{\beta }\int_{0}^{\infty }\frac{\left( \frac{x}{%
\alpha }\right) ^{\beta }\log \frac{x}{\alpha }}{1+\left( \frac{x}{\alpha }%
\right) ^{\beta }}\frac{\beta ^{\tau +1}\alpha ^{\beta \left( \tau +1\right)
}x^{\left( \beta -1\right) \left( \tau +1\right) }}{\left( x^{\beta }+\alpha
^{\beta }\right) ^{^{2\left( \tau +1\right) }}}dx \\
&=&-4\frac{\beta ^{\tau -2}}{\alpha ^{\tau }}\int_{0}^{\infty }\log t%
\frac{t^{\frac{\beta \tau -\tau +\beta }{\beta }}}{\left( t+1\right) ^{2\tau
+3}}dt.
\end{eqnarray*}%

By Lemma \ref{lema2} with $m(\beta )=\frac{\beta \tau -\tau +\beta }{\beta }$
and $s=2\tau +3,$ we have,%
\begin{equation*}
\int_{0}^{\infty }\log t\frac{t^{\frac{\beta \tau -\tau +\beta }{%
\beta }}}{\left( t+1\right) ^{2\tau +3}}dt=B\left( \frac{\tau \beta +\beta
+\tau }{\beta },\frac{\tau \beta +2\beta -\tau }{\beta }\right) \left\{ \Psi
\left( \frac{\tau \beta +2\beta -\tau }{\beta }\right) -\Psi \left( \frac{%
\tau \beta +\beta +\tau }{\beta }\right) \right\} ,
\end{equation*}%
and
\begin{equation}
N_{3}=-4\frac{\beta ^{\tau -2}}{\alpha ^{\tau }}B\left( \frac{\tau \beta
+\beta +\tau }{\beta },\frac{\tau \beta +2\beta -\tau }{\beta }\right)
\left\{ \Psi \left( \frac{\tau \beta +2\beta -\tau }{\beta }\right) -\Psi
\left( \frac{\tau \beta +\beta +\tau }{\beta }\right) \right\} .  \label{N3}
\end{equation}

\underline{In relation to $N_{4}$ we have }%
\begin{eqnarray*}
N_{4} &=&\int_{0}^{\infty }\left( \log \frac{x}{\alpha }\right) ^{2}%
\frac{\beta ^{\tau +1}\alpha ^{\beta \left( \tau +1\right) }x^{\left( \beta
-1\right) \left( \tau +1\right) }}{\left( x^{\beta }+\alpha ^{\beta }\right)
^{^{2\left( \tau +1\right) }}}dx \\
&=&\frac{\beta ^{\tau -2}}{\alpha ^{\tau }}\int_{0}^{\infty }\left(
\log t\right) ^{2}\frac{t^{\frac{\beta \tau -\tau }{\beta }}}{\left(
t+1\right) ^{2\tau +2}}dt.
\end{eqnarray*}

By Lemma \ref{lema3} with $m(\beta )=\frac{\beta \tau -\tau }{\beta }$ and $%
s=2\tau +2,$ we have,%
\begin{eqnarray*}
\int_{0}^{\infty }\left( \log t\right) ^{2}\frac{t^{\frac{\beta \tau
-\tau +\beta }{\beta }}}{\left( t+1\right) ^{2\tau +2}}dt &=&B\left( \frac{%
\tau \beta +\beta +\tau }{\beta },\frac{\tau \beta +\beta -\tau }{\beta }%
\right) \left\{ \left[ \Psi \left( \frac{\tau \beta +\beta -\tau }{\beta }%
\right) -\Psi \left( \frac{\tau \beta +\beta +\tau }{\beta }\right) \right]
^{2}\right. \\
&&\left. \times \Psi ^{\prime }\left( \frac{\tau \beta +\beta -\tau }{\beta }%
\right) +\Psi ^{\prime }\left( \frac{\tau \beta +\beta +\tau }{\beta }%
\right) \right\} ,
\end{eqnarray*}%
and
\begin{eqnarray}
N_{4} &=&\frac{\beta ^{\tau -2}}{\alpha ^{\tau }}B\left( \frac{\tau \beta
+\beta +\tau }{\beta },\frac{\tau \beta +\beta -\tau }{\beta }\right)
\left\{ \left[ \Psi \left( \frac{\tau \beta +\beta -\tau }{\beta }\right)
-\Psi \left( \frac{\tau \beta +\beta +\tau }{\beta }\right) \right]
^{2}\right.  \label{N4} \\
&&\left. +\Psi ^{\prime }\left( \frac{\tau \beta +\beta -\tau }{\beta }%
\right) +\Psi ^{\prime }\left( \frac{\tau \beta +\beta +\tau }{\beta }%
\right) \right\} .  \notag
\end{eqnarray}

\underline{In relation to $N_{5}$ we have }

\begin{eqnarray*}
N_{5} &=&-4\int_{0}^{\infty }\left( \log \frac{x}{\alpha }\right)
^{2}\left( \frac{x}{\alpha }\right) ^{\beta }\frac{\beta ^{\tau +1}\alpha
^{\beta \left( \tau +1\right) }x^{\left( \beta -1\right) \left( \tau
+1\right) }}{\left( x^{\beta }+\alpha ^{\beta }\right) ^{^{2\left( \tau
+1\right) }}\left( 1+\left( \frac{x}{\alpha }\right) ^{\beta }\right) }dx \\
&=&-4\frac{\beta ^{\tau -2}}{\alpha ^{\tau }}\int_{0}^{\infty
}\left( \log t\right) ^{2}\frac{t^{\frac{\beta \tau -\tau +\beta }{\beta }}}{%
\left( t+1\right) ^{2\tau +3}}dt.
\end{eqnarray*}%

By Lemma \ref{lema3} with $m(\beta )=\frac{\beta \tau -\tau +\beta }{\beta }$
and $s=2\tau +3,$ we have,%
\begin{eqnarray}
N_{5} &=&-4\frac{\beta ^{\tau -2}}{\alpha ^{\tau }}B\left( \frac{\tau \beta
+\beta +\tau }{\beta },\frac{\tau \beta +2\beta -\tau }{\beta }\right)
\left\{ \left( \Psi \left( \frac{\tau \beta +2\beta -\tau }{\beta }\right)
-\Psi \left( \frac{\tau \beta +\beta +\tau }{\beta }\right) \right)
^{2}\right.  \label{N5} \\
&&+\left( \Psi ^{\prime }\left( \frac{\tau \beta +2\beta -\tau }{\beta }%
\right) +\Psi ^{\prime }\left( \frac{\tau \beta +\beta +\tau }{\beta }%
\right) \right) .  \notag
\end{eqnarray}%

\underline{In relation to $N_{6}$ we have}%
\begin{eqnarray*}
N_{6} &=&4\int_{0}^{\infty }\frac{\left( \frac{x}{\alpha }\right)
^{2\beta }\left( \log \frac{x}{\alpha }\right) ^{2}}{\left( 1+\left( \frac{x%
}{\alpha }\right) ^{\beta }\right) ^{2}}\frac{\beta ^{\tau +1}\alpha ^{\beta
\left( \tau +1\right) }x^{\left( \beta -1\right) \left( \tau +1\right) }}{%
\left( x^{\beta }+\alpha ^{\beta }\right) ^{^{2\left( \tau +1\right) }}}dx \\
&=&4\frac{\beta ^{\tau -2}}{\alpha ^{\tau }}\int_{0}^{\infty }\left(
\log t\right) ^{2}\frac{t^{\frac{\beta \tau -\tau +2\beta }{\beta }}}{\left(
t+1\right) ^{2\tau +4}}dt.
\end{eqnarray*}%

By Lemma \ref{lema3} with $m(\beta )=\frac{\beta \tau -\tau +2\beta }{\beta }
$ and $s=2\tau +4,$ we have,%
\begin{eqnarray*}
\int_{0}^{\infty }\left( \log t\right) ^{2}\frac{t^{\frac{\beta \tau
-\tau +2\beta }{\beta }}}{\left( t+1\right) ^{2\tau +4}}dt &=&B\left( \frac{%
\tau \beta +\beta +\tau }{\beta },\frac{\tau \beta +3\beta -\tau }{\beta }%
\right) \\
& & \left\{ \left( \Psi \left( \frac{\tau \beta +3\beta -\tau }{\beta }%
\right) -\Psi \left( \frac{\tau \beta +\beta +\tau }{\beta }\right) \right)
^{2}\right. \\
&&+\left( \Psi ^{\prime }\left( \frac{\tau \beta +3\beta -\tau }{\beta }%
\right) +\Psi ^{\prime }\left( \frac{\tau \beta +\beta +\tau }{\beta }%
\right) \right\} .
\end{eqnarray*}%

Therefore,
\begin{eqnarray}
N_{6} &=&4\frac{\beta ^{\tau -2}}{\alpha ^{\tau }}B\left( \frac{\tau \beta
+\beta +\tau }{\beta },\frac{\tau \beta +3\beta -\tau }{\beta }\right)
\left\{ \left[ \Psi \left( \frac{\tau \beta +3\beta -\tau }{\beta }\right)
-\Psi \left( \frac{\tau \beta +\beta +\tau }{\beta }\right) \right]
^{2}\right.  \label{N6} \\
&&+\left( \Psi ^{\prime }\left( \frac{\tau \beta +3\beta -\tau }{\beta }%
\right) +\Psi ^{\prime }\left( \frac{\tau \beta +\beta +\tau }{\beta }%
\right) \right) .  \notag
\end{eqnarray}

For $\tau =0,$ we obtain

\begin{eqnarray*}
I_{F}\left( \beta \right) &=&J_{\tau =0}(\beta
)=N_{1}+N_{2}+N_{3}+N_{4}+N_{5}+N_{6} \\
&=&\frac{1}{\beta ^{2}}\left( 1-2-\frac{2\pi ^{2}}{3}\right) +\frac{4}{3}+8%
\frac{\pi ^{2}}{18}+\frac{2\pi ^{2}}{6} \\
&=&\frac{1}{\beta ^{2}}\left( \frac{3}{9}+8\frac{\pi ^{2}}{18}-\frac{2\pi
^{2}}{6}\right) \\
& = & \frac{3+\pi ^{2}}{9 \beta ^{2}}
\end{eqnarray*}

\subsection{Proof Theorem \protect\ref{Theorem3}}

We have,%
\begin{eqnarray*}
\xi _{\tau }(\beta ) &=&\int_{0}^{\infty } \frac{\partial \log
f_{\beta }(x)}{\partial \beta } f_{\beta }(x)^{\tau +1}dx \\
&=&\int_{0}^{\infty }\left( \frac{1}{\beta }+\log \frac{x}{\alpha }-2%
\frac{\left( \frac{x}{\alpha }\right) ^{\beta }\log \frac{x}{\alpha }}{%
1+\left( \frac{x}{\alpha }\right) ^{\beta }}\right) f_{\beta }(x)^{\tau +1}dx
\\
&=&\frac{1}{\beta }\int_{0}^{\infty }\frac{\beta ^{\tau +1}\alpha
^{\left( \tau +1\right) \beta }x^{\left( \tau +1\right) \left( \beta
-1\right) }}{\left( x^{\beta }+\alpha ^{\beta }\right) ^{2\left( \tau
+1\right) }}dx \\
&&+\int_{0}^{\infty }\left( \log \frac{x}{\alpha }\right) \frac{\beta ^{\tau
+1}\alpha ^{\left( \tau +1\right) \beta }x^{\left( \tau +1\right) \left(
\beta -1\right) }}{\left( x^{\beta }+\alpha ^{\beta }\right) ^{2\left( \tau
+1\right) }}dx \\
&&-2\int_{0}^{\infty }\frac{\left( \frac{x}{\alpha }\right) ^{\beta
}\left( \log \frac{x}{\alpha }\right) }{1+\left( \frac{x}{\alpha }\right) ^{\beta }}\frac{%
\beta ^{\tau +1}\alpha ^{\left( \tau +1\right) \beta }x^{\left( \tau
+1\right) \left( \beta -1\right) }}{\left( x^{\beta }+\alpha ^{\beta
}\right) ^{2\left( \tau +1\right) }}dx \\
&=&S_{1}+S_{2}+S_{3}.
\end{eqnarray*}

\underline{In relation to $S_{1}$ we have }%
\begin{eqnarray*}
S_{1} &=&\frac{1}{\beta }\int_{0}^{\infty }\frac{\beta ^{\tau
+1}\alpha ^{\left( \tau +1\right) \beta }x^{\left( \tau +1\right) \left(
\beta -1\right) }}{\left( x^{\beta }+\alpha ^{\beta }\right) ^{2\left( \tau
+1\right) }}dx \\
&=&\frac{\beta ^{\tau -1}}{\alpha ^{\tau }}\int_{0}^{\infty }\frac{%
t^{\frac{\beta \tau -\tau }{\beta }}}{\left( t+1\right) ^{2\tau +2}}dt.
\end{eqnarray*}

By Lemma \ref{lema1} with $m(\beta )=\frac{\beta \tau -\tau }{\beta }$ and $%
s=2\tau +2,$ we have,%
\begin{equation*}
\int_{0}^{\infty }\frac{t^{\frac{\beta \tau -\tau }{\beta }}}{\left(
t+1\right) ^{2\tau +2}}dt=B\left( \frac{\tau \beta +\beta +\tau }{\beta },%
\frac{\tau \beta +\beta -\tau }{\beta }\right) .
\end{equation*}%

Therefore,%
\begin{equation*}
S_{1}=\frac{\beta ^{\tau -1}}{\alpha ^{\tau }}B\left( \frac{\tau \beta
+\beta +\tau }{\beta },\frac{\tau \beta +\beta -\tau }{\beta }\right) .
\end{equation*}

\underline{In relation to $S_{2}$ we have }%
\begin{eqnarray*}
S_{2} &=&\int_{0}^{\infty }\log \frac{x}{\alpha }\frac{\beta ^{\tau
+1}\alpha ^{\left( \tau +1\right) \beta }x^{\left( \tau +1\right) \left(
\beta -1\right) }}{\left( x^{\beta }+\alpha ^{\beta }\right) ^{2\left( \tau
+1\right) }}dx \\
&=&\frac{\beta ^{\tau -1}}{\alpha ^{\tau }}\int_{0}^{\infty }\log t%
\frac{t^{\frac{\beta \tau -\tau }{\beta }}}{\left( t+1\right) ^{2\tau +2}}dt.
\end{eqnarray*}%

By Lemma \ref{lema2} with $m(\beta )=\frac{\beta \tau -\tau }{\beta }$ and $%
s=2\tau +2,$ we have,%
\begin{equation*}
\int_{0}^{\infty }\log t\frac{t^{\frac{\beta \tau -\tau }{\beta }}}{%
\left( t+1\right) ^{2\tau +2}}dt=B\left( \frac{\tau \beta +\beta +\tau }{%
\beta },\frac{\tau \beta +\beta -\tau }{\beta }\right) \left\{ \Psi \left(
\frac{\tau \beta +\beta -\tau }{\beta }\right) -\Psi \left( \frac{\tau \beta
+\beta +\tau }{\beta }\right) \right\} ,
\end{equation*}%
and
\begin{equation*}
S_{2}=\frac{\beta ^{\tau -1}}{\alpha ^{\tau }}B\left( \frac{\tau \beta
+\beta +\tau }{\beta },\frac{\tau \beta +\beta -\tau }{\beta }\right)
\left\{ \Psi \left( \frac{\tau \beta +\beta -\tau }{\beta }\right) -\Psi
\left( \frac{\tau \beta +\beta +\tau }{\beta }\right) \right\} .
\end{equation*}

\underline{In relation to $S_{3}$ we have }%
\begin{eqnarray*}
S_{3} &=&-2\int_{0}^{\infty }\frac{\left( \frac{x}{\alpha }\right)
^{\beta }\log \frac{x}{\alpha }}{1+\left( \frac{x}{\alpha }\right) ^{\beta }}%
\frac{\beta ^{\tau +1}\alpha ^{\left( \tau +1\right) \beta }x^{\left( \tau
+1\right) \left( \beta -1\right) }}{\left( x^{\beta }+\alpha ^{\beta
}\right) ^{2\left( \tau +1\right) }}dx \\
&=&-2\frac{\beta ^{\tau -1}}{\alpha ^{\tau }}\int_{0}^{\infty }\log t%
\frac{t^{\frac{\tau \left( \beta -1\right) +\beta }{\beta }}}{\left(
t+1\right) ^{2\tau +3}}dt.
\end{eqnarray*}

By Lemma \ref{lema2} with $m(\beta )=\frac{\tau \left( \beta -1\right)
+\beta }{\beta }$ and $s=2\tau +3,$ we have,%
\begin{equation*}
\int_{0}^{\infty }\log t\frac{t^{\frac{\tau \left( \beta -1\right)
+\beta }{\beta }}}{\left( t+1\right) ^{2\tau +3}}dt=B\left( \frac{\tau \beta
+\beta +\tau }{\beta },\frac{\tau \beta +2\beta -\tau }{\beta }\right)
\left\{ \Psi \left( \frac{\tau \beta +2\beta -\tau }{\beta }\right) -\Psi
\left( \frac{\tau \beta +\beta +\tau }{\beta }\right) \right\} ,
\end{equation*}%
and
\begin{equation*}
S_{3}=-2\frac{\beta ^{\tau -1}}{\alpha ^{\tau }}B\left( \frac{\tau \beta
+\beta +\tau }{\beta },\frac{\tau \beta +2\beta -\tau }{\beta }\right)
\left\{ \Psi \left( \frac{\tau \beta +2\beta -\tau }{\beta }\right) -\Psi
\left( \frac{\tau \beta +\beta +\tau }{\beta }\right) \right\} .
\end{equation*}%

But%
\begin{equation*}
B\left( \frac{\tau \beta +\beta +\tau }{\beta },\frac{\tau \beta +2\beta
-\tau }{\beta }\right) =\frac{\tau \beta +\beta -\tau }{\beta \left( 2\tau
+2\right) }B\left( \frac{\tau \beta +\beta +\tau }{\beta },\frac{\tau \beta
+\beta -\tau }{\beta }\right) .
\end{equation*}%
and
\begin{equation*}
\Psi \left( \frac{\tau \beta +2\beta
-\tau }{\beta }\right) = \Psi \left( \frac{\tau \beta +\beta
-\tau }{\beta }\right) + {\beta \over \tau \beta +\beta
-\tau }
\end{equation*}%

Hence,

\begin{eqnarray*}
S_3 & = & - \frac{\beta ^{\tau -1}}{\alpha ^{\tau }} B\left( \frac{\tau \beta
+\beta +\tau }{\beta },\frac{\tau \beta +\beta -\tau }{\beta }\right) \\
& & \left[ {1\over \tau +1} + \frac{\tau \beta +\beta -\tau }{\beta \left( \tau
+1\right) }\left( \Psi \left( \frac{\tau \beta +\beta -\tau }{\beta }\right) -\Psi
\left( \frac{\tau \beta +\beta +\tau }{\beta }\right) \right) \right]
\end{eqnarray*}

Therefore,%
\begin{eqnarray*}
\xi _{\tau }(\beta ) &=&S_{1}+S_{2}+S_{3.} \\
&=&\frac{\beta ^{\tau -1}}{\alpha ^{\tau }}{\tau \over \tau +1} B\left( \frac{\tau \beta
+\beta +\tau }{\beta },\frac{\tau \beta +\beta -\tau }{\beta }\right) \\
& & \left\{ 1+\Psi \left( \frac{\tau \beta +\beta -\tau }{\beta }\right) -\Psi
\left( \frac{\tau \beta +\beta +\tau }{\beta }\right) \right\}
\end{eqnarray*}%

\subsection{Proof Theorem \protect\ref{Theorem4}}

We have,
\begin{equation*}
\text{ }J_{\tau }^{12}\left( \alpha ,\beta \right) =\int_{0}^{\infty
}\left( \frac{\partial \log f_{\alpha ,\beta }(x)}{\partial \beta }\right)
\left( \frac{\partial \log f_{\alpha ,\beta }(x)}{\partial \beta }\right)
f_{\beta }(x)^{\tau +1}dx.
\end{equation*}%

But
\begin{equation*}
\frac{\partial \log f_{\alpha }(x)}{\partial \alpha }=\frac{\beta }{\alpha }%
-2\frac{\beta \alpha ^{\beta -1}}{\left( x^{\beta }+\alpha ^{\beta }\right) }%
\text{ and }\frac{\partial \log f_{\beta }(x)}{\partial \beta }=\frac{1}{%
\beta }+\log \frac{x}{\alpha }-2\frac{\left( \frac{x}{\alpha }\right)
^{\beta }\log \frac{x}{\alpha }}{1+\left( \frac{x}{\alpha }\right) ^{\beta }}%
.
\end{equation*}%

Therefore,%
\begin{eqnarray*}
\text{ }J_{\tau }^{12}\left( \alpha ,\beta \right) &=&\frac{1}{\alpha }%
\int_{0}^{\infty }\frac{\beta ^{\tau +1}\alpha ^{\left( \tau
+1\right) \beta }x^{\left( \tau +1\right) \left( \beta -1\right) }}{\left(
x^{\beta }+\alpha ^{\beta }\right) ^{2\left( \tau +1\right) }}dx \\
&&+\frac{\beta }{\alpha }\int_{0}^{\infty }\log \frac{x}{\alpha }%
\frac{\beta ^{\tau +1}\alpha ^{\left( \tau +1\right) \beta }x^{\left( \tau
+1\right) \left( \beta -1\right) }}{\left( x^{\beta }+\alpha ^{\beta
}\right) ^{2\left( \tau +1\right) }}dx \\
&&-2\frac{\beta }{\alpha }\int_{0}^{\infty }\frac{\left( \frac{x}{%
\alpha }\right) ^{\beta }\log \frac{x}{\alpha }}{1+\left( \frac{x}{\alpha }%
\right) ^{\beta }}\frac{\beta ^{\tau +1}\alpha ^{\left( \tau +1\right) \beta
}x^{\left( \tau +1\right) \left( \beta -1\right) }}{\left( x^{\beta }+\alpha
^{\beta }\right) ^{2\left( \tau +1\right) }}dx \\
&&-2\alpha ^{\beta -1}\int_{0}^{\infty }\frac{1}{\left( x^{\beta
}+\alpha ^{\beta }\right) }\frac{\beta ^{\tau +1}\alpha ^{\left( \tau
+1\right) \beta }x^{\left( \tau +1\right) \left( \beta -1\right) }}{\left(
x^{\beta }+\alpha ^{\beta }\right) ^{2\left( \tau +1\right) }}dx \\
&&-2\beta \alpha ^{\beta -1}\int_{0}^{\infty }\log \frac{x}{\alpha }%
\frac{1}{\left( x^{\beta }+\alpha ^{\beta }\right) }\frac{\beta ^{\tau
+1}\alpha ^{\left( \tau +1\right) \beta }x^{\left( \tau +1\right) \left(
\beta -1\right) }}{\left( x^{\beta }+\alpha ^{\beta }\right) ^{2\left( \tau
+1\right) }}dx \\
&&+4\beta \alpha ^{\beta -1}\int_{0}^{\infty }\log \frac{x}{\alpha }%
\frac{\left( \frac{x}{\alpha }\right) ^{\beta }}{\left( x^{\beta }+\alpha
^{\beta }\right) }\frac{1}{1+\left( \frac{x}{\alpha }\right) ^{\beta }}\frac{%
\beta ^{\tau +1}\alpha ^{\left( \tau +1\right) \beta }x^{\left( \tau
+1\right) \left( \beta -1\right) }}{\left( x^{\beta }+\alpha ^{\beta
}\right) ^{2\left( \tau +1\right) }}dx \\
&=&B_{1}+B_{2}+B_{3}+B_{4}+B_{5}+B_{6}.
\end{eqnarray*}%

\underline{In relation to $B_{1}$ we have }%
\begin{eqnarray*}
B_{1} &=&\frac{1}{\alpha }\int_{0}^{\infty }\frac{\beta ^{\tau
+1}\alpha ^{\left( \tau +1\right) \beta }x^{\left( \tau +1\right) \left(
\beta -1\right) }}{\left( x^{\beta }+\alpha ^{\beta }\right) ^{2\left( \tau
+1\right) }}dx \\
&=&\frac{\beta ^{\tau }}{\alpha ^{\tau +1}}\int_{0}^{\infty }\frac{%
t^{\frac{\tau \left( \beta -1\right) }{\beta }}}{\left( t+1\right) ^{2\tau
+2}}dt.
\end{eqnarray*}%

By Lemma \ref{lema1} with $m(\beta )=\frac{\tau \left( \beta -1\right) }{%
\beta }$ and $s=2\tau +2,$ we have,%
\begin{equation}
B_{1}=\frac{\beta ^{\tau }}{\alpha ^{\tau +1}}B\left( \frac{\tau \beta
+\beta +\tau }{\beta },\frac{\tau \beta +\beta -\tau }{\beta }\right) .
\label{B1}
\end{equation}%

\underline{In relation to $B_{2}$ we have }%
\begin{eqnarray*}
B_{2} &=&\frac{\beta }{\alpha }\int_{0}^{\infty }\log \frac{x}{%
\alpha }\frac{\beta ^{\tau +1}\alpha ^{\left( \tau +1\right) \beta
}x^{\left( \tau +1\right) \left( \beta -1\right) }}{\left( x^{\beta }+\alpha
^{\beta }\right) ^{2\left( \tau +1\right) }}dx \\
&=&\frac{\beta ^{\tau }}{\alpha ^{\tau +1}}\int_{0}^{\infty }\log t%
\frac{t^{\frac{\tau \left( \beta -1\right) }{\beta }}}{\left( t+1\right)
^{2\tau +2}}dt.
\end{eqnarray*}%

By Lemma \ref{lema2} with $m(\beta )=\frac{\tau \left( \beta -1\right) }{%
\beta }$ and $s=2\tau +2,$ we have,%
\begin{equation}
B_{2}=\frac{\beta ^{\tau }}{\alpha ^{\tau +1}}B\left( \frac{\tau \beta
+\beta +\tau }{\beta },\frac{\tau \beta +\beta -\tau }{\beta }\right)
\left\{ \Psi \left( \frac{\tau \beta +\beta -\tau }{\beta }\right) -\Psi
\left( \frac{\tau \beta +\beta +\tau }{\beta }\right) \right\} .  \label{B2}
\end{equation}%

\underline{In relation to $B_{3}$ we have }%
\begin{eqnarray*}
B_{3} &=&-2\frac{\beta }{\alpha }\int_{0}^{\infty }\frac{\left(
\frac{x}{\alpha }\right) ^{\beta }\log \frac{x}{\alpha }}{1+\left( \frac{x}{%
\alpha }\right) ^{\beta }}\frac{\beta ^{\tau +1}\alpha ^{\left( \tau
+1\right) \beta }x^{\left( \tau +1\right) \left( \beta -1\right) }}{\left(
x^{\beta }+\alpha ^{\beta }\right) ^{2\left( \tau +1\right) }}dx \\
&=&\frac{\beta ^{\tau }}{\alpha ^{\tau +1}}\int_{0}^{\infty }\log t%
\frac{t^{\frac{\tau \left( \beta -1\right) +\beta }{\beta }}}{\left(
t+1\right) ^{2\tau +3}}dt.
\end{eqnarray*}%

By Lemma \ref{lema2} with $m(\beta )=\frac{\tau \left( \beta -1\right)
+\beta }{\beta }$ and $s=2\tau +3,$ we have,%
\begin{equation}
B_{3}=-2\frac{\beta ^{\tau }}{\alpha ^{\tau +1}}B\left( \frac{\tau \beta
+\beta +\tau }{\beta },\frac{\tau \beta +2\beta -\tau }{\beta }\right)
\left\{ \Psi \left( \frac{\tau \beta +2\beta -\tau }{\beta }\right) -\Psi
\left( \frac{\tau \beta +\beta +\tau }{\beta }\right) \right\} .  \label{B3}
\end{equation}%

\underline{In relation to $B_{4}$ we have }

\begin{eqnarray*}
B_{4} &=&-2\alpha ^{\beta -1}\int_{0}^{\infty }\frac{1}{\left(
x^{\beta }+\alpha ^{\beta }\right) }\frac{\beta ^{\tau +1}\alpha ^{\left(
\tau +1\right) \beta }x^{\left( \tau +1\right) \left( \beta -1\right) }}{%
\left( x^{\beta }+\alpha ^{\beta }\right) ^{2\left( \tau +1\right) }}dx \\
&=&-2\frac{\beta ^{\tau }}{\alpha ^{\tau +1}}\int_{0}^{\infty }\frac{%
t^{\frac{\tau \left( \beta -1\right) }{\beta }}}{\left( t+1\right) ^{2\tau
+3}}dt.
\end{eqnarray*}%

By Lemma \ref{lema1} with $m(\beta )=\frac{\tau \left( \beta -1\right)
+\beta }{\beta }$ and $s=2\tau +3,$ we have,%
\begin{equation}
B_{4}=-2\frac{\beta ^{\tau }}{\alpha ^{\tau +1}}B\left( \frac{\tau \beta
+2\beta +\tau }{\beta },\frac{\tau \beta +\beta -\tau }{\beta }\right) .
\label{B4}
\end{equation}%

\underline{In relation to $B_{5}$ we have }%
\begin{eqnarray*}
B_{5} &=&-2\beta \alpha ^{\beta -1}\int_{0}^{\infty }\log \frac{x}{%
\alpha }\frac{1}{\left( x^{\beta }+\alpha ^{\beta }\right) }\frac{\beta
^{\tau +1}\alpha ^{\left( \tau +1\right) \beta }x^{\left( \tau +1\right)
\left( \beta -1\right) }}{\left( x^{\beta }+\alpha ^{\beta }\right)
^{2\left( \tau +1\right) }}dx \\
&=&-2\frac{\beta ^{\tau }}{\alpha ^{\tau +1}}\int_{0}^{\infty }\log t%
\frac{t^{\frac{\tau \left( \beta -1\right) }{\beta }}}{\left( t+1\right)
^{2\tau +3}}dt.
\end{eqnarray*}%

By Lemma \ref{lema2} with $m(\beta )=\frac{\tau \left( \beta -1\right) }{%
\beta }$ and $s=2\tau +3,$ we have,%
\begin{equation}
B_{5}=-2\frac{\beta ^{\tau }}{\alpha ^{\tau +1}}B\left( \frac{\tau \beta
+2\beta +\tau }{\beta },\frac{\tau \beta +\beta -\tau }{\beta }\right)
\left\{ \Psi \left( \frac{\tau \beta +\beta -\tau }{\beta }\right) -\Psi
\left( \frac{\tau \beta +2\beta +\tau }{\beta }\right) \right\} .  \label{B5}
\end{equation}%

\underline{In relation to $B_{6}$ we have }%
\begin{eqnarray*}
B_{6} &=&4\beta \alpha ^{\beta -1}\int_{0}^{\infty }\log \frac{x}{%
\alpha }\frac{\left( \frac{x}{\alpha }\right) ^{\beta }}{\left( x^{\beta
}+\alpha ^{\beta }\right) }\frac{1}{1+\left( \frac{x}{\alpha }\right)
^{\beta }}\frac{\beta ^{\tau +1}\alpha ^{\left( \tau +1\right) \beta
}x^{\left( \tau +1\right) \left( \beta -1\right) }}{\left( x^{\beta }+\alpha
^{\beta }\right) ^{2\left( \tau +1\right) }}dx \\
&=&4\frac{\beta ^{\tau }}{\alpha ^{\tau +1}}\int_{0}^{\infty }\log t%
\frac{t^{\frac{\tau \left( \beta -1\right) +\beta }{\beta }}}{\left(
t+1\right) ^{2\tau +4}}dt.
\end{eqnarray*}%

By Lemma \ref{lema2} with $m(\beta )=\frac{\tau \left( \beta -1\right)
+\beta }{\beta }$ and $s=2\tau +4,$ we have,%
\begin{equation}
B_{6}=4\frac{\beta ^{\tau }}{\alpha ^{\tau +1}}B\left( \frac{\tau \beta
+2\beta +\tau }{\beta },\frac{\tau \beta +2\beta -\tau }{\beta }\right)
\left\{ \Psi \left( \frac{\tau \beta +2\beta -\tau }{\beta }\right) -\Psi
\left( \frac{\tau \beta +2\beta +\tau }{\beta }\right) \right\} .  \label{B6}
\end{equation}

\newpage

\section{Supplementary material: Complementary simulation results}

In the simulation section, we provided the asymptotic results for different values of $n$ and $\beta =2.5$. Here, we present the results for other combinations of $n$ and $\beta $ in order to show that the conclusions are the same when varying $\beta .$

\begin{table}[h]
\begin{center}
\begin{tabular}{|c|cccc|}
\hline
& Bias & RMSE & $\hat{\alpha }$ & $\hat{\beta }$ \\ \hline
MLE & 0.70552 & 0.70018 & 1.06539 & 1.72721 \\
$DPD_{0.1}$ & 0.71532 & 0.71051 & 1.05879 & 1.73198 \\
$DPD_{0.2}$ & 0.75039 & 0.75728 & 1.05021 & 1.76130 \\
$DPD_{0.3}$ & 0.79667 & 0.82319 & 1.04092 & 1.80127 \\
$DPD_{0.4}$ & 0.84772 & 0.89570 & 1.03146 & 1.84625 \\
$DPD_{0.5}$ & 0.89645 & 0.96194 & 1.02225 & 1.88999 \\
$DPD_{0.6}$ & 0.93881 & 1.01797 & 1.01397 & 1.92845 \\
$DPD_{0.7}$ & 0.97577 & 1.06488 & 1.00673 & 1.96250 \\
$DPD_{0.8}$ & 1.00405 & 1.09993 & 1.00026 & 1.98885 \\
$DPD_{0.9}$ & 1.02790 & 1.12795 & 0.99404 & 2.01224 \\
$DPD_{1.0}$ & 1.04767 & 1.14980 & 0.98857 & 2.03121 \\
RM          & 0.73494 & 0.74710 & 1.03998 & 1.60020 \\
SM          & 1.04031 & 0.93286 & 1.24624 & 1.00385 \\
HL          & 0.92904 & 0.78697 & 1.06757 & 0.90906 \\ \hline
\end{tabular}%
\end{center}
\par
\label{tablasimsincont10-1}
\caption{Results for $n=10$ and $\protect\beta = 1.5.$}
\end{table}

\begin{table}
\begin{center}
\begin{tabular}{|c|cccc|}
\hline
& Bias & RMSE & $\hat{\alpha }$ & $\hat{\beta }$ \\ \hline
MLE & 0.86117 & 0.99160 & 1.02635 & 2.88258 \\
$DPD_{0.1}$ & 0.87040 & 1.00759 & 1.02438 & 2.88073 \\
$DPD_{0.2}$ & 0.91207 & 1.07743 & 1.02174 & 2.91365 \\
$DPD_{0.3}$ & 0.97792 & 1.19125 & 1.01872 & 2.97222 \\
$DPD_{0.4}$ & 1.06071 & 1.33208 & 1.01545 & 3.04856 \\
$DPD_{0.5}$ & 1.14342 & 1.46068 & 1.01236 & 3.12639 \\
$DPD_{0.6}$ & 1.22792 & 1.58548 & 1.00932 & 3.20736 \\
$DPD_{0.7}$ & 1.30270 & 1.68709 & 1.00638 & 3.27959 \\
$DPD_{0.8}$ & 1.36963 & 1.77269 & 1.00368 & 3.34440 \\
$DPD_{0.9}$ & 1.42582 & 1.84076 & 1.00110 & 3.40043 \\
$DPD_{1.0}$ & 1.47734 & 1.90194 & 0.99867 & 3.45194 \\
RM          & 0.90345 & 1.09261 & 1.01045 & 2.66286 \\
SM          & 1.27049 & 1.21828 & 1.12259 & 1.67069 \\
HL          & 1.21681 & 1.13798 & 1.02703 & 1.51554 \\ \hline
\end{tabular}%
\end{center}
\par
\label{tablasimsincont10-2}
\caption{Results for $n=10$ and $\protect\beta = 2.5.$}
\end{table}

\begin{table}
\begin{center}
\begin{tabular}{|c|cccc|}
\hline
& Bias & RMSE & $\hat{\alpha }$ & $\hat{\beta }$ \\ \hline
MLE & 1.43374 & 1.89623 & 1.00663 & 5.75789 \\
$DPD_{0.1}$ & 1.45352 & 1.93043 & 1.00619 & 5.75820 \\
$DPD_{0.2}$ & 1.53770 & 2.07945 & 1.00559 & 5.82961 \\
$DPD_{0.3}$ & 1.67607 & 2.33313 & 1.00490 & 5.95832 \\
$DPD_{0.4}$ & 1.85094 & 2.63925 & 1.00410 & 6.12844 \\
$DPD_{0.5}$ & 2.03453 & 2.93733 & 1.00331 & 6.30980 \\
$DPD_{0.6}$ & 2.20377 & 3.18600 & 1.00257 & 6.47981 \\
$DPD_{0.7}$ & 2.35084 & 3.39102 & 1.00182 & 6.62937 \\
$DPD_{0.8}$ & 2.48455 & 3.56323 & 1.00113 & 6.76695 \\
$DPD_{0.9}$ & 2.60501 & 3.70961 & 1.00042 & 6.89455 \\
$DPD_{1.0}$ & 2.71562 & 3.84173 & 0.99975 & 7.01138 \\
RM          & 1.52836 & 2.12342 & 0.99862 & 5.34594 \\
SM          & 2.15292 & 2.31282 & 1.04958 & 3.36417 \\
HL & 2.14313 & 2.21691 & 1.00699 & 3.03342 \\ \hline
\end{tabular}%
\end{center}
\par
\label{tablasimsincont10-3}
\caption{Results for $n=10$ and $\protect\beta = 5.$}
\end{table}

\begin{table}
\begin{center}
\begin{tabular}{|c|cccc|}
\hline
& Bias & RMSE & $\hat{\alpha }$ & $\hat{\beta }$ \\ \hline
MLE & 2.75441 & 3.83907 & 1.00244 & 11.55387 \\
$DPD_{0.1}$ & 2.78746 & 3.89702 & 1.00224 & 11.54996 \\
$DPD_{0.2}$ & 2.95156 & 4.20011 & 1.00199 & 11.69283 \\
$DPD_{0.3}$ & 3.21926 & 4.69697 & 1.00171 & 11.94575 \\
$DPD_{0.4}$ & 3.55924 & 5.29268 & 1.00143 & 12.27509 \\
$DPD_{0.5}$ & 3.90653 & 5.85557 & 1.00115 & 12.62050 \\
$DPD_{0.6}$ & 4.23545 & 6.34660 & 1.00085 & 12.95335 \\
$DPD_{0.7}$ & 4.55107 & 6.78338 & 1.00061 & 13.27909 \\
$DPD_{0.18}$ & 4.81986 & 7.14146 & 1.00037 & 13.56289 \\
$DPD_{0.9}$ & 5.05875 & 7.43550 & 1.00015 & 13.81764 \\
$DPD_{1.0}$ & 5.26077 & 7.66508 & 0.99993 & 14.03665 \\
RM & 2.93315 & 4.29091 & 0.99830 & 10.69892 \\
SM & 4.16657 & 4.64748 & 1.02252 & 6.70247 \\
HL & 4.14191 & 4.42085 & 1.00259 & 6.09889 \\ \hline
\end{tabular}%
\end{center}
\par
\label{tablasimsincont10-4}
\caption{Results for $n=10$ and $\protect\beta = 10.$}
\end{table}

\begin{table}
\begin{center}
\begin{tabular}{|c|cccc|}
\hline
& Bias & RMSE & $\hat{\alpha }$ & $\hat{\beta }$ \\ \hline
MLE & 0.40790 & 0.37874 & 1.02502 & 1.58135 \\
$DPD_{0.1}$ & 0.41409 & 0.38430 & 1.02317 & 1.58186 \\
$DPD_{0.2}$ & 0.42995 & 0.40036 & 1.02065 & 1.58936 \\
$DPD_{0.3}$ & 0.44886 & 0.42109 & 1.01779 & 1.59966 \\
$DPD_{0.4}$ & 0.46808 & 0.44388 & 1.01478 & 1.61104 \\
$DPD_{0.5}$ & 0.48622 & 0.46717 & 1.01178 & 1.62253 \\
$DPD_{0.6}$ & 0.50263 & 0.48936 & 1.00886 & 1.63353 \\
$DPD_{0.7}$ & 0.51732 & 0.50943 & 1.00612 & 1.64370 \\
$DPD_{0.8}$ & 0.53043 & 0.52651 & 1.00359 & 1.65302 \\
$DPD_{0.9}$ & 0.54170 & 0.54110 & 1.00129 & 1.66111 \\
$DPD_{1.0}$ & 0.55217 & 0.55504 & 0.99907 & 1.66927 \\
RM & 0.43073 & 0.40264 & 1.02557 & 1.54105 \\
SM & 0.73254 & 0.62951 & 1.03458 & 1.00712 \\
HL & 0.80669 & 0.68330 & 1.02878 & 0.88301 \\ \hline
\end{tabular}%
\end{center}
\par
\label{tablasimsincont25-1}
\caption{Results for $n=25$ and $\protect\beta = 1.5.$}
\end{table}

\begin{table}
\begin{center}
\begin{tabular}{|c|cccc|}
\hline
& Bias & RMSE & $\hat{\alpha }$ & $\hat{\beta }$ \\ \hline
MLE & 0.48028 & 0.50752 & 1.00983 & 2.63727 \\
$DPD_{0.1}$ & 0.48763 & 0.51759 & 1.00919 & 2.63747 \\
$DPD_{0.2}$ & 0.50554 & 0.54230 & 1.00828 & 2.64756 \\
$DPD_{0.3}$ & 0.52864 & 0.57471 & 1.00722 & 2.66305 \\
$DPD_{0.4}$ & 0.55465 & 0.61302 & 1.00605 & 2.68205 \\
$DPD_{0.5}$ & 0.58227 & 0.65862 & 1.00485 & 2.70358 \\
$DPD_{0.6}$ & 0.60767 & 0.69806 & 1.00369 & 2.72405 \\
$DPD_{0.7}$ & 0.63134 & 0.73463 & 1.00257 & 2.74390 \\
$DPD_{0.8}$ & 0.65404 & 0.77150 & 1.00151 & 2.76356 \\
$DPD_{0.9}$ & 0.67422 & 0.80133 & 1.00055 & 2.78127 \\
$DPD_{1.0}$ & 0.69244 & 0.82689 & 0.99964 & 2.79784 \\
RM & 0.51604 & 0.55662 & 1.01006 & 2.57113 \\
SM & 0.99659 & 0.95598 & 1.01341 & 1.67301 \\
HL & 1.14041 & 1.07281 & 1.01168 &  1.47365 \\ \hline
\end{tabular}%
\end{center}
\par
\label{tablasimsincont25-2}
\caption{Results for $n=25$ and $\protect\beta = 2.5.$}
\end{table}

\begin{table}
\begin{center}
\begin{tabular}{|c|cccc|}
\hline
& Bias & RMSE & $\hat{\alpha }$ & $\hat{\beta }$ \\ \hline
MLE & 0.78984 & 0.96720 & 1.00458 & 5.27681 \\
$DPD_{0.1}$ & 0.79833 & 0.98166 & 1.00426 & 5.27431 \\
$DPD_{0.2}$ & 0.82748 & 1.02653 & 1.00389 & 5.29165 \\
$DPD_{0.3}$ & 0.86827 & 1.08952 & 1.00350 & 5.32083 \\
$DPD_{0.4}$ & 0.91568 & 1.16530 & 1.00310 & 5.35779 \\
$DPD_{0.5}$ & 0.96642 & 1.25211 & 1.00269 & 5.40009 \\
$DPD_{0.6}$ & 1.01671 & 1.33428 & 1.00230 & 5.44362 \\
$DPD_{0.7}$ & 1.06463 & 1.41110 & 1.00192 & 5.48675 \\
$DPD_{0.8}$ & 1.11075 & 1.48487 & 1.00158 & 5.52953 \\
$DPD_{0.9}$ & 1.15258 & 1.55327 & 1.00124 & 5.56955 \\
$DPD_{1.0}$ & 1.18951 & 1.61166 & 1.00093 & 5.60614 \\
RM & 0.85734 & 1.06339 & 1.00456 & 5.13022 \\
SM & 1.80725 & 1.89398 & 1.00532 & 3.33789 \\
HL & 2.11688 & 2.12979 & 1.00568 & 2.94310 \\ \hline
\end{tabular}%
\end{center}
\par
\label{tablasimsincont25-3}
\caption{Results for $n=25$ and $\protect\beta = 5.0.$}
\end{table}

\begin{table}
\begin{center}
\begin{tabular}{|c|cccc|}
\hline
& Bias & RMSE & $\hat{\alpha }$ & $\hat{\beta }$ \\ \hline
MLE & 1.50633 & 1.93866 & 1.00132 & 10.56004 \\
$DPD_{0.1}$ & 1.52754 & 1.96981 & 1.00125 & 10.55993 \\
$DPD_{0.2}$ & 1.58713 & 2.06369 & 1.00117 & 10.60041 \\
$DPD_{0.3}$ & 1.66935 & 2.19889 & 1.00108 & 10.66605 \\
$DPD_{0.4}$ & 1.76512 & 2.36284 & 1.00098 & 10.74774 \\
$DPD_{0.5}$ & 1.86848 & 2.53665 & 1.00089 & 10.83841 \\
$DPD_{0.6}$ & 1.97319 & 2.72036 & 1.00080 & 10.93398 \\
$DPD_{0.7}$ & 2.07177 & 2.88695 & 1.00071 & 11.02681 \\
$DPD_{0.8}$ & 2.16179 & 3.03810 & 1.00063 & 11.11428 \\
$DPD_{0.9}$ & 2.24256 & 3.17238 & 1.00055 & 11.19575 \\
$DPD_{1.0}$ & 2.31232 & 3.28251 & 1.00048 & 11.26843 \\
RM & 1.64705 & 2.14634 & 1.00140 & 10.28970 \\
SM & 3.49014 & 3.77015 & 1.00158 & 6.68764 \\
HL & 4.13418 & 4.24730 & 1.00187 & 5.89885 \\ \hline
\end{tabular}%
\end{center}
\par
\label{tablasimsincont25-4}
\caption{Results for $n=25$ and $\protect\beta = 10.$}
\end{table}

\begin{table}
\begin{center}
\begin{tabular}{|c|cccc|}
\hline
& Bias & RMSE & $\hat{\alpha }$ & $\hat{\beta }$ \\ \hline
MLE & 0.28085 & 0.25512 & 1.01466 & 1.53923 \\
$DPD_{0.1}$ & 0.28551 & 0.25997 & 1.01354 & 1.53968 \\
$DPD_{0.2}$ & 0.29539 & 0.26986 & 1.01212 & 1.54306 \\
$DPD_{0.3}$ & 0.30639 & 0.28085 & 1.01061 & 1.54724 \\
$DPD_{0.4}$ & 0.31668 & 0.29137 & 1.00908 & 1.55149 \\
$DPD_{0.5}$ & 0.32585 & 0.30090 & 1.00760 & 1.55551 \\
$DPD_{0.6}$ & 0.33393 & 0.30934 & 1.00620 & 1.55921 \\
$DPD_{0.7}$ & 0.34096 & 0.31676 & 1.00490 & 1.56258 \\
$DPD_{0.8}$ & 0.34722 & 0.32338 & 1.00372 & 1.56567 \\
$DPD_{0.9}$ & 0.35297 & 0.32948 & 1.00261 & 1.56864 \\
$DPD_{1.0}$ & 0.35840 & 0.33532 & 1.00160 & 1.57151 \\
RM & 0.29520 & 0.26826 & 1.01176 & 1.50809 \\
SM & 0.72610 & 0.62477 & 1.04661 & 0.93438 \\
HL & 0.75229 & 0.65215 & 1.01450 & 0.88029 \\ \hline
\end{tabular}%
\end{center}
\par
\label{tablasimsincont50-1}
\caption{Results for $n=50$ and $\protect\beta = 1.5.$}
\end{table}

\begin{table}
\begin{center}
\begin{tabular}{|c|cccc|}
\hline
& Bias & RMSE & $\hat{\alpha }$ & $\hat{\beta }$ \\ \hline
MLE & 0.32909 & 0.33628 & 1.00623 & 2.56531 \\
$DPD_{0.1}$ & 0.33391 & 0.34216 & 1.00580 & 2.56618 \\
$DPD_{0.2}$ & 0.34529 & 0.35621 & 1.00524 & 2.57183 \\
$DPD_{0.3}$ & 0.35938 & 0.37378 & 1.00463 & 2.57971 \\
$DPD_{0.4}$ & 0.37397 & 0.39244 & 1.00399 & 2.58854 \\
$DPD_{0.5}$ & 0.38798 & 0.41086 & 1.00336 & 2.59760 \\
$DPD_{0.6}$ & 0.40103 & 0.42831 & 1.00275 & 2.60644 \\
$DPD_{0.7}$ & 0.41301 & 0.44449 & 1.00217 & 2.61487 \\
$DPD_{0.8}$ & 0.42391 & 0.45940 & 1.00164 & 2.62283 \\
$DPD_{0.9}$ & 0.43381 & 0.47280 & 1.00114 & 2.63015 \\
$DPD_{1.0}$ & 0.44288 & 0.48480 & 1.00068 & 2.63692 \\
RM & 0.35286 & 0.36535 & 1.00431 & 2.51693 \\
SM & 1.03277 & 0.98680 & 1.02396 & 1.56206 \\
HL & 1.10907 & 1.05263 & 1.00620 & 1.46905 \\ \hline
\end{tabular}%
\end{center}
\par
\label{tablasimsincont50-2}
\caption{Results for $n=50$ and $\protect\beta = 2.5.$}
\end{table}

\begin{table}
\begin{center}
\begin{tabular}{|c|cccc|}
\hline
& Bias & RMSE & $\hat{\alpha }$ & $\hat{\beta }$ \\ \hline
MLE & 0.53098 & 0.63556 & 1.00084 & 5.13463 \\
$DPD_{0.1}$ & 0.53836 & 0.64555 & 1.00068 & 5.13623 \\
$DPD_{0.2}$ & 0.55749 & 0.67100 & 1.00050 & 5.14656 \\
$DPD_{0.3}$ & 0.58207 & 0.70445 & 1.00030 & 5.16158 \\
$DPD_{0.4}$ & 0.60873 & 0.74133 & 1.00011 & 5.17910 \\
$DPD_{0.5}$ & 0.63588 & 0.77884 & 0.99992 & 5.19774 \\
$DPD_{0.6}$ & 0.66201 & 0.81520 & 0.99973 & 5.21653 \\
$DPD_{0.7}$ & 0.68635 & 0.84935 & 0.99956 & 5.23485 \\
$DPD_{0.8}$ & 0.70862 & 0.88063 & 0.99940 & 5.25229 \\
$DPD_{0.9}$ & 0.72872 & 0.90867 & 0.99926 & 5.26866 \\
$DPD_{1.0}$ & 0.74682 & 0.93360 & 0.99913 & 5.28395 \\
RM & 0.58109 & 0.69556 & 1.00003 & 5.04133 \\
SM & 1.92898 & 1.96102 & 1.00908 & 3.12027 \\
HL & 2.09928 & 2.09405 & 1.00089 & 2.93963 \\ \hline
\end{tabular}%
\end{center}
\par
\label{tablasimsincont50-3}
\caption{Results for $n=50$ and $\protect\beta = 5.0.$}
\end{table}

\begin{table}
\begin{center}
\begin{tabular}{|c|cccc|}
\hline
& Bias & RMSE & $\hat{\alpha }$ & $\hat{\beta }$ \\ \hline
MLE & 1.00061 & 1.26915 & 1.00048 & 10.25602 \\
$DPD_{0.1}$ & 1.01559 & 1.28824 & 1.00044 & 10.25445 \\
$DPD_{0.2}$ & 1.05135 & 1.33645 & 1.00041 & 10.27116 \\
$DPD_{0.3}$ & 1.09804 & 1.40064 & 1.00037 & 10.29863 \\
$DPD_{0.4}$ & 1.14868 & 1.47253 & 1.00034 & 10.33251 \\
$DPD_{0.5}$ & 1.20054 & 1.54677 & 1.00031 & 10.36990 \\
$DPD_{0.6}$ & 1.25113 & 1.61986 & 1.00028 & 10.40870 \\
$DPD_{0.7}$ & 1.29852 & 1.68959 & 1.00026 & 10.44744 \\
$DPD_{0.8}$ & 1.34226 & 1.75446 & 1.00024 & 10.48508 \\
$DPD_{0.9}$ & 1.38170 & 1.81327 & 1.00023 & 10.52091 \\
$DPD_{1.0}$ & 1.41712 & 1.86547 & 1.00022 & 10.55462 \\
RM & 1.11054 & 1.39728 & 1.00001 & 10.07130 \\
SM & 3.78958 & 3.92204 & 1.00466 & 6.24218 \\
HL & 4.15225 & 4.19937 & 1.00049 & 5.86731 \\ \hline
\end{tabular}%
\end{center}
\par
\label{tablasimsincont50-4}
\caption{Results for $n=50$ and $\protect\beta = 10.$}
\end{table}

\begin{table}
\begin{center}
\begin{tabular}{|c|cccc|}
\hline
& Bias & RMSE & $\hat{\alpha }$ & $\hat{\beta }$ \\ \hline
MLE & 0.22559 & 0.20261 & 1.00858 & 1.52665 \\
$DPD_{0.1}$ & 0.23003 & 0.20670 & 1.00781 & 1.52749 \\
$DPD_{0.2}$ & 0.23863 & 0.21468 & 1.00689 & 1.53003 \\
$DPD_{0.3}$ & 0.24772 & 0.22345 & 1.00592 & 1.53297 \\
$DPD_{0.4}$ & 0.25623 & 0.23179 & 1.00494 & 1.53589 \\
$DPD_{0.5}$ & 0.26365 & 0.23927 & 1.00400 & 1.53860 \\
$DPD_{0.6}$ & 0.27005 & 0.24584 & 1.00312 & 1.54107 \\
$DPD_{0.7}$ & 0.27562 & 0.25160 & 1.00230 & 1.54330 \\
$DPD_{0.8}$ & 0.28054 & 0.25668 & 1.00155 & 1.54535 \\
$DPD_{0.9}$ & 0.28496 & 0.26126 & 1.00086 & 1.54726 \\
$DPD_{1.0}$ & 0.28902 & 0.26547 & 1.00022 & 1.54908 \\
RM & 0.23900 & 0.21600 & 1.00890 & 1.50996 \\
SM & 0.67792 & 0.59172 & 1.01045 & 0.94681 \\
HL & 0.72630 & 0.64081 & 1.00874 & 0.87998 \\ \hline
\end{tabular}%
\end{center}
\par
\label{tablasimsincont75-1}
\caption{Results for $n=75$ and $\protect\beta = 1.5.$}
\end{table}

\begin{table}
\begin{center}
\begin{tabular}{|c|cccc|}
\hline
& Bias & RMSE & $\hat{\alpha }$ & $\hat{\beta }$ \\ \hline
MLE & 0.26557 & 0.26803 & 1.00389 & 2.54115 \\
$DPD_{0.1}$ & 0.27000 & 0.27299 & 1.00360 & 2.54085 \\
$DPD_{0.2}$ & 0.27886 & 0.28332 & 1.00325 & 2.54332 \\
$DPD_{0.3}$ & 0.28899 & 0.29529 & 1.00286 & 2.54706 \\
$DPD_{0.4}$ & 0.29927 & 0.30739 & 1.00245 & 2.55142 \\
$DPD_{0.5}$ & 0.30906 & 0.31896 & 1.00204 & 2.55603 \\
$DPD_{0.6}$ & 0.31804 & 0.32973 & 1.00162 & 2.56064 \\
$DPD_{0.7}$ & 0.32625 & 0.33960 & 1.00122 & 2.56513 \\
$DPD_{0.8}$ & 0.33368 & 0.34861 & 1.00084 & 2.56944 \\
$DPD_{0.9}$ & 0.34047 & 0.35687 & 1.00048 & 2.57356 \\
$DPD_{1.0}$ & 0.34668 & 0.36449 & 1.00014 & 2.57749 \\
RM & 0.28623 & 0.29136 & 1.00393 & 2.51036 \\
SM & 1.00466 & 0.96143 & 1.00495 & 1.57055 \\
HL & 1.10116 & 1.05026 & 1.00395 & 1.46391 \\ \hline
\end{tabular}%
\end{center}
\par
\label{tablasimsincont75-2}
\caption{Results for $n=75$ and $\protect\beta = 2.5.$}
\end{table}

\begin{table}
\begin{center}
\begin{tabular}{|c|cccc|}
\hline
& Bias & RMSE & $\hat{\alpha }$ & $\hat{\beta }$ \\ \hline
MLE & 0.42703 & 0.50897 & 1.00090 & 5.08768 \\
$DPD_{0.1}$ & 0.43322 & 0.51721 & 1.00080 & 5.08817 \\
$DPD_{0.2}$ & 0.44721 & 0.53591 & 1.00069 & 5.09441 \\
$DPD_{0.3}$ & 0.46505 & 0.55990 & 1.00056 & 5.10381 \\
$DPD_{0.4}$ & 0.48415 & 0.58620 & 1.00042 & 5.11491 \\
$DPD_{0.5}$ & 0.50338 & 0.61282 & 1.00029 & 5.12673 \\
$DPD_{0.6}$ & 0.52169 & 0.63846 & 1.00016 & 5.13862 \\
$DPD_{0.7}$ & 0.53871 & 0.66234 & 1.00004 & 5.15015 \\
$DPD_{0.8}$ & 0.55421 & 0.68404 & 0.99993 & 5.16109 \\
$DPD_{0.9}$ & 0.56818 & 0.70342 & 0.99983 & 5.17132 \\
$DPD_{1.0}$ & 0.58074 & 0.72055 & 0.99974 & 5.18084 \\
RM & 0.46836 & 0.56012 & 1.00091 & 5.03146 \\
SM & 1.88233 & 1.90110 & 1.00106 & 3.15576 \\
HL & 2.09985 & 2.08990 & 1.00089 & 2.93212 \\ \hline
\end{tabular}%
\end{center}
\par
\label{tablasimsincont75-3}
\caption{Results for $n=75$ and $\protect\beta = 5.$}
\end{table}

\begin{table}
\begin{center}
\begin{tabular}{|c|cccc|}
\hline
& Bias & RMSE & $\hat{\alpha }$ & $\hat{\beta }$ \\ \hline
MLE & 0.81532 & 1.02538 & 1.00026 & 10.17362 \\
$DPD_{0.1}$ & 0.82764 & 1.04127 & 1.00026 & 10.17328 \\
$DPD_{0.2}$ & 0.85563 & 1.07762 & 1.00025 & 10.18411 \\
$DPD_{0.3}$ & 0.89010 & 1.12444 & 1.00024 & 10.20128 \\
$DPD_{0.4}$ & 0.92709 & 1.17575 & 1.00023 & 10.22205 \\
$DPD_{0.5}$ & 0.96456 & 1.22778 & 1.00022 & 10.24462 \\
$DPD_{0.6}$ & 1.00013 & 1.27808 & 1.00022 & 10.26774 \\
$DPD_{0.7}$ & 1.03306 & 1.32515 & 1.00021 & 10.29052 \\
$DPD_{0.8}$ & 1.06322 & 1.36812 & 1.00021 & 10.31242 \\
$DPD_{0.9}$ & 1.09036 & 1.40662 & 1.00021 & 10.33314 \\
$DPD_{1.0}$ & 1.11443 & 1.44069 & 1.00021 & 10.35259 \\
RM & 0.90410 & 1.13202 & 1.00026 & 10.05669 \\
SM & 3.71061 & 3.80447 & 1.00058 & 6.30922 \\
HL & 4.15315 & 4.18142 & 1.00026 & 5.86290 \\ \hline
\end{tabular}%
\end{center}
\par
\label{tablasimsincont75-4}
\caption{Results for $n=75$ and $\protect\beta = 10.$}
\end{table}

\begin{table}
\begin{center}
\begin{tabular}{|c|cccc|}
\hline
& Bias & RMSE & $\hat{\alpha }$ & $\hat{\beta }$ \\ \hline
MLE & 0.19593 & 0.17578 & 1.00715 & 1.52121 \\
$DPD_{0.1}$ & 0.19897 & 0.17874 & 1.00654 & 1.52169 \\
$DPD_{0.2}$ & 0.20573 & 0.18501 & 1.00582 & 1.52344 \\
$DPD_{0.3}$ & 0.21313 & 0.19192 & 1.00507 & 1.52549 \\
$DPD_{0.4}$ & 0.21998 & 0.19846 & 1.00433 & 1.52752 \\
$DPD_{0.5}$ & 0.22597 & 0.20429 & 1.00362 & 1.52942 \\
$DPD_{0.6}$ & 0.23112 & 0.20938 & 1.00294 & 1.53115 \\
$DPD_{0.7}$ & 0.23559 & 0.21383 & 1.00231 & 1.53273 \\
$DPD_{0.8}$ & 0.23950 & 0.21776 & 1.00172 & 1.53418 \\
$DPD_{0.9}$ & 0.24303 & 0.22131 & 1.00118 & 1.53555 \\
$DPD_{1.0}$ & 0.24626 & 0.22458 & 1.00067 & 1.53685 \\
RM & 0.20642 & 0.18569 & 1.00603 & 1.50542 \\
SM & 0.68084 & 0.59860 & 1.02299 & 0.92929 \\
HL & 0.71402 & 0.63639 & 1.00736 & 0.87946 \\ \hline
\end{tabular}%
\end{center}
\par
\label{tablasimsincont100-1}
\caption{Results for $n=100$ and $\protect\beta = 1.5.$}
\end{table}

\begin{table}
\begin{center}
\begin{tabular}{|c|cccc|}
\hline
& Bias & RMSE & $\hat{\alpha }$ & $\hat{\beta }$ \\ \hline
MLE & 0.22742 & 0.22921 & 1.00262 & 2.53247 \\
$DPD_{0.1}$ & 0.23106 & 0.23338 & 1.00241 & 2.53353 \\
$DPD_{0.2}$ & 0.23835 & 0.24227 & 1.00215 & 2.53658 \\
$DPD_{0.3}$ & 0.24682 & 0.25275 & 1.00186 & 2.54046 \\
$DPD_{0.4}$ & 0.25548 & 0.26342 & 1.00157 & 2.54465 \\
$DPD_{0.5}$ & 0.26389 & 0.27363 & 1.00129 & 2.54888 \\
$DPD_{0.6}$ & 0.27175 & 0.28309 & 1.00102 & 2.55299 \\
$DPD_{0.7}$ & 0.27900 & 0.29171 & 1.00077 & 2.55689 \\
$DPD_{0.8}$ & 0.28556 & 0.29950 & 1.00053 & 2.56055 \\
$DPD_{0.9}$ & 0.29153 & 0.30655 & 1.00032 & 2.56397 \\
$DPD_{1.0}$ & 0.29700 & 0.31298 & 1.00012 & 2.56718 \\
RM & 0.24559 & 0.25133 & 1.00259 & 2.50775 \\
SM & 1.01447 & 0.97193 & 1.01121 & 1.55085 \\
HL & 1.09114 & 1.04550 & 1.00277 & 1.46494 \\ \hline
\end{tabular}%
\end{center}
\par
\label{tablasimsincont100-2}
\caption{Results for $n=100$ and $\protect\beta = 2.5.$}
\end{table}

\begin{table}
\begin{center}
\begin{tabular}{|c|cccc|}
\hline
& Bias & RMSE & $\hat{\alpha }$ & $\hat{\beta }$ \\ \hline
MLE & 0.36844 & 0.43431 & 1.00030 & 5.06549 \\
$DPD_{0.1}$ & 0.37232 & 0.44030 & 1.00024 & 5.06584 \\
$DPD_{0.2}$ & 0.38300 & 0.45529 & 1.00016 & 5.07067 \\
$DPD_{0.3}$ & 0.39717 & 0.47472 & 1.00007 & 5.07793 \\
$DPD_{0.4}$ & 0.41270 & 0.49589 & 0.99999 & 5.08649 \\
$DPD_{0.5}$ & 0.42850 & 0.51711 & 0.99990 & 5.09562 \\
$DPD_{0.6}$ & 0.44379 & 0.53734 & 0.99982 & 5.10485 \\
$DPD_{0.7}$ & 0.45786 & 0.55603 & 0.99975 & 5.11385 \\
$DPD_{0.8}$ & 0.47062 & 0.57299 & 0.99968 & 5.12247 \\
$DPD_{0.9}$ & 0.48213 & 0.58821 & 0.99962 & 5.13062 \\
$DPD_{1.0}$ & 0.49251 & 0.60182 & 0.99957 & 5.13828 \\
RM & 0.40334 & 0.47948 & 0.99994 & 5.01636 \\
SM & 1.93207 & 1.93817 & 1.00447 & 3.10069 \\
HL & 2.09938 & 2.08758 & 1.00038 & 2.92855 \\ \hline
\end{tabular}%
\end{center}
\par
\label{tablasimsincont100-3}
\caption{Results for $n=100$ and $\protect\beta = 5.$}
\end{table}

\begin{table}
\begin{center}
\begin{tabular}{|c|cccc|}
\hline
& Bias & RMSE & $\hat{\alpha }$ & $\hat{\beta }$ \\ \hline
MLE & 0.69407 & 0.86490 & 1.00005 & 10.13583 \\
$DPD_{0.1}$ & 0.70174 & 0.87637 & 1.00001 & 10.13524 \\
$DPD_{0.2}$ & 0.72326 & 0.90545 & 0.99996 & 10.14356 \\
$DPD_{0.3}$ & 0.75158 & 0.94370 & 0.99992 & 10.15692 \\
$DPD_{0.4}$ & 0.78316 & 0.98577 & 0.99987 & 10.17306 \\
$DPD_{0.5}$ & 0.81507 & 1.02828 & 0.99982 & 10.19055 \\
$DPD_{0.6}$ & 0.84533 & 1.06918 & 0.99977 & 10.20839 \\
$DPD_{0.7}$ & 0.87336 & 1.10731 & 0.99973 & 10.22596 \\
$DPD_{0.8}$ & 0.89879 & 1.14210 & 0.99969 & 10.24285 \\
$DPD_{0.9}$ & 0.92162 & 1.17335 & 0.99965 & 10.25886 \\
$DPD_{1.0}$ & 0.94200 & 1.20116 & 0.99961 & 10.27390 \\
RM & 0.76712 & 0.95667 & 0.99990 & 10.03177 \\
SM & 3.81978 & 3.88194 & 1.00186 & 6.19626 \\
HL & 4.15390 & 4.17153 & 1.00008 & 5.86007 \\ \hline
\end{tabular}%
\end{center}
\par
\label{tablasimsincont100-4}
\caption{Results for $n=100$ and $\protect\beta = 10.$}
\end{table}


\begin{thebibliography}{99}
\bibitem{abta16} Abbas, K. and Tang, Y. (2016). Objective Bayesian analysis for
log-logistic distribution. \emph{Communications in Statistics - Simulation
and Computation,} 45, 2782--2791.

\bibitem{acsi86} Acreman, A.C. and Sinclair, C.D. (1986). Classification of drainage basins according to their physical characteristics; an application for flood frequency in Scotland. \emph{Journal of Hydrology}, 84(3), 365--380.
%

\bibitem{alsi66} Ali, S. M. and Silvey, S. D. (1966). A general class of
coefficient of divergence of one distribution from another. \emph{Journal
of the Royal Statistical Society}, 28(1), 131--142.

\bibitem{asma03} Ashkar, F., and Mahdi, S. (2003). Comparison of two testing
methods for the log-logistic distribution. \emph{Water Resources Research},
39, 12-17.

\bibitem{bama87} Balakrishnan, N. and Malik, H. J. (1987). Moments of order
statistics from truncated loglogistic distribution. \emph{Journal of
Statistical Planning and Inference,} 17, 251--267.

\bibitem{bachghpa22} Basu, A., Chakraborty, S., Ghosh, A. and Pardo, L. (2022).
Robust density power divergence based tests in multivariate analysis: a
comparative overview of different approaches. \emph{Journal of Multivariate Analysis}, 188,
Paper No. 104846.

\bibitem{baghmamapa17} Basu, A., Ghosh, A., Mandal, A., Mart\'{\i}n, N. and Pardo,
L. (2017). A Wald-type test statistic for testing linear hypothesis in
logistic regression models based on minimum density power divergence
estimator. \emph{Electronic Journal of Statistics,} 11(2), 2741--2772.

\bibitem{baghmamapa21} Basu, A., Ghosh, A., Mandal, A., Mart\'{\i}n, N. and Pardo,
L.(2021). Robust Wald-type tests in GLM with random design based on minimum
density power divergence estimators. \emph{Statistical Methods and Applications}, 30(3), 973--1005.

\bibitem{baghmapa18} Basu, A., Ghosh, A., Mart\'{\i}n, N. and Pardo, L. (2018). Robust
Wald-type tests for non-homogeneous observations based on the minimum
density power divergence estimator. \emph{Metrika}, 81(5), 493--522.

\bibitem{baghmapa22} Basu, A., Ghosh, A., Mart\'{\i}n, N. and Pardo, L. (2022). A
robust generalization of the Rao test. \emph{Journal of Business and Economic Statistics}, 40(2),
868--879.

\bibitem{bahahjjo98} Basu, A., Harris, I. R., Hjort, N. L. and Jones, M. C.
(1998). Robust and efficient estimation by minimizing a density power
divergence. \emph{Biometrika,} 85, 549--559.

\bibitem{bamamapa15} Basu, A., Mandal, A., Mart\'{\i}n, N. and Pardo, L. (2015). Robust tests
for the equality of two normal means based on the density power divergence. \emph{Metrika}, 78(5), 611--634.

\bibitem{bamamapa18} Basu, A., Mandal, A., Mart\'{\i}n, N. and Pardo, L. (2018). Testing composite hypothesis based on the density power divergence. \emph{Sankhya
B}, 80(2), 222--262.

\bibitem{bamamapa19} Basu, A., Mandal, A., Mart\'{\i}n, N. and Pardo, L. (2019). A
robust Wald-type test for testing the equality of two means from log-normal
samples. \emph{Methodology and Computing in Applied Probability,} 21(1), 85--107.

\bibitem{bamamapa16} Basu, A., Mandal, A., Mart\'{\i}n, N. and Pardo, L. (2016).
Generalized Wald-type tests based on minimum density power divergence
estimators. \emph{Statistics}, 50(1), 1--26.


\bibitem{bur42} Burr, I. W. (1942). Cumulative frequency functions. \emph{%
Annals of Mathematical Statistics}, 13, 215--232.

\bibitem{csi63} Csisz\'{a}r, I. (1963). Eine Informationstheoretische
Ungleichung und ihre Anwendung auf den Bewis der Ergodizit\"{a}t on
Markhoffschen Ketten. \emph{Publications of the Mathematical Institute of
the Hungarian Academy of Sciences}, 8, 84--108.


\bibitem{fis61} Fisk, P.R. (1961). The graduation of income distributions. \emph{%
Econometrica}, 29, 171--185.


\bibitem{ghmamapa16} Ghosh, A., Mandal, A., Mart\'{\i}n, N. and Pardo, L. (2016).
Influence analysis of robust Wald-type tests. \emph{Journal of Multivariate Analysis}, 147,
102--126.

\bibitem{ghmabapa18} Ghosh, A., Mart\'{\i}n, N., Basu, A. and Pardo, L. (2018). A new
class of robust two-sample Wald-type tests. \emph{International Journal of Biostatistics,} 14(2),
20170023.


\bibitem{ghbapa21} Ghosh, A., Basu, A. and Pardo, L. (2021). Robust Wald-type
tests under random censoring. \emph{Statistics in Medicine}, 40(5), 1285--1305.

\bibitem{guaklv99} Gupta, R.C., Akman, O. and Lvin, S. (1999). A study of
log-logistic model in survival analysis. \emph{Biometrical Journal, }41,
431-443.

\bibitem{hechqi20} He, X., Chen, W. and Qian, W. (2020). Maximum likelihood
estimators of the parameters of the log-logistic distribution. \emph{%
Statistical Papers}, 61, 1875--1892.

\bibitem{hechya21} He, X., Chen, W. and Yang, R. (2021) Modified best linear
unbiased estimator of the shape parameter of log-logistic distribution.
\emph{Journal of Statistical Computation and Simulation}, 91, 2, 383-395.

\bibitem{kasr02} Kantam, R.R.L. and Srinavasa, G. (2002). Log-logistic distribution; modified maximum likelihood estimation.
\emph{Gujarat Statistical Review}, 29(1-2(, 25--36.



\bibitem{mawapa23} Ma, Z. , Wang, M. and Park, C. (2023). Robust Explicit
Estimation of the Log-Logistic Distribution with Applications. \emph{Journal
of Statistical Theory and Practice,} 17-21.

\bibitem{mar20} Mart\'{\i}n, N. (2020). Rao's Score Tests on Correlation
Matrices. arXiv:2012.14238.

\bibitem{mar20b} Mart\'{\i}n, N. (2020). Robust and efficient Breusch-Pagan
test-statistic: an application of the beta-score Lagrange multipliers test
for non-identically distributed individuals. arXiv:2301.07245.

\bibitem{mijo73} Mielke, P. W. and Johnson, E. S. (1973). Three parameter Kappa
distribution maximum likelihood estimates and likelihood ratio tests. \emph{%
Monthly Weather Review,} 101, 701--709.


\bibitem{papazo02} Pardo, J. A., Pardo, L. and Zografos, K. (2002). Minimum
Phi-divergence etimators with constraints in multiomial populations. \emph{Journal
of Statistical Planning and Inference}, 104, 221--237.

\bibitem{redowa18} Reath, J., Dong, J. and Wang, M. (2018). Improved parameter
estimation of the log-logistic distribution with applications. \emph{%
Computational Statistics}, 33(1), 339--356.


\bibitem{shmitr88} Shoukri, M. M., Mian, I. U. M. and Tracy, D. (1988). Sampling
properties of estimators of log-logistic distribution with application to
Canadian precipitation data. \emph{Canadian Journal of Statistics}, 16,
223--236.

\bibitem{zhchtswa21} Zheng, X., Chiang, J.-Y., Tsai, T.-R. and Wang, S. (2021). Estimating the failure rate of the log-logistic distribution by smooth adaptive and bias-corrections methods. \emph{Computers $\& $ Industrial Engeneering}, 156, 107188.
\end{thebibliography}
\end{document}